\documentclass[twoside, a4paper]{article}
\usepackage{amsmath}			
\usepackage{amssymb}			
\usepackage{amsthm}			
\usepackage{mathrsfs}			
\usepackage{dsfont}			
\usepackage{geometry}			
\usepackage{enumitem}			
\usepackage{stmaryrd}			
\usepackage{longtable}			
\usepackage{verbatim}			
\usepackage[all]{xy}			
\usepackage{fancyhdr}          	 	
\usepackage[colorlinks]{hyperref}	
\usepackage[UKenglish]{babel}     	

\numberwithin{equation}{section}

\title{Convergence to the boundary for random walks on discrete quantum groups and monoidal categories}
\author{Bas P. A. Jordans}
\date{\today}

\newtheorem{Thm}{Theorem}[section]
\newtheorem{Prop}[Thm]{Proposition}
\newtheorem{Lem}[Thm]{Lemma}
\newtheorem{Cor}[Thm]{Corollary}
\newtheorem{Conj}[Thm]{Conjecture}

\theoremstyle{definition}
\newtheorem{Assum}[Thm]{Assumption}
\newtheorem{Def}[Thm]{Definition}
\newtheorem{Exam}[Thm]{Example}
\newtheorem{Rem}[Thm]{Remark}
\newtheorem{Not}[Thm]{Notation}

\newenvironment{pf}[1][]{\noindent\vspace{0pt plus 1pt minus 1pt}{\it Proof#1. }}{\mbox{}\hfill$\boxtimes$\vspace{6pt plus 1pt minus 2 pt}}
\newenvironment{sketch}[1][]{\noindent\vspace{0pt plus 1pt minus 1pt}{\it Sketch of proof#1. }}{\mbox{}\hfill$\boxtimes$\vspace{6pt plus 1pt minus 2 pt}}

\newlist{enumeraterm}{enumerate}{1}%
\setlist[enumeraterm]{label={\rm{(\roman*)}}}%

\hypersetup{
citecolor = {blue},
}

\newcommand{\Nat}{\mathbb{N}}
\newcommand{\Int}{\mathbb{Z}}

\newcommand{\Rea}{\mathbb{R}}
\newcommand{\Com}{\mathbb{C}}

\newcommand{\Acal}{\mathcal{A}}
\newcommand{\Bcal}{\mathcal{B}}
\newcommand{\Ccal}{\mathcal{C}}
\newcommand{\Dcal}{\mathcal{D}}
\newcommand{\Ecal}{\mathcal{E}}
\newcommand{\Fcal}{\mathcal{F}}

\newcommand{\Hcal}{\mathcal{H}}

\newcommand{\Mcal}{\mathcal{M}}
\newcommand{\Ncal}{\mathcal{N}}

\newcommand{\Pcal}{\mathcal{P}}

\newcommand{\Ucal}{\mathcal{U}}

\newcommand{\sufrak}{\mathfrak{su}}
\newcommand{\Cscr}{\mathscr{C}}

\newcommand{\alg}{\textrm{alg}}

\newcommand{\Cstar}{C$^*$-}
\newcommand{\eps}{\varepsilon}
\newcommand{\half}{\frac{1}{2}}
\newcommand{\hyph}{-\penalty0\hskip0pt\relax} 			

\newcommand{\ra}{\rightarrow}

\newcommand{\sstar}{\textrm{s}^*\textrm{-}}
\newcommand{\SUqt}{\textrm{SU}_q(2)}
\newcommand{\ten}{\otimes}
\newcommand{\unit}{\mathds{1}}
\newcommand{\vtr}{\vartriangleright}

\newcommand{\RN}[1]{
  \textup{\uppercase\expandafter{\romannumeral#1}}%
}

\DeclareMathOperator*{\bds}{l^\infty-\bigoplus} 	
\DeclareMathOperator*{\ds}{\bigoplus}			
\DeclareMathOperator{\End}{End}				
\DeclareMathOperator{\Hom}{Hom}				
\DeclareMathOperator{\im}{Im}				
\DeclareMathOperator{\Irr}{Irr}				
\DeclareMathOperator{\nat}{Nat}				
\DeclareMathOperator{\Ob}{Ob}				
\DeclareMathOperator{\Rep}{Rep}				
\DeclareMathOperator{\sgn}{sgn}				
\DeclareMathOperator{\supp}{supp}			
\DeclareMathOperator{\Span}{span}			
\DeclareMathOperator{\Tr}{Tr}				
\DeclareMathOperator{\tr}{tr}				
\DeclareMathOperator*{\vds}{c_0-\bigoplus} 		

\begin{document}
\thispagestyle{plain}
\begin{center}
{\Large Convergence to the boundary for random walks on discrete quantum groups and monoidal categories}

\bigskip
{\large Bas P.\ A.\ Jordans\footnote{E-mail address: bpjordan@math.uio.no\\
Department of Mathematics, University of Oslo, P.O. Box 1053 Blindern, 0316 Oslo, Norway. \\
The research leading to these results has received funding from the European Research Council under the European
Union’s Seventh Framework Programme (FP/2007-2013) / ERC Grant Agreement no. 307663 (PI: S. Neshveyev).\\
Date: \today}
}
\end{center}

\bigskip

\begin{abstract}
We study the problem of convergence to the boundary in the setting of random walks on discrete quantum groups. Convergence to the boundary is established for random walks on $\widehat{\textrm{SU}_q(2)}$. Furthermore, we will define the Martin boundary for random walks on C$^*$-tensor categories and give a formulation for convergence to the boundary for such random walks. These categorical definitions are shown to be compatible with the definitions in the quantum group case. This implies that convergence to the boundary for random walks on quantum groups is stable under monoidal equivalence.
\end{abstract}

\section{Introduction}
When studying random walks several natural questions arise. Among these are ``what is the asymptotic behaviour of the random walk as the time tends to infinity?'' and ``can we describe all invariant functions on the space?'' To be precise, given a discrete Markov chain $(X,P)$, consider the space of all infinite paths $\Omega:=X^\Nat$ with associated coordinate projections $X_n\colon\Omega\ra X$, can we describe the behaviour of $X_n$ as $n\ra\infty$? And can we find all functions $h\colon X\ra\Com$ such that $Ph=h$? It turns out that the answers to these two questions are related.
Martin \cite{Martin41} defined a compactification $\tilde M(X,P)$ of $X$ with respect to $P$ (nowadays called the Martin compactification) and a boundary $M(X,P):=\tilde M(X,P)\setminus X$. He proved that every positive harmonic function can be represented by an integral over this boundary, hereby partially answering the second question. This result is the probabilistic analogue of the theorem that any analytic function on the disc can be represented by an integral over the circle.
Some years later Doob \cite{Doob59} and Hunt \cite{Hunt60} independently showed existence of a measurable function $X_\infty\colon\Omega\ra M(X,P)$ such that the coordinate maps $X_n$ converge to $X_\infty$. Note that this function $X_\infty$ takes values in the Martin boundary, explaining the terminology ``convergence to the boundary''. Moreover they strengthened Martin's result and proved by means of this convergence result that positive harmonic functions can uniquely be represented on a smaller subset of the Martin boundary, the so-called minimal boundary. Let $\nu^\unit$ be the measure on the Martin boundary which represents the constant function $1$. It can be shown that if $h$ is a bounded harmonic function, then the corresponding representing measure $\nu^h$ is absolutely continuous with respect to $\nu^\unit$. The Martin boundary together with this measure $\nu^\unit$ is called the Poisson boundary. It describes all positive bounded harmonic functions.\\
In the early nineties Biane \cite{Biane91, Biane92a, Biane92b, Biane94} started the study of non-commutative random walks on duals of compact groups. His idea was to work on the group von Neumann algebra $L(G)$ and act with operators of the form $P_\varphi:=(\varphi\ten\iota)\Delta$, where $\Delta$ is the comultiplication given by $\Delta(\lambda_g)=\lambda_g\ten\lambda_g$. These operators $P_\varphi$ form the analogues of Markov operators used in Markov chains on discrete spaces.
Izumi \cite{Izumi02} continued this study and defined random walks on discrete quantum groups. His main motivation was to study the (non-)minimality of actions of compact quantum groups on von Neumann algebras. He considered actions on infinite tensor products and proved that the relative commutant of the action can be described as the space of harmonic elements of a Markov operator $P$. The space of $P$-harmonic elements is called the Poisson boundary and forms the non-commutative analogue of the Poisson boundary for classical random walks. Neshveyev and Tuset \cite{NeshveyevTuset04} built further on this story and defined the Martin boundary for noncommutative random walks on discrete quantum groups. In that paper they proved the important result that any positive harmonic element can be represented by a linear functional on the boundary. But the problem of convergence to the boundary remained open. In another overview paper \cite{NeshveyevTuset03} they gave a conjecture of what convergence to the boundary should correspond to in the quantum world.

The problem with proving boundary convergence in the noncommutative setting is that the ``commutative proof'' is very hard to translate. Classically stopping times and Martingale convergence theorems are used to obtain almost everywhere convergence. However, it is not clear how to formulate such stopping times, up- and downcrossings in a noncommutative way. In this paper we prove that the conjecture of convergence to the boundary as proposed in \cite{NeshveyevTuset03} holds for $\SUqt$. Our approach is very computational, but it shows that there is exponential fast convergence, which is a lot faster than what occurs classically.

There are not many examples of quantum groups for which the Martin boundary has been computed. The difficulty is that it is very hard to compute how the Martin kernel acts on matrix units of the discrete dual $\hat G$. One can say that $\SUqt$ is the main example for which the Martin boundary has been identified \cite{NeshveyevTuset04}. Other cases are based on this computation of $\SUqt$. For example, the free orthogonal quantum groups \cite{VaesVergnioux07}, \cite{VaesVanderVennet08} which are monoidally equivalent to $\SUqt$. A theorem by De Rijdt and Vander Vennet \cite{DeRijdtVanderVennet10} generalises this and gives a concrete method to compute the Poisson and Martin boundaries for a quantum group out of the Poisson and Martin boundaries of a monoidally equivalent quantum group. Later Neshveyev and Yamashita \cite{NeshveyevYamashita14c} proved an extended Tannaka--Krein duality between certain algebras with $G$-$\hat G$ actions and \Cstar tensor categories related to $\Rep(G)$. With this correspondence they were able to define a categorical version of the Poisson boundary. The results of both papers indicate that a similar result should hold for the Martin boundary too. In this paper we will show that this is indeed the case. We give a definition of a categorical Martin boundary and show that under the correspondence of \cite{NeshveyevYamashita14c} the Martin boundary of a random walk on a discrete quantum group can be reconstructed from the categorical Martin boundary of the random walk on the representation category. A natural question is then whether also convergence to the boundary can be presented in a categorical way. Fortunately the answer to this question appears to be positive, yielding a broader range of examples for which convergence to the boundary holds.

\bigskip
The paper is structured in the following way. We start with the preliminaries in which we build up the theory about compact and discrete quantum groups, \Cstar tensor categories and noncommutative random walks. In section \ref{sec_convergence_SUQt} we establish convergence to the boundary for random walks on $\widehat{\SUqt}$. In section \ref{sec_monoidal_equivalence} we will briefly outline how the methods of \cite{DeRijdtVanderVennet10} can be extended to show that convergence to the boundary is invariant under monoidal equivalence of quantum groups. We proceed by defining the Martin boundary and convergence to the boundary in a categorical way. In the final section we establish compatibility of this categorical description with the theory for random walks on discrete quantum groups.

\bigskip
\noindent{\bf Acknowledgements} I would like to thank my supervisor prof. S. Neshveyev for his guidance and advice when working on the problems discussed in this paper.

\section{Preliminaries}
\subsection{Compact and discrete quantum groups}
There are many good presentations of the theory of compact and discrete quantum groups see for example \cite{NeshveyevTuset13}, \cite{KustermansTuset99}, \cite{VanDaele96}. Here we give a brief overview of the concepts will we need.

\begin{Def}
A {\it compact quantum group} is a pair $G=(C(G),\Delta)$ consisting of a unital \Cstar algebra $C(G)$ and a $*$-homomorphism $\Delta\colon C(G)\ra C(G)\ten C(G)$ called the {\it comultiplication}, satisfying coassociativity $(\iota\ten\Delta)\Delta=(\Delta\ten\iota)\Delta$ and the cancellation property
\[
\Span\{(a\ten 1)\Delta(b)\,:\, a,\,b\in C(G)\} \textrm{ and } \Span\{(1\ten a)\Delta(b)\,:\, a,\,b\in C(G)\}
\]
are norm-dense in $C(G)\ten C(G)$.
\end{Def}
The tensor product $\ten$ of \Cstar algebras always indicates the minimal tensor product unless mentioned otherwise.

There exists a unique state $h\colon C(G)\ra\Com$, called the {\it Haar state} which satisfies
\[
(h\ten \iota)\Delta(a) = h(a)1 = (\iota\ten h)\Delta(a).
\]
This is called right invariance, respectively left invariance of the Haar state. Throughout the paper we will assume that the Haar state is faithful, so we are dealing with {\it reduced quantum groups}.

Often Sweedler's sumless notation will be useful. If $x=\sum_i x^{(1)}_i\ten x^{(2)}_i$ we write for example $x_{12}=  x_{1,2}=x\ten 1$ or $x_{13}=x_{1,3}=\sum_i x^{(1)}_i\ten 1
\ten x^{(2)}_i$.

\begin{Def}
A {\it unitary representation} of $G$ on a Hilbert space $\Hcal$ is a unitary element $U\in M(C(G)\ten B_0(\Hcal))$, such that
\begin{align}\label{def_rep_eq1}
(\Delta\ten\iota)U=U_{13}U_{23}.
\end{align}
Here $B_0(\Hcal)$ denotes the compact operators on $\Hcal$.
\end{Def}
We will only deal with unitary representations and therefore we will usually omit the prefix ``unitary''. If $U$ and $V$ are unitary representations, denote the space of {\it intertwiners}
\[
\Hom(U,V):=\{T\in B(\Hcal_U,\Hcal_V)\,:\, (\iota\ten T) U=V(\iota\ten T)\}.
\]
A representation is called {\it irreducible} if $\Hom(U,U) = \Com 1_{\Hcal_U}$. Two representations $U$ and $V$ are {\it (unitarily) equivalent} if $\Hom(U,V)$ contains a unitary intertwiner, in which case we write $U\cong V$. The set of equivalence classes of irreducible representations is denoted by $\Irr(G)$. We fix representatives $U_s$ for every $s\in\Irr(G)$. The equivalence class of the trivial representation, the representation on $\Com$ given by $1\in C(G)\ten B(\Com)$, is denoted by $0$.\footnote{Note that in literature the label $0$ is sometimes used for the $0$-dimensional representation} The (classical) dimension of $U$ is by definition the dimension of $\Hcal_U$ and is denoted by $\dim(U)$. Sometimes, due to the leg numbering, we put the $s$ as a superscript, so we write for example $U_{12}^s:= (U_s)_{12}$.

Every irreducible representation is finite dimensional and every finite dimensional representation decomposes into a direct sum of irreducible representations. So in most cases it suffices to deal with finite dimensional representations.
For two representations $U$ and $V$ define the tensor product representation on $\Hcal_U\ten\Hcal_V$ by $U\times V:= U_{12}V_{13}$. Note that in general $U\times V \not\cong V\times U$. For $s,\,t\in\Irr(G)$ we write $s\ten t$ when we want to indicate the tensor product $U_s\times U_t$. For example, $\Hom(s\ten t,r)=\Hom(U_s\times U_t, U_r)$.

\bigskip
Let $\Hcal$ be a Hilbert space. The inner product is denoted $\langle\cdot,\cdot\rangle$ and is antilinear in the first component. Identify the dual space $\Hcal^*$ with the complex conjugate Hilbert space $\bar\Hcal$. So if $\xi\in\Hcal$, then $\bar\xi\in\bar\Hcal$ acts as $\bar\xi(\zeta)=\langle\xi,\zeta\rangle$. Consider the map $j\colon B(\Hcal)\ra B(\bar\Hcal)$, $j(T)\bar\xi:=\overline{T^*\xi}$. So under the identification $j$ maps $T$ to its dual operator.
\begin{Def}\label{Def_conjugate_rep}
Let $U$ be a finite dimensional representation. The {\it contragredient representation} $U^c$ on $\Hcal_U^*$ is defined by
\[
U^c:= (\iota\ten j)(U^{-1})\in C(G)\ten B(\Hcal^*).
\]
Then $U^c\in C(G)\ten B(\Hcal_U^*)$ is invertible and satisfies $(\Delta\ten\iota)(U^c)=U^c_{13}U^c_{23}$, although it need not be unitary. If $U$ is irreducible there exists a unique positive invertible operator $\rho_U\in\Hom(U,U^{cc})$ such that
\[
\Tr(x\rho_U)=\Tr(x\rho_U^{-1}),\qquad\textrm{for all } x\in\Hom(U,U).
\]
The {\it conjugate representation} $\bar U$ of $U$ is defined as
\[
\bar U:= (1\ten j(\rho_U)^{\half})U^c(1\ten j(\rho_U)^{-\half})\in C(G)\ten B(\bar\Hcal_U)
\]
this is again a unitary representation. The {\it quantum dimension} of $U$ equals $\Tr(\rho_U)$ and is denoted by $d_U$. Note that $\dim(U)\leq d_U$ and equality holds if and only if $\rho_U=\iota_U$.
\end{Def}

Frobenius reciprocity holds and can be formulated as
\[
\Hom(U_1,U_3\times \bar U_2)\cong\Hom(U_1\times U_2, U_3)\cong \Hom(U_2, \bar U_1\times U_3),
\]
where $U_1$, $U_2$ and $U_3$ are representations. In particular if $U$ is irreducible, then so is $\bar U$. We write $\bar s\in\Irr(G)$ for the unique representative satisfying $\bar U_s=U_{\bar s}$.

\begin{Not}
Given a collection $\{X_i\}_{i\in I}$ of Banach spaces for some index set $I$, we use the following conventions:
\begin{align*}
\prod_{i\in I} X_i &:= \{(x_i)_i\,:\, x_i\in X_i\};\\
\bigoplus_{i\in I} X_i &:= \big\{(x_i)_i \in \prod_{i\in I} X_i \,:\, x_i\neq 0 \textrm{ for at most finitely many } i\big\};\\
\vds_{i\in I} X_i &:= \big\{(x_i)_i \in \prod_{i\in I} X_i \,:\, (\|x_i\|)_i \in c_0(I)\big\};\\
\bds_{i\in I} X_i &:= \big\{(x_i)_i \in \prod_{i\in I} X_i \,:\, \sup_{i\in I}\{\|x_i\|\} <\infty\big\}.
\end{align*}
Note that the last two algebras are Banach spaces, with respect to the norm $\|(x_i)_i\|:=\sup_{i}\|x_i\|$. The norm-closure of $\bigoplus_{i\in I} X_i$ in $\bds_{i\in I} X_i$ equals $\vds_i X_i$. Furthermore if all $X_i=A_i$ are \Cstar algebras, then $\vds_{i\in I} A_i$ and $\bds_{i\in I} A_i$ are \Cstar algebras. Moreover if all $A_i$ are unital for the multiplier algebra the identity $M(\vds_i A_i)=\bds_i A_i$ holds.
\end{Not}

\begin{Not}
The matrix coefficients of $G$ are defined as
\[
\Com[G]:=\{(\iota\ten\omega)(U)\,:\, U \textrm{ f.d. representation},\, \omega\in (B(\Hcal_U))^*\}.
\]
If $\xi,\zeta\in\Hcal_U$, denote the matrix unit $m_{\xi,\zeta}:=\langle\zeta,\cdot\rangle\xi\in B(\Hcal_U)$. For each $s\in\Irr(G)$ we fix once and for all an orthonormal basis $\{\xi_i^s\}_{i=1}^{\dim(s)}$ in $\Hcal_s$ such that $\rho_s$ acts diagonal with respect to this basis. Abbreviate $m_{ij}^s:=m_{\xi_i^s,\xi_j^s}\in B(\Hcal_s)$, so $m_{ij}(\xi_k)=\delta_{jk}\xi_i$. Write $U_s=\sum_{i,j} u_{ij}^s\ten m_{ij}^s = (u_{ij}^s)_{i,j}$ and $\rho_s=\sum_{ij}(\rho_s)_{ij}\,m_{ij}^s=\sum_i (\rho_s)_{ii}\,m_{ii}^s$. Note that \eqref{def_rep_eq1} reads as $\Delta(u_{ij}^s)=\sum_{k=1}^{\dim(U_s)}u_{ik}^s\ten u_{kj}^s$. The matrix coefficients satisfy the following orthogonality relations:
\begin{align}
h((u_{kl}^t)^*u_{ij}^s) = \delta_{st}\delta_{jl}\frac{(\rho_s^{-1})_{ik}}{d_s}, & & h(u_{kl}^t (u_{ij}^s)^*) = \delta_{st}\delta_{ki}\frac{(\rho_s)_{jl}}{d_s}. \label{orthogonality_relations_eq1}
\end{align}
\end{Not}

\begin{Def}
A {\it Hopf $*$-algebra} is a pair $(\Acal,\Delta)$ consisting of a unital $*$-algebra $\Acal$ and a unital $*$-homomorphism $\Delta\colon\Acal\ra\Acal\ten\Acal$, the {\it comultiplication}, satisfying $(\Delta\ten\iota)\Delta=(\iota\ten\Delta)\Delta$, together with linear maps $S\colon\Acal\ra \Acal$ and $\eps\colon\Acal\ra\Com$ such that the identities
\begin{equation}\label{def_Hopf_algebra_eq1}
(\eps\ten\iota)\Delta(a)=(\iota\ten\eps)\Delta(a)=a \quad \textrm{and} \quad m(S\ten\iota)\Delta(a) = m(\iota\ten S)\Delta(a) = \eps(a)1,
\end{equation}
hold for all $a\in\Acal$. Here $m\colon\Acal\ten\Acal\ra\Acal$ is the multiplication map. $S$ is called the {\it antipode} and $\eps$ the {\it counit}.
\end{Def}

The maps $S$ and $\eps$ are uniquely determined by \eqref{def_Hopf_algebra_eq1}. Moreover $\eps$ is a $*$-homomorphism, $S$ is an anti-homomorphism and they satisfy the following relations:
\begin{align*}
\eps(S(a))= \eps(a); && \Delta S= (S\ten S)\sigma\Delta; && S(S(a^*)^*)=a, && (a\in\Acal),
\end{align*}
where $\sigma$ is the flip map.

\bigskip
Note that $(\Com[G],\Delta)$ is a Hopf $*$-algebra with maps $S$, $\eps$ defined by the identities
\[
(S\ten\iota)(U)=U^* \quad \textrm{and} \quad (\eps\ten\iota)(U)=1.
\]
Moreover $\Com[G]$ is dense in $C(G)$.

\bigskip
Consider the space of linear functionals $\Com[G]^*$. For every finite dimensional representation $U$ of $G$ we get a representation $\pi_U$ of $\Com[G]^*$ on $\Hcal_U$ by
\begin{equation}\label{def_dual_eq1}
\pi_U\colon\Com[G]^*\ra B(\Hcal_U),\qquad \omega\mapsto (\omega\ten\iota)U.
\end{equation}
The collection $\{\pi_{U_s}\}_{s\in\Irr(G)}$ defines an isomorphism $\Com[G]^*\cong \prod_s B(\Hcal_s)=:\Ucal(G)$. Similarly $(\Com[G]^{\ten n})^*\cong \prod_{s_1,\ldots,s_n} B(\Hcal_{s_1}\ten\ldots\ten \Hcal_{s_n})=:\Ucal(G^n)$. Define a unital $*$-morphism $\hat\Delta\colon\Com[G]^*\ra(\Com[G]^{\ten 2})^*$ dual to the multiplication on $G$
\[
\hat\Delta(\omega)(a\ten b):=\omega(ab),\qquad \textrm{for } \omega\in\Com[G]^* \textrm{ and } a,b\in \Com[G].
\]
Then $(\Com[G]^*,\hat\Delta)$ satisfies the axioms of a Hopf $*$-algebra (when understood properly) with counit $\hat\eps(\omega)=\omega(1)$ and antipode $\hat S(\omega)=\omega(S(\cdot))$ whenever $\omega\in\Com[G]^*$. Via the isomorphisms above this leads to a map $\hat\Delta\colon\Ucal(G)\ra\Ucal(G^2)$. We will generally use $\Ucal(G)$ instead of $\Com[G]^*$. Let $\pi_s\colon\Ucal(G)\ra B(\Hcal_s)$ denote the projection on the matrix block corresponding to $s$. Equivalently one can define the comultiplication $\hat\Delta(a)$ for $a\in\Ucal(G)$ by
\begin{equation}\label{notation_CQG_eq1}
(\pi_s\ten\pi_t)(\hat\Delta(a))T=T\pi_r(a), \qquad\textrm{for all } T\in\Hom(r,s\ten t).
\end{equation}

Another way goes via the multiplicative unitary. For this consider the GNS representation $(\Hcal_h,\xi_h,\pi_h)$ associated to the Haar state $h$. We will write $L^2(G)$ for $\Hcal_h$. If $x\in C(G)$ write $\Lambda(x)=x\xi_h\in L^2(G)$ for the corresponding vector in the GNS construction. Since $h$ is faithful, we can identify $C(G)$ with $\pi_h(C(G))$ inside $B(L^2(G))$. Denote by $L^\infty(G)$ the von Neumann algebra generated by $\pi_h(C(G))$ in $B(L^2(G))$. By density of $C(G)$ inside $L^2(G)$ the map
\begin{equation}\label{def_multiplicative_unitary_eq1}
(\xi\ten a\xi_h)\mapsto \Delta(a)(\xi\ten\xi_h)
\end{equation}
extends to a unitary operator in $B(L^2(G)\ten L^2(G))$. The adjoint $W$ of \eqref{def_multiplicative_unitary_eq1} is called the {\it multiplicative unitary} and defines a unitary (not necessarily finite dimensional) representation of $G$ on $L^2(G)$, the so-called {\it left regular representation}.
Given $\xi,\,\zeta\in\Hcal$ denote the functional $\omega_{\xi,\zeta}\colon B(\Hcal)\ra\Com$, $\omega_{\xi,\zeta}(T):=\langle\xi,T\zeta\rangle$. If $U_s$ is an irreducible representation, then for every $\zeta\in\Hcal_s$ the map
\[
\theta_\zeta\colon\Hcal_s\ra L^2(G),\qquad \xi\mapsto (d_s)^\half(\Lambda\ten\omega_{\rho_s^{-\half}\zeta,\xi})(U_s^*)
\]
intertwines $U_s$ and $W$. From the orthogonality relations it follows that $\theta_\zeta$ is isometric if $\zeta$ is a unit vector. Moreover if $\zeta,\,\zeta'\in\Hcal_s$ are two orthogonal vectors the corresponding images of $\theta_\zeta$ and $\theta_{\zeta'}$ are orthogonal. Thus by picking orthonormal bases we obtain a canonical inclusion $\bar\Hcal_s\ten\Hcal_s\hookrightarrow L^2(G)$, it corresponds to the space of matrix coefficients of $U_s$. By identifying $B(\Hcal_s)\cong \bar\Hcal_s\ten\Hcal_s$ we obtain an inclusion $\ds_s B(\Hcal_s)\hookrightarrow L^2(G)$. Taking all irreducible representations exhausts $L^2(G)$. This gives the Peter--Weyl decomposition for compact quantum groups. Observe that $W$ can be expressed as
\begin{equation}\label{multiplicative_unitary_eq1}
W=\bigoplus_{s\in\Irr(G)} U_s = \sum_s \sum_{i,j} u_{ij}^s\ten m_{ij}^s.
\end{equation}
The multiplicative unitary satisfies the {\it pentagon equation} $W_{12}W_{13}W_{23}=W_{23}W_{12}$.

\begin{Not}
Let $G=(C(G),\Delta)$ be a compact quantum group. The discrete quantum group dual to $G$ is the virtual object indicated by $(\hat G,\hat\Delta)$. We write
\begin{align*}
c_{00}(\hat{G})&:= \ds_{s\in\Irr(G)} B(\Hcal_s), & c_0(\hat{G})&:= \vds_{s\in \Irr(G)} B(\Hcal_s), \\
l^\infty(\hat G)&:= \bds_{s\in \Irr(G)} B(\Hcal_s), &\Ucal(G)&:= \prod_{s\in \Irr(G)} B(\Hcal_s).
\end{align*}
The isomorphism \eqref{def_dual_eq1} shows that $\hat G$ is indeed dual to the compact quantum group $G$.
The multiplicative unitary $W$ encodes all information of the quantum group $G$. Namely its matrix coefficients span both $\Com[G]$ and $c_{00}(\hat G)$ and the comultiplications are given by
\begin{align}
\Delta(a)&=W^*(1\ten a)W\in C(G)\ten C(G); \notag\\
\hat\Delta(x)&=W(x\ten 1)W^*\in l^\infty(\hat G)\bar\ten l^\infty(\hat G), \label{def_dual_eq3}
\end{align}
for $a\in C(G)$ and $x\in l^\infty(\hat G)$. We will not need it, but also the antipode and counit can be expressed in terms of $W$ (cf.\ \cite[Prop.\ 11.37]{KlimykSchmudgen97}). Recall the elements $\rho_s$ defined in Definition \ref{Def_conjugate_rep}. We write $\rho\in\Ucal(G)$ for the element that satisfies $\pi_s(\rho)=\rho_s$ for every $s$.

The comultiplication $\hat\Delta$ has unique right- and left-invariant weights, denoted $\hat\psi$ and respectively $\hat\phi$. Invariance means that they satisfy
\begin{align}
\hat\psi((\iota\ten\omega)\hat\Delta(x))=\hat\psi(x)\omega(1); && \hat\phi((\omega\ten\iota)\hat\Delta(y)) = \hat\phi(y)\omega(1), \label{Haar_weights_eq1}
\end{align}
for all $x,y\in l^\infty(\hat G)^+$ such that $\hat\psi(x)<\infty$, $\hat\phi(y)<\infty$ and all $\omega$ positive normal linear functional on $l^\infty(\hat G)$. These weights can be written down explicitly as
\begin{align}
\hat\psi(x) = \sum_{s\in\Irr(G)} d_s\Tr(\pi_s(x\rho^{-1})); \qquad \hat\phi(y) = \sum_{s\in\Irr(G)} d_s\Tr(\pi_s(y\rho)). \label{Haar_weights_eq2}
\end{align}
Note that $\hat\psi$ and $\hat\phi$ are unbounded on $c_0(\hat G)$. The modular groups are given by
\begin{align}
\sigma_t^{\hat\psi}(x) = \rho^{-it}x\rho^{it}; && \sigma_t^{\hat\phi}(x) = \rho^{it}x\rho^{-it}, && (x\in c_{00}(\hat G)). \label{Haar_weights_eq3}
\end{align}
\end{Not}

\begin{Def}
Let $G$, $\hat G$ be a compact respectively discrete quantum group and $B$ a unital \Cstar algebra. A {\it left action} of $G$ (resp.\ $\hat G$) on $B$ is a unital $*$-homomorphism $\alpha\colon B\ra C(G)\ten B$ (resp.\ $\alpha\colon B\ra M(c_0(\hat G)\ten B)$) satisfying $(\iota\ten\alpha)\alpha=(\Delta\ten\iota)\alpha$ (resp.\ $(\iota\ten\alpha)\alpha=(\hat\Delta\ten\iota)\alpha$) and $\alpha(B)(C(G)\ten1)$ is norm dense in $C(G)\ten B$ (resp.\ $\alpha(B)(c_0(\hat G)\ten 1)$ is norm dense in $c_0(\hat G)\ten B$). {\it Right actions} are defined similarly.
\end{Def}

\begin{Def}
Let $G$, $\hat G$ be a compact respectively discrete quantum group and $N$ a von Neumann algebra. A {\it right action} of $G$ (resp.\ $\hat G$) on $N$ is an injective normal unital $*$-homomorphism $\beta\colon N\ra N\bar\ten L^\infty(G)\ten B$ (resp.\ $\beta\colon N\ra N\bar\ten l^\infty(\hat G)$) satisfying $(\beta\ten\iota)\beta=(\iota\ten\Delta)\beta$ (resp.\ $(\beta\ten\iota)\beta=(\iota\ten\hat\Delta)\beta$). {\it Left actions} are defined similarly.
\end{Def}

For us the most important actions are the adjoint actions.
\begin{Not}
The multiplicative unitary $W$ defines left and right {\it adjoint actions} of $G$ and $\hat G$ by
\begin{align*}
\alpha_l\colon l^\infty(\hat G)&\ra M(C(G)\ten l^\infty(\hat G)),&	 x&\mapsto W^*(1\ten x)W;\\
\alpha_r\colon C(G)&\ra M(C(G)\ten l^\infty(\hat G)),&	 a&\mapsto W(a\ten 1)W^*.
\end{align*}
\end{Not}

For $s_1,\ldots,s_n,t\in \Irr(G)$ denote by $m_{s_1\ldots,s_n}^t:=\dim(\Hom(U_t,U_{s_1}\times\cdots\times U_{s_n}))$ the {\it multiplicity} of the representation $U_t$ in $U_{s_1}\ten\ldots\ten U_{s_n}$.

\begin{Lem}\label{estimates_multiplicities}
For $s_1,\ldots,s_n,t\in \Irr(G)$ the multiplicities satisfy the inequalities
\begin{align}
m_{s_1,\ldots,s_n}^t &\leq \frac{\dim(U_{s_1})\cdots\dim(U_{s_n})}{\dim(U_t)} \label{estimates_multiplicities_eq1};\\
\sum_{t\in \Irr(G)} m_{s_1,\ldots,s_n}^t &\leq \dim(U_{s_1})\cdots\dim(U_{s_n}). \label{estimates_multiplicities_eq2}
\end{align}
\end{Lem}
\begin{pf}
For $s_1,\ldots,s_n\in\Irr(G)$ it holds that $\bigoplus_r m_{s_1,\ldots,s_n}^r U_r \cong U_{s_1} \times \ldots \times U_{s_n}$. Thus
\[
\sum_{r\in\Irr(G)} m_{s_1,\ldots,s_n}^r\dim(U_r) = \dim(U_{s_1})\cdots \dim(U_{s_n}).
\]
From which \eqref{estimates_multiplicities_eq1} follows immediately. As each representation has dimension $\geq 1$ also \eqref{estimates_multiplicities_eq2} is clear from this identity.
\end{pf}

For KMS states we follow \cite{BratteliRobinson97}. The following result is easy to prove using the elementary properties of the modular conjugation $J$.

\begin{Lem}\label{lem_inner_prod}
Let $\eta\colon A\ra\Com$ be a KMS state on a \Cstar algebra $A$. Then $\eta$ extends to a normal faithful state on the von Neumann algebra $M:=\pi_\eta(A)''$. Let $\sigma^\eta = \sigma$ be the associated modular group, $J_\eta$ the modular conjugation and $\xi_\eta$ the GNS vector. Then $\eta\colon M\ra\Com$ is a $\sigma$-KMS state and
\begin{enumeraterm}
\item the sesquilinear form $(\cdot,\cdot)_\eta\colon A\times A\ra\Com$, $(a,b)_\eta:= \eta(b\sigma_{-\frac{i}{2}}(a^*))= \langle \xi_\eta,bJ_\eta a\xi_\eta\rangle$ is a semi-inner product;
\item the sesquilinear form $(\cdot,\cdot)_\eta\colon M\times M\ra\Com$, $(a,b)_\eta:= \eta(b\sigma_{-\frac{i}{2}}(a^*))= \langle \xi_\eta,bJ_\eta a\xi_\eta\rangle$ is an inner product;
\item the linear functionals $(\cdot,c)_\eta\colon M\ra\Com$ and $(c,\cdot)_\eta\colon M\ra\Com$ are positive if $c\in M$ is positive;
\end{enumeraterm}
\end{Lem}
This form above is a modification of the well-known semi-inner product $\langle \cdot,\cdot\rangle_\eta\colon A\times A\ra\Com$, $\langle a,b\rangle_\eta:= \eta(ab^*)$. The modular group ensures that positivity result in (iii) holds. Note that $\sigma_{-\frac{i}{2}}(b^*)$ need not be defined for all $b\in A$, but the inner product $\langle \xi_\eta,bJ_\eta a\xi_\eta\rangle$ is well-defined for all $a,b\in A$.

We will freely use the results of the following lemma.
\begin{Lem}\label{Lem_q-numbers}
The $q$-numbers defined by
\[
[n]_q:=\frac{q^n-q^{-n}}{q-q^{-1}}
\]
satisfy the following identities:
\begin{enumerate}[label=(\roman*)]
\item $[n]_q= q^{n-1} + q^{n-3}+\ldots + q^{-n+1}$;
\item $[m+n]_q = q^n[m]_q + q^{-m}[n]_q = q^{-n}[m]_q+q^m[n]_q$;
\item $\frac{[m]_q}{[n]_q} = q^{n-m}(1+O(q^{2m})+O(q^{2n}))$ as $m,\,n\ra\infty$.
\end{enumerate}
\end{Lem}
\begin{pf}
(i) and (ii) are trivial. For (iii)
\[
\frac{[m]_q}{[n]_q}q^{m-n}-1 = \frac{q^{2m}-1}{q^{2n}-1} - \frac{q^{2n}-1}{q^{2n}-1} = \frac{q^{2m}-q^{2n}}{q^{2n}-1} = O(q^{2m}) + O(q^{2n}),
\]
which proves the lemma.
\end{pf}

\subsection{\texorpdfstring{\Cstar}{C*} tensor categories}
For the definitions and elementary results regarding \Cstar (tensor) categories we refer to \cite[Ch.\ 2]{NeshveyevTuset13}. We will follow their notation and conventions.

\begin{Assum}
We will always assume that our \Cstar categories are small and closed under finite direct sums and subobjects. Moreover, \Cstar tensor categories are in addition assumed to be strict and rigid.
\end{Assum}

Notice that simplicity of the unit object is not assumed. To be clear we fix some notation.

\begin{Not}
Let $\Ccal$ be a strict \Cstar tensor category with simple unit $\unit$. By $\Ob(\Ccal)$ we denote the objects of $\Ccal$ and by $\Hom_\Ccal(U,V)$ the set of morphisms between $U$ and $V$. Let $\Irr(\Ccal)$ denote the equivalence classes of the simple objects in $\Ccal$. For each $s\in\Irr(\Ccal)$ fix a representative $U_s$. For each $U_s$ fix a conjugate object $\bar U_s$. By Frobenius reciprocity $\bar U_s$ is again simple, thus isomorphic to $U_t$ for some $t\in\Irr{\Ccal}$. We define a map $\Irr(\Ccal)\ra\Irr(\Ccal)$, $s\mapsto \bar{s}$, where $\bar{s}$ is defined by the identity $\bar U_s=U_{\bar s}$. $0\in\Irr(\Ccal)$ indicates the unit object, thus $U_0=\unit$. The multiplicity of $U_t$ in $U\ten V$ is denoted by $m_{U,V}^t=\dim(\Hom_\Ccal(U_t,U\ten V))$. Thus $U\ten V \cong \bigoplus_t m_{U,V}^t U_t$, where $m_{U,V}^t U_t$ means the direct sum of $m_{U,V}^t$ copies of $U_t$.
For each object $U$ we let $(R_U,\bar R_U)$ be a standard solution of the conjugate equations. Let $d_U:=\|R_U\|^2$ be the {\it intrinsic dimension} of $U$. The {\it normalized categorical traces} are given by
\[
\tr_U\colon\End_\Ccal(U)\ra\Com,\qquad T\mapsto d_U^{-1}R_U^*(\iota_{\bar U}\ten T)R_U = d_U^{-1}\bar R_U^*(T\ten\iota_{\bar U})\bar R_U.
\]
The adjective ``normalized'', refers to the constant $d_U^{-1}$, to make sure that $\tr_U(\iota_U)=1$. If $\Ccal=\Rep(G)$, then it can be shown that $R_s:= R_{U_s}$ and $\bar R_s:=\bar R_{U_s}$ are given by
\begin{align}
R_s(1):=\sum_i \bar\xi_i^s\ten\rho_s^{-\half}\xi_i^s & &\bar R_s(1):=\sum_i \rho_s^\half\xi_i^s\ten\bar\xi_i^s. \label{Not_conjugate_equations_eq1}
\end{align}
In particular if $U$ is multiplicity free it holds that $\varphi_U(x)=\tr_U(x)$ for $x\in\End_{\Rep(G)}(U)\subset l^\infty(\hat G)$ (see \cite[Rem.\ 2.2.15]{NeshveyevTuset13}). This also implies that the intrinsic dimension of a representation equals the quantum dimension.
\end{Not}

\subsection{Infinite tensor products}\label{subsec_infinite_tensor_products}
For the construction of infinite tensor products of von Neumann algebras we follow \cite[\textsection~XIV.1]{Takesaki03III}. This construction makes use of infinite tensor products of \Cstar algebras. Given a sequence of \Cstar algebras $(A_n)_n$ with $*$-homomorphisms $\pi_n\colon A_n\ra A_{n+1}$, there exists an {\it inductive limit}
\[
A :=\lim_{\longrightarrow}(A_n,\pi_n),
\]
with inclusion maps $\iota_n\colon A_n\ra A$. This \Cstar algebra $A$ can be described by the following universal property. If $B$ is a \Cstar algebra with maps $\alpha_n\colon A_n\ra\ B$ such that for each $n\in\Nat$ the diagram
\[
\xymatrix{A_n \ar[d]_{\pi_n} \ar[dr]^{\alpha_n}  \\
A_{n+1}\ar[r]_{\alpha_{n+1}} & B}
\]
commutes, then there exists a unique map $\alpha\colon A\ra B$ such that for each $n\in\Nat$ the following diagram commutes
\[
\xymatrix{A_n \ar[d]_{\iota_n} \ar[dr]^{\alpha_n}  \\
A\ar[r]_{\alpha} & B}.
\]
The {\it \Cstar algebraic tensor product} of a sequence $(A_n)_n$ of \Cstar algebras is defined as the inductive limit
\[
\bigotimes_{n=1}^\infty A_n :=\lim_{\longrightarrow} (A_n',\pi_n),
\]
where $A_n':= A_1\ten\cdots\ten A_n$ with the minimal tensor product and connecting maps $\pi_n(x):=x\ten 1$. \\
Given a sequence $(M_n)_n$ of von Neumann algebras with normal states $\omega_n\colon M_n\ra\Com$, define the state $\omega=\bigotimes_n\omega_n$ on the \Cstar algebra $A:= \bigotimes_n M_n$ by
\[
\omega(x_1\ten\cdots\ten x_n\ten 1\ten 1\ten\cdots) := \omega_1(x_1)\cdots\omega_n(x_n)
\]
and extension to $A$. The GNS-construction applied to $\omega$ gives a cyclic representation $(\pi_\omega,\Hcal_\omega,\xi_\omega)$ of $A$. Define the {\it von Neumann algebraic infinite tensor product} of $(M_n,\omega_n)_n$ as
\[
M:=\bigotimes_{n=1}^\infty (M_n,\omega_n):=(\pi_\omega(A))''.
\]
The commutant is taken in $B(\Hcal_\omega)$. The functional $\omega$ extends to a state on $M$, which is faithful if each $\omega_i$ is faithful on $\pi_{\omega_i}(M_i)$. This construction depends heavily on the choice of the sequence of states $(\omega_n)$, different choices can give non-isomorphic algebras see for example \cite[Thm.\ XVIII.1.1]{Takesaki03III}. For us however, tensor products of the form $\bigotimes_{-\infty}^{-1}$ are relevant, these are defined in a similar manner.

The formulation of the following lemma is from \cite[Lem.\ 3.4.]{Izumi02}, the proof can be found in \cite[Lem.\ 2]{Connes75}.

\begin{Lem}[Noncommutative martingale convergence theorem]\label{lem_martingale_convergence}
Suppose $M$ is a von Neumann algebra with a normal state $\varphi$. By construction of the infinite tensor product for any $n\in\Nat$ there is the embedding
\[
i_n\colon\bigotimes_{-n}^{-1} (M,\varphi)\hookrightarrow \bigotimes_{-\infty}^{-1} (M,\varphi), \qquad x\mapsto \cdots\ten1\ten1\ten x.
\]
Define slice maps
\[
E_n'\colon \bigotimes_{-\infty}^{-1} (M,\varphi)\ra \bigotimes_{-n}^{-1} (M,\varphi), \qquad x\mapsto (\cdots\ten\varphi\ten\varphi\ten \iota^{\ten n})(x),
\]
Then $E_n:=i_n \circ E_n'$ is the unique $\varphi$-preserving conditional expectation onto $i_n(\bigotimes_{-n}^{-1} (M,\varphi))$. The maps $E_n$ satisfy:
\begin{enumerate}[label=(\roman*)]
\item for every $x\in \bigotimes_{-\infty}^{-1} (M,\varphi)$ it holds that $x=\sstar\lim_n E_n(x)$;
\item if the sequence $(x_n)_n\subset \bigotimes_{-\infty}^{-1} (M,\varphi)$ satisfies $E_n(x_{n+1})=x_n$, then there exists a unique $x\in \bigotimes_{-\infty}^{-1} (M,\varphi)$ such that $x_n=E_n(x)$.
\end{enumerate}
\end{Lem}

\begin{Not}\label{Def_adjoint_on_path_space}
Recall the left adjoint action $\alpha_l\colon l^\infty(\hat G)\ra M(C(G)\ten l^\infty(\hat G))$. This action can be extended to a left action of $G$ on $\bigotimes_{i=-n}^{-1}  l^\infty(\hat G)$ by
\begin{equation}\label{extended_adjoint_action_eq1}
\alpha_l\colon\bigotimes_{-n}^{-1}l^\infty(\hat G)\ra L^\infty(G)\bar\ten\Big(\bigotimes_{-n}^{-1} l^\infty(\hat G)\Big), \qquad x\mapsto W_{1,n+1}^*\cdots W_{1,2}^*(1\ten x) W_{1,2} \cdots W_{1,n+1}.
\end{equation}
Moreover for $x\in \bigotimes_{-\infty}^{-1}(l^\infty(\hat G),\varphi)$ the following limit can be shown to exist in norm (see \cite[\textsection~3]{Izumi02})
\[
\lim_{n\ra\infty}  W_{1,-1}^*\cdots W_{1,-n}^*(1\ten x) W_{1,-n} \cdots W_{1,-1}.
\]
defining an action on the infinite tensor product. The leg numbering here is different than elsewhere in this paper, the $1$ refers to $L^\infty(G)$ and the $-j$ to the $-j$th component of $\bigotimes_{-\infty}^{-1}$. Again this action is denoted by $\alpha_l$.
\end{Not}

\subsection{Noncommutative random walks}
Here we review the basic properties of noncommutative random walks on discrete quantum groups. For classical random walks take a look in \cite{Woess00}.

\begin{Def}
A {\it discrete Markov chain} consists of a pair $(X,P)$ where $X$ is a discrete space and $P=\{p(x,y)\}_{x,y\in X}$ is a matrix which satisfies the properties $p(x,y)\in [0,1]$ for all $x,y\in X$ and $\sum_{y\in X} p(x,y)=1$ for all $x\in X$. We say that $P$ defines a (classical) random walk on $X$. The scalars $p(x,y)$ are the transition probabilities that the random walk jumps from $x$ to $y$.
\end{Def}

\begin{Def}
Suppose that $\{p(x,y)\}_{x,y\in X}$ defines a random walk on $X$. Let
\begin{align*}
p^n(x,y):=\begin{cases} \delta_{x,y} & \textrm{if } n=0;\\
p(x,y), & \textrm{if }n=1;\\
\sum_{z \in X}p^{n-1}(x,z)p(z,y), &\textrm{if } n>1. \end{cases}
\end{align*}
This is the probability that the random walk is at $y$ after $n$ steps when started at $x$. The random walk is called {\it transient} if $\sum_{n=1}^\infty p^n(x,y) <\infty$ for all $x,y\in X$. It is called {\it irreducible} if for all $x,y\in X$ there exists an $n\in\Nat$ such that $p^n(x,y)>0$.
\end{Def}

Transience means that for any pair $(x,y)$ the expected number of times the random walk hits $y$ when starting at $x$ is finite. Irreducibility exactly means that every point can eventually be reached from any other point with positive probability.

\begin{Def}
Let $U$ be a finite dimensional representation of $G$. Define the state
\[
\varphi_U\colon B(\Hcal_U)\ra\Com,\qquad T\mapsto \frac{\Tr(T\pi_U(\rho^{-1}))}{d_U}.
\]
Write $\varphi_s:=\varphi_{U_s}$. Obviously $\varphi_s\circ\pi_s$ defines a state on $l^\infty(\hat G)$, which we again denote by $\varphi_s$. Denote $\Cscr:=\overline{\Span\{\varphi_s\,:\, s\in\Irr(G)\}}$, where we take the norm-closure. Clearly, for any $\varphi\in\Cscr$ there exists a finite (complex) measure $\mu$ on $\Irr(G)$ such that $\varphi=\sum_{s\in\Irr(G)} \mu(s)\varphi_s$. In that case we write $\varphi=\varphi_\mu$. Note that $\varphi_\mu$ is a state if and only if $\mu$ is a probability measure.
\end{Def}

The orthogonality relations imply that
\[
\varphi_s(x)1_s = (h\ten\iota)(U_s^*(1\ten x)U_s),
\]
so $\varphi_s$ is a $\alpha_l|_{B(\Hcal_s)}$-invariant state. In fact, it is the unique $\alpha_l|_{B(\Hcal_s)}$-invariant state (\cite[\textsection~1.4]{NeshveyevTuset04}). We get
\begin{equation}\label{phi_invariant_eq1}
\varphi_s(x) = \varphi_s(x)\varphi_s(1_s) = (h\ten\varphi_s)\big(W^*(1\ten x)W\big).
\end{equation}
Note that \eqref{Haar_weights_eq2} translates to $\hat\psi(x) = \sum_s d_s^2\varphi_s(x)$.

\begin{Def}
Given a normal linear functional $\varphi$ on $l^\infty(\hat G)$, define the {\it Markov operator} associated to $\varphi$ by $P_\varphi:= (\varphi\ten\iota)\hat\Delta$.
\end{Def}

The operator $P_\varphi$ is a completely positive if $\varphi$ is positive. In addition the operator satisfies $P_\varphi(Z(l^\infty(\hat G)))\subset Z(l^\infty(\hat G))$ if and only if $\varphi\in\Cscr$ (see \cite[Prop.\ 2.1]{NeshveyevTuset04}). Therefore we focus on Markov operators defined by states of the form $\varphi_\mu$ for probability measures $\mu$. We write $P_\mu:=P_{\varphi_\mu}$.

In the literature, there is another quite common convention. Namely noncommutative random walks defined by the states $\psi_s:=d_s^{-1}\Tr(\pi_s(\cdot\,\rho))$ and $\tilde P_\mu:=\sum_s \mu(s)(\iota\ten\psi_s)\hat\Delta$. Thus slicing in the right leg of the comultiplication with a different state. Of course all results that hold for $P_\mu$ also hold for $\tilde P_\mu$ and conversely, but one has to be aware on how to translate them. In this paper we will only work with $P_\mu$.

Let $I=\ds_s I_s$ be the identity in $l^\infty(\hat G)$. Since $\hat\Delta(I)=I\ten I$ and $\varphi_\mu$ is a state, it follows that $P_\mu(I)=I$. Define scalars $p_\mu(s,t)\in[0,1]$ by $P_\mu(I_t)I_s=p_\mu(s,t)I_s$, then it follows that $\sum_{t\in\Irr(G)} p_\mu(s,t)=1$ and $p_\mu(s,t)>0$. Therefore $\{p_\mu(s,t)\}_{s,t\in\Irr(G)}$ defines a discrete Markov chain on $\Irr(G)$. Define the measure $\bar{\mu}$ by $\bar{\mu}(s):=\mu(\overline{s})$. Write $\check \varphi_s:=\varphi_{\bar s}$ and extend this linearly to the states $\varphi_\mu$. If $\varphi$ and $\psi$ are two functionals on $l^\infty(\hat G)$, the product is defined by $\varphi\psi:=(\varphi\ten\psi)\hat\Delta$. For properties of this Markov chain see \cite[Lem.\ 2.4, Cor. 2.5]{NeshveyevTuset04}. In particular using induction one easily shows that $p_\mu^n(s,t) = p_{\varphi_\mu^n}(s,t)$.

\begin{Def}\label{Def_quantum_generating_transient}
The measure $\mu$ or the operator $P_\mu$ is called {\it generating} if for all $s\in\Irr(G)$ there exists an $n\in\Nat$ such that $\varphi_\mu^n(I_s)>0$. $P_\mu$ or $\mu$ is called {\it transient} if the random walk with kernel $\{p_\mu(s,t)\}_{s,t\in\Irr(G)}$ is transient.
\end{Def}

Transience almost automatically holds for ``true'' quantum groups as is shown by the following result.

\begin{Lem}[\protect{\cite[Thm.\ 2.6]{NeshveyevTuset04}}]\label{constants_cnr}
Suppose that $\mu$ is a probability measure on $\Irr(G)$. Define constants $c_{n,r}(\mu)$ by the identity $\varphi_\mu^n = \sum_{r\in\Irr(G)} c_{n,r}(\mu)\varphi_r$, then
\begin{equation}\label{constants_cnr_eq1}
c_{n,r}(\mu) = \sum_{t_1,\ldots,t_n\in\Irr(G)} \mu(t_1)\cdots \mu(t_n) m_{t_1,\ldots,t_n}^r \frac{d_r}{d_{t_1}\cdots d_{t_n}}.
\end{equation}
Put $\lambda:= \sum_r \mu(r)\,\frac{\dim(U_r)}{d_r}$. The following inequalities hold for any $n\in\Nat$
\begin{align*}
\sum_{r\in\Irr(G)} c_{n,r}(\mu)d_r^{-1}&\leq \lambda^n; \\
p_\mu^n(s,t)&\leq \frac{d_t}{d_s}\,\frac{\dim(U_s)}{\dim(U_t)}\,\lambda^n,\qquad \textrm{for any } s,t\in\Irr(G).
\end{align*}
In particular if there exists $s\in\supp(\mu)$ such that $\dim(U_s)<d_s$, then
\begin{align*}
\lim_{n\ra\infty} \sum_{r\in\Irr(G)} c_{n,r}(\mu)d_r^{-1} =0; && \sum_{n=1}^\infty p_\mu^n(s,t)<\infty,
\end{align*}
in which case $\mu$ is transient.
\end{Lem}

The statements regarding $p_\mu^n$ are from \cite[Thm.\ 2.6]{NeshveyevTuset04}. The identities of $c_{n,r}$ are an easy consequence of the computations of the proof of that theorem.

\medskip
\begin{pf}[ of Lemma \ref{constants_cnr}]
We prove \eqref{constants_cnr_eq1} by induction on $n$. If $n=1$, this is trivially true as $m_{t_1}^r=\delta_{t_1,r}$. For $n>1$ we have
\begin{align*}
\varphi_\mu^{n+1} &= \varphi_\mu^n\varphi_\mu = \sum_{r,t} c_{n,r}(\mu) \mu(t)\varphi_r\varphi_t\\
&= \sum_{r,t}\Big(\sum_{t_1,\ldots,t_n} \mu(t_1)\cdots \mu(t_n) m_{t_1,\ldots,t_n}^r \frac{d_r}{d_{t_1}\cdots d_{t_n}}\Big)\mu(t)\Big(\sum_s m_{r,t}^s \frac{d_s}{d_rd_t} \varphi_s\Big) \displaybreak[2]\\
&= \sum_s \Big(\sum_{r,t_{n+1}}\sum_{t_1,\ldots,t_n} \mu(t_1)\cdots \mu(t_n) \mu(t_{n+1})  \frac{d_s}{d_{t_1}\cdots d_{t_n}d_{t_{n+1}}} m_{t_1,\ldots,t_n}^r m_{r,t_{n+1}}^s\Big) \varphi_s\\
&= \sum_s \Big(\sum_{t_1,\ldots,t_{n+1}} \mu(t_1)\cdots \mu(t_{n+1})  \frac{d_s}{d_{t_1}\cdots d_{t_{n+1}}} m_{t_1,\ldots,t_n,t_{n+1}}^s\Big) \varphi_s,
\end{align*}
which completes the induction. Use the estimates of Lemma \ref{estimates_multiplicities} to obtain
\begin{align*}
\sum_r c_{n,r}(\mu)d_r^{-1} &= \sum_{t_1,\ldots,t_n} \mu(t_1)\cdots \mu(t_n) \Big(\sum_r m_{t_1,\ldots,t_n}^r\Big) \frac{1}{d_{t_1}\cdots d_{t_n}} \\
&\leq \sum_{t_1,\ldots,t_n} \mu(t_1)\cdots \mu(t_n) \frac{\dim(U_{t_1})\cdots\dim(U_{t_n})}{d_{t_1}\cdots d_{t_n}}\\
&= \Big(\sum_t \mu(t)\,\frac{\dim(U_t)}{d_t}\Big)^n.
\end{align*}
For the estimate of $p_\mu^n$ observe that \eqref{constants_cnr_eq1} implies that
\begin{align*}
p_\mu^n(s,t) &= p_{\varphi_\mu^n}(s,t) = \sum_r c_{n,r}(\mu)\,\frac{d_t}{d_rd_s}m_{r,s}^t\\
&= \sum_{r_1,\ldots,r_n} \mu(r_1)\cdots\mu(r_n)\frac{1}{d_{r_1}\cdots d_{r_n}}\,\frac{d_t}{d_s}\, m_{r_1,\ldots,r_n,s}^t \displaybreak[1]\\
&\leq \sum_{r_1,\ldots,r_n} \mu(r_1)\cdots\mu(r_n)\frac{1}{d_{r_1}\cdots d_{r_n}}\,\frac{d_t}{d_s}\, \frac{\dim(U_{r_1})\cdots\dim(U_{r_n})\dim(U_s)}{\dim(U_t)}  = \frac{d_t}{d_s}\,\frac{\dim(U_s)}{\dim(U_t)}\,\lambda^n.
\end{align*}
The last part follows from the observation that if $d_s>\dim(U_s)$ for some $s\in \supp(\mu)$, then $0<\sum_t \mu(t)\,\frac{\dim(U_t)}{d_t}<1$.
\end{pf}

\begin{Def}[\protect{\cite{NeshveyevTuset04}}]
Suppose that $\mu$ is transient then the following operator, the {\it Green kernel}, makes sense
\[
G_\mu\colon c_{00}(\hat G)\ra l^{\infty}(\hat G),\qquad x\mapsto \sum_{n=0}^\infty P_\mu^n(x).
\]
If in addition $\mu$ is generating, the {\it Martin kernel} of $\mu$ given by
\[
K_{\bar{\mu}}\colon c_{00}(\hat G)\ra l^{\infty}(\hat G), \qquad x\mapsto G_{\bar{\mu}}(x)(G_{\bar{\mu}}(I_0))^{-1}
\]
is well-defined. Here $I_0$ denotes the identity in the trivial representation $\Hcal_0=\Com$. By definition, the {\it Martin compactification} $\tilde M(\hat G,\mu)$ of $\hat G$ with respect to $\mu$ is the \Cstar subalgebra of $l^\infty(\hat G)$ generated by $c_0(\hat G)$ and $K_{\bar{\mu}}\big(c_{00}(\hat G)\big)$. The {\it Martin boundary} $M(\hat G,\mu)$ of $\hat G$ is the quotient \Cstar algebra $\tilde M(\hat G,\mu) / c_0(\hat G)$.
\end{Def}

\begin{Def}[\protect{\cite[\textsection 2.5]{Izumi02}}]\label{Def_Poisson_boundary}
Let $\varphi$ be a state on $l^\infty(\hat G)$. An element $h\in l^\infty(\hat G)$ is called {\it $\varphi$-harmonic} if $P_\varphi(h)=h$. If $x\geq 0$ and $P_\varphi(x)\leq x$, then $x$ is called {\it $\varphi$-superharmonic}. Define the {\it Poisson boundary} $H^\infty(\hat G,\varphi):=\{h\in l^\infty(\hat G)\,:\, P_\varphi(h)=h\}$. This is a von Neumann algebra with product $x\cdot y:= \sstar\lim_n P_\varphi^n(xy)$. If $\varphi=\varphi_\mu$ for some $\mu$ probability measure on $\Irr(G)$ we write $H^\infty(\hat G,\mu)$ for the Poisson boundary.
\end{Def}

Note that in the literature different manifestations of this product can be found. However, all products are the same, because $H^\infty(l^\infty(\hat G),\mu)$ is an operator system in $l^\infty(\hat G)$. By a result of Choi and Effros \cite{ChoiEffros77}, it admits at least one product turning it into a von Neumann algebra. On the other hand a \Cstar algebra has, up to complete order isomorphism, a unique product. The nontrivial part is therefore to find this product. Originally Izumi \cite[\textsection 2.5]{Izumi02} used an ultra-filter and Ces\`{a}ro-summation to define the product (see also \cite[Appendix]{Izumi12} for a discussion on the product on $H^\infty(\hat G,\varphi)$).

\begin{Prop}[\protect{\cite[Thm.\ 3.3]{NeshveyevTuset04}}]
For any $\mu$-superharmonic element $x\in l^\infty(\hat G)$ there exists a bounded positive linear functional $\omega_x\colon \tilde M(\hat G,\mu)\ra\Com$ such that $(y,x)_{\hat\psi}=\omega_x(K_{\bar\mu}(y))$ for all $y\in c_{00}(\hat G)$.

Conversely if $\omega\colon \tilde M(\hat G,\mu)\ra\Com$ is a bounded positive linear functional, then there exists a unique superharmonic element $x_\omega\in l^\infty(\hat G)$ such that $(y,x_\omega)_{\hat\psi}=\omega (K_{\bar\mu}(y))$ for all $y\in c_{00}(\hat G)$. If $x_\omega$ is $\mu$-harmonic, then $\omega|_{c_0(\hat G)}=0$. Moreover, if $\supp(\mu)$ is finite, then $x_\omega$ is harmonic if and only if $\omega|_{c_0(\hat G)}=0$.
\end{Prop}

We say that the above functional $\omega_x$ represents the element $x$.

\subsection{General results on convergence to the boundary}\label{subsec_general_results_boundary_convergence}
Discrete Markov chains are known to converge to the boundary. For random walks on discrete quantum groups this is still an open problem. In this section we prove some results regarding the convergence of paths in the quantum setting.

\bigskip
Recall that a net $(x_i)_i\subset B(\Hcal)$ is said to converge to $x$ in the {\it strong topology} if $\|(x_i-x)\xi\|\ra 0$ for all $\xi\in \Hcal$. It converges in {\it strong$^*$ topology} if both $x_i\ra x$ and $x_i^*\ra x^*$ in strong topology. We write $\textrm{s-}\lim_n x_i=x$ and $\sstar\lim_n x_i=x$ for the strong and strong$^*$ limits.

If $M\subset B(\Hcal)$ is a von Neumann algebra, then the strong and strong$^*$ topologies on $M$ in general depend on the embedding of $M$ into $B(\Hcal)$. The restriction of the topologies to bounded sets are independent of such an embedding.
Note that if $(x_i)_i\subset M$ is a bounded net and $\varphi\colon M\ra\Com$ is a faithful state on $M$, then $\textrm{s-}\lim_i x_i= x$ if and only if $\lim_i\varphi((x_i-x)^*(x_i-x))=0$.

\begin{Not}
Let $\varphi$ be a state on $(l^\infty(\hat G),\hat \Delta)$. Denote for $n\geq 1$ the unital $*$-homomorphisms
\[
j_n\colon l^\infty(\hat G)\ra\bigotimes_{-\infty}^{-1} (l^\infty(\hat G),\varphi), \qquad x\mapsto \cdots\ten 1\ten1\ten\hat\Delta^{n-1}(x).
\]
An element $x$ is called {\it $\varphi$-regular} if $\textrm{s}^*\textrm{-}\lim_n j_n(x)$ exists in $\bigotimes_{-\infty}^{-1} (l^\infty(\hat G),\varphi)$. Denote $R_\varphi:=\{x\in l^\infty(\hat G)\,:\, x \textrm{ is $\varphi$-regular}\}$. For $x\in R_\varphi$ define $j_\infty(x):= \textrm{s}^*\textrm{-}\lim_n j_n(x)$. Write $\varphi^\infty:=\cdots\ten\varphi\ten\varphi$ and observe that $\varphi^\infty\circ j_n=\varphi^n$.
\end{Not}

This algebra $\bigotimes_{-\infty}^{-1} (l^\infty(\hat G),\varphi)$ turns out to be the proper generalisation of the space of functions on paths of the random walk. The maps $j_n$ correspond to the duals of the coordinate maps.

\begin{Prop}[\protect{Izumi \cite{Izumi02} and Neshveyev--Tuset \cite{NeshveyevTuset03}}] \label{prop_reg_algebra}
Let $(l^\infty(\hat G),\hat \Delta)$ be a discrete quantum group and $\mu$ be a probability measure on $\Irr(G)$. Assume that the random walk defined by $\varphi=\varphi_\mu$ is generating and transient. The following holds:
\begin{enumeraterm}
\item $R_\varphi$ is a \Cstar subalgebra of $l^\infty(\hat G)$ that contains $c_0(\hat G)$ and $H^\infty(\hat G,\varphi)$.
\end{enumeraterm}
Write $\theta:=j_\infty|_{H^\infty(\hat G,\varphi)}\colon H^\infty(\hat G,\varphi)\ra \bigotimes_{-\infty}^{-1} (l^\infty(\hat G),\varphi)$. Then
\begin{enumeraterm}[resume]
\item $j_\infty\colon R_\varphi\ra \bigotimes_{-\infty}^{-1} (l^\infty(\hat G),\varphi)$ is a $*$-homomorphism onto $\theta(H^\infty(\hat G,\varphi))$;
\item $c_0(\hat G)\subset \ker(j_\infty)$;
\item $\theta\colon H^\infty(\hat G,\varphi) \ra \theta(H^\infty(\hat G,\varphi))$ is a $*$-isomorphism;
\item $\theta^{-1}\circ j_\infty(x) = \sstar\lim_n P_\varphi^n(x)$, for any $x\in R_\varphi$. In particular $\hat\eps\circ\theta^{-1}\circ j_\infty = \lim_n \varphi^n$.
\end{enumeraterm}
Denote $\theta_0:=\theta^{-1}\circ j_\infty\colon R_\varphi\ra H^\infty(\hat G,\varphi)$. This is a $*$-homomorphism by (ii) and (iv).
\end{Prop}

No proof of this result is given in \cite{NeshveyevTuset03} and only some parts are in \cite[Thm.\ 3.6]{Izumi02}, so we give a proof for completeness. Note that part (iv) does not follow from part (ii) since the Poisson boundary is equipped with different product than $l^\infty(\hat G)$ (cf.\ Definition \ref{Def_Poisson_boundary}).

\medskip
\begin{pf}[ of Proposition \ref{prop_reg_algebra}]
(i) Clearly $R_\varphi$ is a linear space which is closed under the involution $*$. Each $j_n$ is a $*$-homomorphism, so it is norm-decreasing, i.e.\ $\|j_n(x)\|\leq \|x\|$. Therefore if $x,y\in R_\varphi$ and $\xi\in\Hcal_{\varphi^\infty}$
\begin{align*}
\|(j_n(xy)-&j_\infty(x)j_\infty(y))\xi\| = \|(j_n(x)j_n(y)-j_\infty(x)j_\infty(y))\xi\|\\
&\leq \|(j_n(x)j_n(y)-j_n(x)j_\infty(y))\xi\| + \|(j_n(x)j_\infty(y)-j_\infty(x)j_\infty(y))\xi\|\\
&\leq \sup_m \|j_m(x)\| \|(j_n(y)-j_\infty(y))\xi\| + \|(j_n(x)-j_\infty(x))j_\infty(y)\xi\|,
\end{align*}
which tends to $0$ as $n\ra\infty$, thus $xy\in R_\varphi$ and $j_\infty(xy)=j_\infty(x)j_\infty(y)$. To show that $R_\varphi$ is closed in norm, let $(x_n)_n\subset R_\varphi$ be a bounded sequence with $\lim_n\|x_n-x\|= 0$. Given $\eps>0$ and $\xi\in\Hcal_{\varphi^\infty}$, find $m$ such that $\|x-x_m\|<\eps$. Then there exists $N\in\Nat$ such that for all $n,n'\geq N$ it holds $\|(j_n(x_m)-j_{n'}(x_m))\xi\|\leq \eps$. This gives for $n,n'\geq N$
\begin{align}
\|(j_n(x)-j_{n'}(x))\xi\| &\leq \|(j_n(x)-j_{n}(x_m))\xi\| + \|(j_n(x_m)-j_{n'}(x_m))\xi\| + \|(j_{n'}(x_m)-j_{n'}(x))\xi\| \notag\\
&\leq \|x-x_m\|\|\xi\| + \eps + \|x_m-x\|\|\xi\| \leq \eps+2\|\xi\|\eps. \label{prop_reg_algebra_eq4}
\end{align}
Thus $R_\varphi$ is closed in norm and hence a \Cstar algebra.\\
We claim that $\sstar\lim_nj_n(x)=0$ for every $x\in c_0(\hat G)$. This would show that $c_0(\hat G)\subset R_\varphi$ and $c_0(\hat G)\subset\ker(j_\infty)$. First consider $I_s\in B(\Hcal_s)$, then
\begin{align*}
\varphi_\mu^\infty(j_n(I_s)) &= (\varphi_\mu\ten\cdots\ten\varphi_\mu)(\hat\Delta^{n-1}(I_s)) =(\varphi_\mu\ten\cdots\ten\varphi_\mu\ten\hat\eps)(\hat\Delta^n(I_s))\\
&=\hat\eps(P_\mu^n(I_s)) = \sum_{t\in\Irr(G)} \hat\eps(P_\mu^n(I_s)I_t) = \sum_{t\in\Irr(G)} \hat\eps(I_t) p_\mu^n(t,s)= p_\mu^n(0,s).
\end{align*}
By assumption $P_\mu$ is transient, thus $\sum_n p_\mu^n(0,s) = g_\mu(0,s)<\infty$ and hence $\lim_n p_\mu^n(0,s)=0$. If $x\in c_{00}(\hat G)$ is finitely supported, then $x^*x$ is dominated by a linear combination $\sum_{i=1}^m c_iI_{s_i}$ for some $s_i\in\Irr(G)$ and $c_i>0$. Therefore
\[
\lim_n\varphi^\infty\big(j_n(x)^*j_n(x)\big)\leq\sum_i c_i\lim_n p_\mu^n(0,s_i)=0.
\]
Since $R_\varphi$ is a \Cstar algebra and $c_{00}(\hat G)$ is dense in $c_0(\hat G)$ the claim follows.\\
The statement that $H^\infty(\hat G,\varphi)\subset R_\varphi$ is contained in \cite[Thm.\ 3.6]{Izumi02}. Izumi only deals with finitely supported measures $\mu$, but the proof also applies to our case. Recall the conditional expectations $E_n$ and maps $E_n'$ introduced in Lemma \ref{lem_martingale_convergence}. By coassociativity one has for $m>n$ and any $x\in l^\infty(\hat G)$
\begin{align*}
E'_n(j_m(x)) &= (\cdots\ten\varphi\ten\varphi\ten\iota^{\ten n})(\cdots\ten1\ten1\ten\hat{\Delta}^{m-1}(x))\\
&= (\underbrace{\varphi\ten\cdots\ten\varphi}_{m-n}\ten\iota^{\ten n})((\iota^{\ten m-n}\ten\hat{\Delta}^{n-1})\hat{\Delta}^{m-n}(x)) = \hat{\Delta}^{n-1}(P_\varphi^{m-n}(x))
\end{align*}
and thus
\begin{equation}\label{prop_reg_algebra_eq2}
E_n(j_m(x))=j_n(P_\varphi^{m-n}(x)).
\end{equation}
Therefore if $h\in H^\infty(\hat G,\varphi)$ is a harmonic element, $E_n(j_{n+1}(h)) = j_n(h)$. In which case the noncommutative Martingale convergence theorem (cf.\ Lemma \ref{lem_martingale_convergence}) implies that there exists a unique element $y\in \bigotimes_{-\infty}^{-1} (l^\infty(\hat G),\varphi)$ such that $E_n(y)=j_n(h)$ for all $n$ and that the sequence $(E_n(y))_n$ converges to $y$ in strong$^*$ topology. In other words $\sstar\lim_n j_n(h)=y$, thus $H^\infty(\hat G,\varphi)\subset R_\varphi$, concluding part (i).

We now deal with (iv). Clearly $\theta$ is surjective on its image, linear and preserves the involution. For an element $x\in R_\varphi$ the limit $\sstar\lim_n P_\varphi^n(x)$ exists by \eqref{prop_reg_algebra_eq2} and is harmonic. Therefore by part (i) it is again in $R_\varphi$. Thus if $h_1,h_2\in H^\infty(\hat G,\varphi)$ are harmonic, by the noncommutative martingale convergence theorem and \eqref{prop_reg_algebra_eq2}
\begin{align*}
j_\infty(h_1)j_\infty(h_2) &= j_\infty(h_1h_2)\\
&= \sstar\lim_n E_n j_\infty(h_1h_2)\\
&= \sstar\lim_n \sstar\lim_m E_n j_m(h_1h_2)\\
&= \sstar\lim_n \sstar\lim_m j_n(P_\varphi^{m-n}(h_1h_2))\\
&= \sstar\lim_n j_n(h_1\cdot h_2)\\
&= j_\infty(h_1\cdot h_2),
\end{align*}
thus $j_\infty\colon H^\infty(\hat G,\varphi)\ra \bigotimes_{-\infty}^{-1} (l^\infty(\hat G),\varphi)$ is a $*$-homomorphism. It remains to show that $j_\infty$ is isometric. Since $j_\infty$ is a $*$-homomorphism it is norm-decreasing, so $\|h\|\geq \|j_\infty(h)\|$. It thus suffices to prove the reverse inequality. Denote $U:=\bigoplus_{s\in\supp{\mu}} U_s$, so if $x\in \bigotimes_{-n}^{-1}l^\infty(\hat G)$ then
\[
\pi_U^{\ten n}(x)I_{s_1\ten\cdots\ten s_n} = \begin{cases}
x I_{s_1\ten\cdots\ten s_n} &\textrm{if } s_1,\ldots,s_n\in\supp(\mu);\\
0&\textrm{otherwise}.
\end{cases}
\]
Since the measure $\mu$ is assumed to be generating we obtain
\begin{align*}
\|h\| &= \sup\{\|\pi_s(h)\|\,:\, s\in\Irr(G)\}\\
&= \sup\{\|\pi_s(h)\|\,:\, s\in\supp(\mu^{*n}),\, n\in\Nat\} \displaybreak[2]\\
&= \sup\{\|\pi_{U}^{\ten n}(\hat\Delta^n(h))\|\,:\, n\in\Nat\}\\
&= \sup\{\|j_n(h)\|\,:\, n\in\Nat\}.
\end{align*}
From \eqref{prop_reg_algebra_eq2} we conclude $j_n(h)=E_n(j_\infty(h))$. Since $E_n$ is a conditional expectation we get
\[
\|h\|=\sup_{n\in\Nat}\|j_n(h)\| = \sup_{n\in\Nat}\|E_n(j_\infty(h))\| \leq \|j_\infty(h)\|.
\]
Hence $j_\infty$ is isometric and thus a $*$-isomorphism onto its image, which proves (iv).

We prove (ii) and (v) simultaneously. Clearly $j_\infty$ is a $*$-homomorphism. By (i) it is immediate that $\theta(H^\infty(\hat G,\varphi))\subset j_\infty(R_\varphi)$.
Let $x\in R_\varphi$. By the noncommutative martingale convergence theorem and \eqref{prop_reg_algebra_eq2} we get
\begin{align}
j_\infty(x) &= \sstar\lim_n E_n j_\infty(x) \notag\\
&= \sstar\lim_n \sstar\lim_m E_n j_m(x) \notag \displaybreak[2]\\
&= \sstar\lim_n \sstar\lim_m j_n(P_\varphi^{m-n}(x)) \notag \displaybreak[2]\\
&= \sstar\lim_n \sstar\lim_m j_n(P_\varphi^{m}(x)) \notag\\
&= j_\infty(\sstar\lim_m P_\varphi^m(x)). \label{prop_reg_algebra_eq3}
\end{align}
Since $\sstar\lim_n P_\varphi^n(x)$ is harmonic, we obtain $j_\infty(R_\varphi)\subset j_\infty(H^\infty(\hat G,\varphi))$ which establishes (ii). Moreover \eqref{prop_reg_algebra_eq3} can be written as $\theta^{-1}\circ j_\infty(x) = \sstar\lim_m P_\varphi^m(x)$. This gives
\[
\hat\eps \circ \theta^{-1}\circ j_\infty(x) = \lim_m \hat\eps\circ P_\varphi^m(x) = \lim_m\hat\eps(\varphi^m\ten\iota)\hat\Delta(x) = \lim_m\varphi^m(x),
\]
which finishes the proof.
\end{pf}

\begin{Conj}[\protect{\cite{NeshveyevTuset03}}]\label{conjecture_convergence}
Let $\mu$ be a transient and generating probability measure on $\Irr(G)$. The following holds for the random walk defined by $\varphi=\varphi_\mu$ on $(\hat G,\hat\Delta)$:
\begin{enumeraterm}
\item $K_{\check\varphi}(x)\in R_\varphi$, for every $x\in c_{00}(\hat{G})$;
\item If $\nu=\lim_n \varphi^n|_{R_\varphi}=\hat\epsilon\theta_0$, then $\hat\psi(xh) = \nu(K_{\check\varphi}(x)h)$ for every $x\in c_{00}(\hat{G})$ and any harmonic element $h\in H^\infty(\hat{G}, \varphi)$.
\end{enumeraterm}
\end{Conj}

There is no proof known for this conjecture. We therefore say that a random walk on a discrete quantum group $(\hat G,\hat\Delta)$ defined by a probability measure $\mu$ {\it converges to the boundary} if the statements of the conjecture hold for $(\hat G,\mu)$.
This notion of convergence to the boundary is compatible with the classical case, see Proposition \ref{Prop_correspondence_classical_convergence} below.

Consider the modular group $\{\sigma^{\hat\psi}_t\}_t$ of the right Haar weight (see \eqref{Haar_weights_eq3}). It is easy to show that (see \cite[\textsection~3.3]{NeshveyevTuset04}) the following holds:
\begin{align}
\sigma_t^{\hat\psi} P_{\check\varphi} = P_{\check\varphi} \sigma_t^{\hat\psi}, && \sigma_t^{\hat\psi} K_{\check\varphi} = K_{\check\varphi} \sigma_t^{\hat\psi}. \label{modular_group_eq1}
\end{align}
Moreover, if $x\in l^\infty(\hat G)$, then $\|\pi_s(\sigma_t^{\hat\psi}(x))\|\leq\|\pi_s(x)\|$ whenever $t\in\Rea$ and thus $\sigma_t^{\hat\psi}(c_0(\hat G))\subset c_0(\hat G)$. It follows that $\sigma_t^{\hat\psi}$ defines a automorphism group on the Martin compactification and factors through the Martin boundary.

The state $\nu|_{\tilde M(\hat G,\varphi)}$ is a $\sigma^{\hat\psi}$-KMS state representing the unit $1\in H^\infty(\hat G,\varphi)$ (see \cite[Thm.\ 3.10]{NeshveyevTuset04}). So in particular $\sigma^\nu_t=\sigma^{\hat\psi}_t$ for all $t$. Define the map $K_{\check\varphi}^*\colon\pi_\nu(\tilde M(\hat G,\varphi))''\ra l^\infty(\hat G)$ by the identity
\begin{equation}\label{def_Kstar}
(a,K_{\check\varphi}(x))_\nu = (K_{\check\varphi}^*(a),x)_{\hat\psi}, \qquad (a\in\pi_\nu(\tilde M(\hat G,\varphi))'',\, x\in c_{00}(\hat G)),
\end{equation}
where $(\cdot,\cdot)_\nu$ and $(\cdot,\cdot)_{\hat\psi}$ are as in Lemma \ref{lem_inner_prod}. Note that $\hat\psi$ is unbounded on $l^\infty(\hat G)$, but it is well defined on $c_{00}(\hat G)$. Identity \eqref{def_Kstar} determines $K_{\check\varphi}^*$ uniquely because we can take for $x$ all matrix units $m_{ij}^s$. In \cite[Lem.\ 3.9]{NeshveyevTuset04} it is shown that  $\im(K_{\check\varphi}^*)\subset H^\infty(\hat G,\varphi)$ and $K_{\check\varphi}^*$ is normal, unital and completely positive.

\begin{Thm}[\protect{\cite[Thm.\ 6.2]{NeshveyevTuset03}}]\label{Thm_representation_harmonic_elements} If Conjecture \ref{conjecture_convergence} holds for a state $\varphi$ on $(l^\infty(\hat G),\hat \Delta)$, then this random walk has the following properties:
\begin{enumeraterm}
\item for any positive harmonic element $h\in H^\infty(\hat G,\varphi)$ the positive linear functional $(h,\cdot)_\nu$ represents $h$, meaning $(h,x)_{\hat\psi} = (h,K_{\check\varphi}(x))_\nu$ for all $x\in c_{00}(\hat G)$;
\item the map $K_{\check\varphi}^*|_{M(\hat G,\varphi)}\colon M(\hat G,\varphi)\ra H^\infty(\hat G,\varphi)$ equals the map $\theta_0|_{M(\hat G, \varphi)}$. It induces an isomorphism $\pi_\nu(M(\hat G,\varphi))''\cong H^\infty(\hat G,\varphi)$ which respects the action of the compact quantum group $G$ and the dual discrete quantum group $\hat G$.
\end{enumeraterm}
\end{Thm}

Given a discrete group $\Gamma$, consider $C(G):=C_{\textrm{r}}^*(\Gamma)$ and define $\Delta\colon C(G)\ra C(G)\ten C(G)$ by $\Delta(\lambda_\gamma):= \lambda_\gamma\ten \lambda_\gamma$, then $G=(C_{\textrm{r}}^*(\Gamma),\Delta)$ is a compact quantum group. The elements $\{\lambda_\gamma\}_{\gamma\in\Gamma}$ are one-dimensional representations of the compact quantum group $G$. They span a dense subset of $C(G)$. Hence by the orthogonality relations these are all finite dimensional irreducible representations. Thus $\Irr(G)\cong\Gamma$ and $\Com[G]\subset C(G)= C_{\textrm{r}}^*(\Gamma)$ equals the ordinary group algebra of $\Gamma$. It follows that $\Ucal(G)=\Com[G]^*$ equals the functions on $\Gamma$. The discrete quantum group can be identified as $l^\infty(\hat G)=l^\infty(\Gamma)$ with pointwise multiplication and involution. The comultiplication is given by $\hat\Delta(f)(\gamma,\gamma'):=f(\gamma\gamma')\in l^\infty(\Gamma\times\Gamma)\cong l^\infty(\Gamma)\bar\ten l^\infty(\Gamma)$.\\
Let $\mu$ be a probability measure on $\Gamma=\Irr(G)$. Since all irreducible representations of $G$ are one-dimensional, $\rho$ equals the identity and thus $\varphi_\mu\colon l^\infty(\Gamma)\ra\Com$ equals $\varphi_\mu(f):=\sum_{\gamma\in\Gamma} f(\gamma)\mu(\gamma)$. It follows that
\[
P_\mu(f)(x) = \big((\varphi\ten\iota)\hat\Delta(f)\big)(x) =\sum_{y\in\Gamma}\mu(y)f(y^{-1}x)= \sum_{y\in\Gamma}\mu(xy^{-1})f(y).
\]
So this equals the classical random walk on $\Gamma$ with Markov kernel $p_\mu(x,y)=\mu(xy^{-1})$. With this notation we get.

\begin{Prop}\label{Prop_correspondence_classical_convergence}
Suppose that we are given a probability measure $\mu$ on a discrete group $\Gamma$. Assume that $\mu$ defines a generating and transient random walk on $\Gamma$. Let $P_\mu$ be the Markov kernel of the random walk on $l^\infty(\hat G)$. Then this quantum random walk converges to the boundary in the sense of Conjecture \ref{conjecture_convergence}.
\end{Prop}

We will not prove this proposition, because it would require to build up the whole theory of classical random walks and convergence to the boundary thereof. Instead we give a rough sketch.

\medskip
\begin{sketch}[ of Proposition \ref{Prop_correspondence_classical_convergence}]
First identify the quantum Martin compactification $\tilde M(\hat G,\mu)$ with the functions on the classical Martin compactification $\tilde M(\Gamma,\mu)$ via
\begin{align*}
&\kappa\colon C(\tilde M(\Gamma,\mu))\ra l^\infty(\hat G);\\
&1_\gamma\ra I_\gamma, \qquad k_\mu(t,\cdot)\mapsto d_t^{-2}K_{\bar\mu}(I_t).
\end{align*}
This fully describes an identification, because by construction of the Martin compactification the functions $\{1_\gamma\}_{\gamma\in\Gamma}$ and $\{k_\mu(t,\cdot)\}_{t\in\Gamma}$ generate the algebra of continuous functions on the Martin compactification.
The next step is to find out how this identification behaves with respect to the maps $j_n$. Let $\Omega=\{(\omega_n)_{n=0}^\infty\,:\, \omega_n\in\Gamma\}$ be the path space and $X_n\colon\Omega\ra\Gamma$ be the $n$th coordinate projection. It can be shown that $j_n$ corresponds to $X_n$.
Classical convergence to the boundary says that the sequence $X_n$ converges almost everywhere (with respect to some measure on $\Omega$ constructed out of $\mu$) to a Borel measurable function $X_\infty\colon \Omega\ra M(\Gamma,\mu)$ \cite[Thm.\ 24.10]{Woess00}. Using that convergence almost everywhere is stronger than convergence in mean, the identification shows that for any $x\in \tilde M(\hat G,\mu)$ the sequence $(j_n(x))_n$ converges in strong$^*$-topology.

To establish the representation of harmonic elements one first needs to identify the state representing the unit, for this the classical convergence to the boundary is again necessary. If $y\in H^\infty(\hat G,\mu)$ is a general harmonic element, then it corresponds to some harmonic function $h\in H^\infty(\Gamma,\mu)$. Using the classical convergence to the boundary, $h$ can be represented as an integral over the Martin boundary. The identification $\kappa$ then leads to the desired expression for $y$.
\end{sketch}

Convergence to the boundary describes the asymptotic behaviour of the paths in the random walk, but by Theorem \ref{Thm_representation_harmonic_elements} it also gives a natural representation of harmonic elements. Let us take a closer look at this second part of the conjecture.

Recall that $\varphi^\infty$ defines a faithful state on the von Neumann algebra $\bigotimes_{-\infty}^{-1} (l^\infty(\hat G),\varphi)$. Consider the GNS representation $(\pi_{\varphi^\infty},\Hcal_{\varphi^\infty},\xi_{\varphi^\infty})$ of $\theta(H^\infty(\hat G,\varphi))$ defined by $\varphi^\infty$. Since $j_\infty\colon R_\varphi \ra \theta(H^\infty(\hat G,\varphi))$ is surjective and
\[
\nu = \lim_n \varphi^n=\lim_n (\varphi\ten\cdots\ten\varphi)\hat\Delta^{n-1} = \varphi^\infty\circ j_\infty
\]
it follows that $\pi_\nu(R_\varphi)=\pi_{\varphi^\infty}(\theta(H^\infty(\hat G,\varphi)))$ inside $B(\Hcal_\nu)=B(\Hcal_{\varphi^\infty})$. Therefore the morphism $j_\infty\colon R_\varphi\ra\bigotimes_{-\infty}^{-1} (l^\infty(\hat G),\varphi)$ can be extended to a s$^*$-continuous map $j_\infty\colon\pi_{\nu}(R_\varphi)''\ra \theta(H^\infty(\hat G,\varphi))$. We denote this extension by $j_\infty$ and similarly we write $\theta_0:=\theta^{-1}\circ j_\infty\colon\pi_{\nu}(R_\varphi)''\ra H^\infty(\hat G,\varphi)$.

\begin{Lem}\label{lemma_properties_j-theta}
For $x,\,y\in R_\varphi$ the following holds:
\begin{align*}
\theta_0(\theta_0(x))= \theta_0(x) \quad \textrm{and} \quad \theta_0(xy) = \theta_0(x)\cdot\theta_0(y),
\end{align*}
where $\cdot$ is the product in $H^\infty(\hat G,\varphi)$. If in addition $\tilde M(\hat G,\varphi)\subset R_\varphi$, then
\begin{equation}
\nu(K_{\check\varphi}(x)\theta_0(a)) = \nu(K_{\check\varphi}(x)a) \qquad \textrm{for } a\in\pi_\nu(\tilde M(\hat G,\varphi))'' \textrm{ and } x\in c_{00}(\hat G).
\end{equation}
\end{Lem}
\begin{pf}
Let $x$ be $\varphi$-regular. By Proposition \ref{prop_reg_algebra} it follows that $\theta_0(x)\in H^\infty(\hat G,\varphi)$, thus
\[
P_\varphi^m(\theta_0(x))= \theta_0(x), \qquad \textrm{(for every $m\geq 0$)}.
\]
We also know from Proposition \ref{prop_reg_algebra} that $\theta_0 = \sstar\lim_n P_\varphi^n$. Hence
\[
\theta_0\circ \theta_0(x) = \sstar\lim_n P_\varphi^n (\theta_0(x)) = \sstar\lim_n \theta_0(x) = \theta_0(x),
\]
which proves the first identity. The second one follows from Proposition \ref{prop_reg_algebra} part (iv).\\
To prove the last statement, observe that if $a\in\pi_\nu(\tilde M(\hat G,\varphi))''$, then $\theta_0(a)$ is $P_\varphi$-harmonic and $\nu=\hat\eps\circ \theta^{-1}\circ j_\infty=\eps\circ\theta_0$. Hence
\begin{align}
\nu(K_{\check\varphi}(x)\theta_0(a)) &= \hat\epsilon\circ \theta_0(K_{\check\varphi}(x)\theta_0(a)) \notag\\
&= \hat\epsilon\big(\theta_0(K_{\check\varphi}(x))\cdot (\theta_0(\theta_0(a)))\big) \notag \displaybreak[2]\\
&= \hat\epsilon\big(\theta_0(K_{\check\varphi}(x))\cdot (\theta_0(a))\big) \notag\\
&= \hat\epsilon\circ \theta_0(K_{\check\varphi}(x)a)) \notag\\
&= \nu(K_{\check\varphi}(x)a), \label{sufficient_conditions_conj2_eq1}
\end{align}
which proves the result.
\end{pf}

\begin{Prop}\label{sufficient_conditions_conj2} Assume that $\tilde M(\hat G,\varphi)\subset R_\varphi$, so that the random walk on $\hat G$ defined by $\mu$ satisfies the first part of the conjecture. Then the following two conditions are equivalent:
\begin{enumeraterm}
\item $K_{\check\varphi}^*= \theta_0|_{\tilde M(\hat G,\varphi)}$ and $j_\infty\colon\pi_\nu(\tilde M(\hat G,\varphi))'' \ra \theta(H^\infty(\hat G,\varphi))$ is surjective;
\item part (ii) of Conjecture \ref{conjecture_convergence} holds.
\end{enumeraterm}
\end{Prop}
\begin{pf}
Suppose that (i) holds. Let $h\in H^\infty(\hat G,\varphi)$. By assumption there exists an $a\in \pi_\nu(\tilde M(\hat G,\varphi))''$ such that $K_{\check\varphi}^*(a)=\theta^{-1}(j_\infty(a))=h$. Lemma \ref{lemma_properties_j-theta} and the definition of $K_{\check\varphi}^*$ in terms of the inner products \eqref{modular_group_eq1} give
\begin{align*}
\nu(K_{\check\varphi}(x)h) &= \nu(K_{\check\varphi}(x)\theta_0(a))= \nu(K_{\check\varphi}(x)a) \\
&= \nu\big(K_{\check\varphi}(\sigma_{-\frac{i}{2}}^\nu (x))\sigma_{-\frac{i}{2}}^\nu(a^{**})\big)\\
&= \big(a^*,K_{\check\varphi}(\sigma_{-\frac{i}{2}}^\nu(x))\big)_\nu \displaybreak[2]\\
&= \big(K_{\check\varphi}^*(a^*),\sigma_{-\frac{i}{2}}^\nu(x)\big)_{\hat\psi} \displaybreak[2]\\
&= \hat\psi \big(\sigma_{-\frac{i}{2}}^{\hat\psi}(x) \sigma_{-\frac{i}{2}}^{\hat\psi}(K_{\check\varphi}^*(a^*)^*)\big) \\
&= \hat\psi (x K_{\check\varphi}^*(a))= \hat\psi (xh),
\end{align*}
which is exactly the second condition.

To prove the converse implication, note that Theorem \ref{Thm_representation_harmonic_elements} implies that $j_\infty = \theta\circ K_{\check\varphi}^*$ and that $K_{\check\varphi}^*$ is an isomorphism. Hence $j_\infty$ is surjective.
\end{pf}

\begin{Cor}\label{sufficient_conditions2_conj2}
Assume that $\tilde M(\hat G,\varphi)\subset R_\varphi$. In addition, assume that $K_{\check\varphi}\colon c_{00}(\hat G)\ra M(\hat G,\varphi)$ has dense range and $K_{\check\varphi}^*\colon\pi_\nu(M(\hat G,\varphi))''\ra H^\infty(\hat G,\varphi)$ intertwines the right coactions of $(c_0(\hat G),\hat\Delta)$, then part (ii) of Conjecture \ref{conjecture_convergence} holds.
\end{Cor}
\begin{pf}
From \cite[Prop.\ 3.12]{NeshveyevTuset04} it follows that $K_{\check\varphi}^*\colon\pi_\nu(M(\hat G,\varphi))''\ra H^\infty(\hat G,\varphi)$ is an isomorphism, while \cite[Prop.\ 3.11]{NeshveyevTuset04} implies that $K_{\check\varphi}^*=(\nu\ten\iota)\hat\Delta = \theta_0$ which is thus is surjective. Now apply Proposition \ref{sufficient_conditions_conj2}.
\end{pf}

\section{Convergence to the boundary for \texorpdfstring{$\textrm{SU}_{\lowercase{q}}(2)$}{SUq(2)}}
\label{sec_convergence_SUQt}
In Proposition \ref{Prop_correspondence_classical_convergence} we obtained that if $\Gamma$ is a discrete group, a quantum random walk on $l^\infty(\Gamma)$ converges to the boundary in sense of Conjecture \ref{conjecture_convergence}. To obtain a ``true'' quantum example of convergence to the boundary we consider random walks on $\widehat\SUqt$. We shall show that if a probability measure $\mu$ on $\Irr(\SUqt)$ is nice enough, then the random walk on $\widehat\SUqt$ defined by $\mu$ converges to the boundary.

\subsection{\texorpdfstring{$\textrm{SU}_{\lowercase{q}}(2)$}{SUq(2)} and its Martin boundary}
\label{subsec_Martin_SUq2}
Fix $q\in(0,1)$. The \Cstar algebra $C(\SUqt)$  is the universal unital \Cstar algebra generated by $\alpha$ and $\gamma$ such that the matrix
\[
\begin{pmatrix}
\alpha^* & \gamma\\
-q\gamma^* &\alpha
\end{pmatrix}
\]
is unitary. The comultiplication $\Delta$ is defined as
\begin{align*}
\Delta(\alpha):= \alpha\ten\alpha- q\gamma^*\ten\gamma, && \Delta(\gamma):=\gamma\ten\alpha+\alpha^*\ten\gamma.
\end{align*}
The Hopf algebra of matrix coefficients $\Com[\SUqt]$ is the $*$-algebra generated by $\alpha$ and $\gamma$. The irreducible representations can be identified as $\Irr(\SUqt)\cong \half\Int_+$ and the discrete dual equals
\[
l^\infty(\widehat\SUqt) = \bds_{s\in\half\Int_+} B(\Hcal_s),
\]
where $\Hcal_s$ is a $2s+1$-dimensional Hilbert space. This is enough information to be able to state the main theorem of this section.

\begin{Thm}\label{Thm_convergence_SUq2}
Assume that $\mu$ is a probability measure on $\half\Int_+$ that satisfies the following two conditions:
\begin{enumeraterm}
\item there exists $s\in\half\Int_+\setminus\Int_+$ with $s\in\supp(\mu)$;
\item $\sum_{t\in\half\Int_+} \mu(t)(1+q^2)^{2t}<\infty$,
\end{enumeraterm}
then the random walk on $\widehat\SUqt$ defined by $\mu$ converges to the boundary.
\end{Thm}
\begin{pf}
We postpone the proof of the regularity of the elements $K_{\bar\mu}(x)$ for $x\in c_{00}(\widehat{\SUqt})$ to the next subsection. For now assume that the first part of the conjecture holds.

The space $K_{\check\varphi_\mu}(c_{00}(\widehat{\SUqt}))\subset M(\hat G,\mu)$ is dense by \cite[Thm.\ 4.10]{NeshveyevTuset04}. Moreover by \cite[Cor.\ 4.13, Prop.\ 3.11]{NeshveyevTuset04} the mapping $K_{\check\varphi_\mu}^*$ intertwines the right coactions of $(c_0(\hat G),\hat\Delta)$. Now Corollary \ref{sufficient_conditions2_conj2} gives the second part of the conjecture of convergence to the boundary.
\end{pf}

Note that Collins also worked on the representation of harmonic functions of the (classical) random walk on the center of $\widehat{\textrm{SU}_q(N)}$, see \cite[Thm.\ 4.1]{Collins04}.

To prove the regularity of elements in the Martin boundary, we need some knowledge of $\tilde M(\widehat{\SUqt},\mu)$. Here we give a short recap of known results. We follow the conventions of \cite{NeshveyevTuset04}.

The Hopf $*$-algebra $U_q(\sufrak_2)$ is generated by elements $k,e,f$ satisfying the following relations:
\begin{align*}
kk^{-1}&=k^{-1}k=1, & 			ke&=qek, &			kf&=q^{-1}fk;\\
[e,f]&=\frac{k^2-k^{-2}}{q-q^{-1}}; \\
\hat\Delta(k)&=k\ten k, & 		\hat\Delta(e)&=e\ten k^{-1}+ k\ten e, & 	\hat\Delta(f)&=f\ten k^{-1}+ k\ten f; \\
\hat{S}(k)&=k^{-1}, & 			\hat{S}(e)&=-q^{-1}e, &		\hat{S}(f)&=-qf; \\
\hat\eps(k) &=1, &			\hat\eps(e)&=\hat\eps(f)=0; \\
k^*&=k, &				e^*&=f,& 			f^*&=e.
\end{align*}
For $s\in\half\Int_+$ let $\{\xi_i^s\}_{i=-s}^s$ be an orthonormal basis of a $2s+1$ dimensional Hilbert space $\Hcal_s$. There exists a unitary representation of $U_q(\sufrak_2)$ on $\Hcal_s$ given by
\begin{align*}
\pi_s(k)\xi_j^s &= q^{-j}\xi_j^s, & \pi_s(e)\xi_j^s &= ([s+j]_q[s-j+1]_q)^{1/2}\xi_{j-1}^s, \\
\pi_s(f)\xi_j^s &= ([s-j]_q[s+j+1]_q)^{1/2}\xi_{j+1}^s.
\end{align*}
This representation is denoted by $U_s$. Up to equivalence these are all irreducible unitary representations. The direct sum of the representations $\{U_s\}_{s\in\half\Int_+}$ defines a faithful representation of the quantum enveloping algebra $U_q(\sufrak_2)$ in $\Ucal(\SUqt)$. The fusion rules of $U_q(\sufrak_2)$ are the same as for the Lie algebra $\sufrak_2$. These are explicitly known:
\begin{equation}\label{fusion_rules_su2_eq2}
U_t\times U_s \cong U_{t+s} \oplus U_{t+s-1} \oplus \ldots \oplus U_{|t-s|}.
\end{equation}
In particular $U_\half\times U_t\cong U_{t+\half}\oplus U_{t-\half}$. So there exist two orthonormal bases of $\Hcal_\half\ten \Hcal_t\cong \Hcal_{t-\half}\oplus\Hcal_{t+\half}$, namely $\{\xi^\half_i\ten\xi^t_j\}_{i,j}$ and $\{\xi_i^{t-\half}\}_i\cup\{\xi_j^{t+\half}\}_j$. These two bases can be expressed in terms of each other by means of inner products, the so-called {\it Clebsch--Gordan coefficients}.
Denote $C_q(s,t,r;j,k,m):=\langle \xi_j^s\ten \xi_k^t, \xi^r_m\rangle$. From \cite[Eq.\ (3.4.68), (3.4.69)]{KlimykSchmudgen97} it follows that\footnote{In \cite[Eq.\ (6.3.68)]{KlimykSchmudgen97} a factor $+1/2$ should be replaced by $+1$, this is done in our formulas.}
\begin{align*}
C_q\Big(\frac{1}{2},s,s+\frac{1}{2}; \pm\frac{1}{2},j,j\pm\frac{1}{2}\Big) &= q^{\frac{1}{2}(\mp s+j)}\Big(\frac{[s\pm j+1]_q}{[2s+1]_q}\Big)^{\frac{1}{2}}\\
C_q\Big(\frac{1}{2},s,s-\frac{1}{2}; \pm\frac{1}{2},j,j\pm\frac{1}{2}\Big) &= \mp q^{\frac{1}{2}(\pm s+j\pm1)}\Big(\frac{[s\mp j]_q}{[2s+1]_q}\Big)^{\frac{1}{2}}.
\end{align*}
To compute these Clebsch--Gordan coefficients, \cite{KlimykSchmudgen97} makes use of an algebra $\breve{U}_q(\sufrak_2)$. This algebra is different, but isomorphic to $U_q(\sufrak_2)$. This algebra $\breve{U}_q(\sufrak_2)$ considered by Klimyk--Schm\"udgen has generators $E,F,K$. Using the isomorphism given on these generators by
\begin{align*}
\psi\colon\breve{U}_q(\sufrak_2)&\ra U_q(\sufrak_2),\\
\psi(E)&:=f, & \psi(F)&:=e, & \psi(K)&:=k^{-1}.
\end{align*}
we get the above Clebsch--Gordan coefficients from \cite{KlimykSchmudgen97}.

Note that the fusion rules \eqref{fusion_rules_su2_eq2} imply that if $s\in\half\Int_+\setminus\Int_+$, then for any $t\in\half\Int_+$ there exists $n\in\Nat$ such that $U_t$ is a subrepresentation of $U_s^{\ten n}$. So $s$ generates all irreducible representations. However, if $s\in\Int_+\subset\half\Int_+$ and $t\in\half\Int_+$, then $U_t$ is a subrepresentation of $U_s^{\ten n}$ for some $n\in\Nat$ if and only if $t\in\Int_+$. So a probability measure on $\Irr(\SUqt)$ is generating if and only if $\supp(\mu)\cap(\half\Int_+\setminus\Int_+)\neq\emptyset$. A last thing to note is that $\rho=k^{-2}$ (see \cite[Eq.\ (4.6)]{NeshveyevTuset04}).

In the same spirit as Lemma \ref{estimates_multiplicities} we have the following (stronger estimate) for $\SUqt$
\begin{equation}\label{estimates_multiplicities_suqt_eq1}
\sum_{t\in\half\Int_+} m_{\half^{\ten n}\ten s}^t\leq 2^n.
\end{equation}
Indeed, use \eqref{fusion_rules_su2_eq2} to obtain
\[
U_\half\ten U_t\cong \begin{cases}
U_{t-\half}\oplus U_{t+\half}, &\textrm{if } t\geq \half;\\
U_{\half}, &\textrm{if } t=0.
\end{cases}
\]
Hence by induction
\[
\sum_{r\in\half\Int_+} m_{\half^{\ten n+1}\ten s}^r = \sum_{r\in\half\Int_+} \sum_{t\in\half\Int_+} m_{\half^{\ten n}\ten s}^t m_{\half\ten t}^r \leq 2^n 2= 2^{n+1},
\]
as desired.

\begin{Not}\label{def_generators}
In accordance with \cite{NeshveyevTuset04} define $\lambda,\,\tilde\lambda\in\prod_s B(\Hcal_s)$ by
\begin{align*}
\pi_s(\lambda):= \frac{q(q^{2s+1}+q^{-2s-1})}{(q-q^{-1})\sqrt{[2]_q}}\,I_s, && \pi_s(\tilde\lambda):=q^{-2s}I_s.
\end{align*}
Define $\chi_{-1}:= -qfk$, $\chi_0:=\frac{ef-q^2fe}{\sqrt{[2]_q}}$ and $\chi_1:=qek$. Denote $X_j:=\lambda^{-1}\chi_j$ for $j=-1,\,0,\, 1$. Moreover write $\tilde X_{-1}:=\tilde\lambda^{-1}fk$, $\tilde X_0:=\tilde\lambda^{-1} k^2$ and $\tilde X_1:= \tilde{\lambda}^{-1}ek$. It can be shown that $X_j,\,\tilde X_j\in l^\infty(\widehat{\SUqt})$.
Let $\Psi$ be the unital \Cstar algebra in $l^\infty(\widehat \SUqt)$ generated by $c_0(\widehat\SUqt)$ and the elements $X_i$, $i=-1,0,1$. The {\it quantum homogeneous sphere of Podle\'s} \cite{Podles87} is the \Cstar algebra $C(S_{q,0}^2):=\Psi /c_0(\widehat\SUqt)$.
\end{Not}

\begin{Lem}\label{Lem_span_generators}
The elements $\tilde X_i$ and $X_i$ satisfy the following:
\begin{align*}
\Com X_i + c_0(\widehat\SUqt) &= \Com\tilde X_i + c_0(\widehat\SUqt), \qquad (\textrm{for } i=-1,\,1);\\
\Com X_0 + c_0(\widehat\SUqt) &\subset \Com 1 + \Com\tilde X_0 + c_0(\widehat\SUqt);\\
\Com \tilde X_0 + c_0(\widehat\SUqt) &\subset \Com 1 + \Com X_0 + c_0(\widehat\SUqt).
\end{align*}
Here $\Com X:=\{t X\,:\, t\in\Com\}$ denotes the complex linear span. In particular $\Psi$ equals the unital \Cstar subalgebra of $l^\infty(\widehat{\SUqt})$ generated by $\tilde X_i$, $i=-1,\,0,\,1$ and $c_0(\widehat \SUqt)$.
\end{Lem}
\begin{pf}
We have
\begin{align*}
\pi_s(\tilde\lambda^{-1}\lambda) &= q^{2s}\frac{q(q^{2s+1}+q^{-2s-1})}{(q-q^{-1})\sqrt{[2]_q}}\,I_s = \frac{1}{(q-q^{-1})\sqrt{[2]_q}}\,I_s + \frac{q^{4s+2}}{(q-q^{-1})\sqrt{[2]_q}}\,I_s;\\
\pi_s(\tilde\lambda\lambda^{-1}) &= q^{-2s}\frac{(q-q^{-1})\sqrt{[2]_q}}{q(q^{2s+1}+q^{-2s-1})}\,I_s = (q-q^{-1})\sqrt{[2]_q}\,I_s + (q-q^{-1})\sqrt{[2]_q}\frac{-q^{4s+2}}{q^{4s+2}+1}\,I_s.
\end{align*}
Write $c:= (q-q^{-1})\sqrt{[2]_q}$. It follows that $\tilde\lambda^{-1}\lambda = c^{-1}1+ a$ and $\tilde\lambda\lambda^{-1} = c1+ a'$ for some $a,\,a'\in c_0(\widehat\SUqt)$. Therefore if $j=-1,\,1$
\begin{align*}
X_j&= \lambda^{-1}\chi_j = (\lambda^{-1}\tilde\lambda)\tilde\lambda^{-1}\chi_j = (c1+a')\tilde\lambda^{-1}\chi_j = jqc\tilde X_j + jqa'\tilde X_j \in \Com\tilde X_j + c_0(\widehat\SUqt); \displaybreak[2]\\
\tilde X_j&= \tilde\lambda^{-1}jq^{-1}\chi_j = (\tilde\lambda^{-1}\lambda)\lambda^{-1}jq^{-1}\chi_j = (c^{-1}1+a)\lambda^{-1}jq^{-1}\chi_j\\
&= c^{-1}jq^{-1}X_j + jq^{-1}aX_j \in \Com X_j + c_0(\widehat\SUqt).
\end{align*}
Hence for $j=-1,\,1$ the statement follows. For $j=0$ observe first that
\begin{align*}
\pi_s(\chi_0)\xi_i^s &= \pi_s\Big(\frac{ef-q^2fe}{\sqrt{[2]_q}}\Big)\xi_i^s = \frac{1}{\sqrt{[2]_q}}([s-i]_q[s+i+1]_q-q^2[s+i]_q[s-i+1]_q)\xi_i^s \\
&= \frac{1}{\sqrt{[2]_q}(q-q^{-1})^2}\\
&\qquad\times\big((q^{2s+1}+q^{-2s-1}-q^{2i+1}-q^{-2i-1}) - q^2(q^{2s+1}+q^{-2s-1}-q^{-2i+1}-q^{2i-1})\big)\xi_i^s \displaybreak[2]\\
&=\frac{1}{\sqrt{[2]_q}}\Big(\frac{1-q^2}{(q-q^{-1})^2}(q^{2s+1}+q^{-2s-1}) + \frac{1-q^4}{(q-q^{-1})^2}q^{-2i-1}\Big)\xi_i^s\\
&= \pi_s(\lambda)\xi_i^s + \frac{q(q^{-2}-q^2)\sqrt{[2]_q}}{(q^{-1}+q)(q-q^{-1})^2}q^{-2i}\xi_i^s\\
&= \pi_s\Big(-\lambda +\frac{q\sqrt{[2]_q}}{q-q^{-1}} k^2\Big)\xi_i^s.
\end{align*}
Thus $\chi_0=-\lambda+\frac{q\sqrt{[2]_q}}{q-q^{-1}} k^2$. Hence
\begin{align*}
X_0&= -1 + \lambda^{-1}\frac{q\sqrt{[2]_q}}{q-q^{-1}} k^2 = -1 + \frac{q\sqrt{[2]_q}}{q-q^{-1}}\lambda^{-1}\tilde\lambda \tilde\lambda^{-1} k^2 = -1 + \frac{q\sqrt{[2]_q}}{q-q^{-1}} (c1+ a')\tilde X_0; \displaybreak[2]\\
\tilde X_0 &= \tilde\lambda^{-1}k^2 = \Big(\frac{q-q^{-1}}{q\sqrt{[2]_q}}\lambda\tilde\lambda^{-1}\Big)\lambda^{-1}\frac{q\sqrt{[2]_q}}{q-q^{-1}}k^2\\
&= \Big(\frac{q-q^{-1}}{q\sqrt{[2]_q}}(c^{-1}1+a)\Big)\lambda^{-1}\frac{q\sqrt{[2]_q}}{q-q^{-1}}k^2 \displaybreak[2]\\
&= \frac{q-q^{-1}}{q\sqrt{[2]_q}}c^{-1} + \frac{q-q^{-1}}{q\sqrt{[2]_q}}c^{-1}\lambda^{-1}\Big(-\lambda + \frac{q\sqrt{[2]_q}}{q-q^{-1}}k^2\Big) + a\lambda^{-1}k^2\\
&= \frac{q-q^{-1}}{q\sqrt{[2]_q}}c^{-1} + \frac{q-q^{-1}}{q\sqrt{[2]_q}}c^{-1}X_0 + a\lambda^{-1}k^2,
\end{align*}
which proves the statement for $j=0$.
\end{pf}

\begin{Thm}[\protect{\cite[Thm.\ 4.1, 4.10]{NeshveyevTuset04}}]\label{Thm_Martin_SUqt}
Assume that $\mu$ is a generating and transient probability measure on $\half\Int_+$ with finite mean, that is $\sum_{s\in\half\Int_+} \mu(s)s<\infty$. Then the Martin compactification $\tilde M(\widehat\SUqt,\mu)$ of $\widehat{\SUqt}$ with respect to $\mu$ equals $\Psi$ and thus the Martin boundary is isomorphic to $C(S_{q,0}^2)$.
\end{Thm}

\subsection{Regularity}\label{sec_strong_star_convergence}
Let $\mu=\sum_{t\in\half\Int_+} \mu_t\delta_t$ be a probability measure on $\half\Int_+$. We form the infinite tensor product with respect to the state $\varphi_\mu$
\begin{equation}\label{def_inf_ten_prod_SUq2_eq1}
\bigotimes_{n=-\infty}^{-1} (l^\infty(\widehat{\SUqt}),\varphi_\mu).
\end{equation}
Recall that the state $\varphi_\mu^\infty$ given by
\[
\varphi^\infty_\mu(\cdots 1\ten 1\ten x_n\ten\cdots\ten x_1):= \varphi_\mu(x_n)\cdots\varphi_\mu(x_1)
\]
is faithful on the infinite tensor product \eqref{def_inf_ten_prod_SUq2_eq1}.

\begin{Prop}\label{prop_estimates_generators}
There exists a constant $C>0$ such that the elements $\tilde X_i$ satisfy the estimate
\[
\|1_{\half}\ten\pi_s(\tilde X_i)-(\pi_{\half}\ten\pi_s)\hat\Delta(\tilde X_i)\|_{\varphi_{\half}\ten\varphi_s}^2 \leq C\, d_s^{-1}.
\]
for every $s\in\half\Int_+$.
\end{Prop}
\begin{pf} To shorten the notation slightly, write
\begin{align}
a^{+\half}_{\pm,i}:=C_q\Big(\half,s,s+\half;\pm\half,i,i\pm\half\Big); &&
a^{-\half}_{\pm,i}:=C_q\Big(\half,s,s-\half;\pm\half,i,i\pm\half\Big), \label{prop_estimates_generators_eq1}
\end{align}
so that $\xi_{\pm\half}^{\half}\ten\xi_i^s = a^{+\half}_{\pm,i}\xi_{i\pm\half}^{s+\half} + a^{-\half}_{\pm,i}\xi_{i\pm\half}^{s-\half}$. Recall that $\langle \xi_i^s,\xi_j^t\rangle = \delta_{s,t}\delta_{i,j}$.
We compute
\begin{align}
&d_\half d_s(\varphi_\half\ten\varphi_s)\big((1_{\half}\ten\pi_s(\tilde X_0))^*(\pi_\half\ten\pi_s)(\hat\Delta(\tilde X_0))\big) \notag\\
&=(\Tr_\half\ten\Tr_s)\big(q^{2s}(1_\half\ten\pi_s(k^2))(\pi_\half\ten\pi_s)(\hat\Delta(\tilde\lambda^{-1})(k^2\ten k^2)(\rho^{-1}\ten\rho^{-1}))\big) \displaybreak[2] \notag\\
&= q^{2s} (\Tr_\half\ten\Tr_s)\big((1_\half\ten\pi_s(k^2))(\pi_\half\ten\pi_s)(\hat\Delta(\tilde\lambda^{-1})(k^4\ten k^4))\big) \displaybreak[2] \notag\\
&= q^{2s} \sum_{\substack{j\in\{\half,\,-\half\}\\ i\in\{-s,-s+1,\ldots,s\}}} \big\langle(1_\half\ten\pi_s(k^2))(\pi_\half\ten\pi_s)(\hat\Delta(\tilde\lambda^{-1})(k^4\ten k^4))\xi_j^\half\ten\xi_i^s,\xi_j^\half\ten\xi_i^s \big\rangle \displaybreak[2] \notag\\
&= q^{2s} \sum_{\substack{j\in\{\half,\,-\half\}\\ i\in\{-s,-s+1,\ldots,s\}}} q^{-4j-4i} \big\langle(\pi_\half\ten\pi_s)(\hat\Delta(\tilde\lambda^{-1}))\xi_j^\half\ten\xi_i^s, (1_\half\ten\pi_s(k^2))\xi_j^\half\ten\xi_i^s \big\rangle \displaybreak[2] \notag\\
&= q^{2s} \sum_{\pm,\, i} q^{-6i\mp2} \big\langle(\pi_\half\ten\pi_s)(\hat\Delta(\tilde\lambda^{-1}))\xi_{\pm\half}^\half\ten\xi_i^s, \xi_{\pm\half}^\half\ten\xi_i^s \big\rangle \displaybreak[2] \notag\\
&= q^{2s} \sum_{\pm,\, i} q^{-6i\mp2} \Big\langle a^{+\half}_{\pm,i}\,\pi_{s+\half}(\tilde\lambda^{-1})\xi_{i\pm\half}^{s+\half} + a^{-\half}_{\pm,i}\, \pi_{s-\half}(\tilde\lambda^{-1})\xi_{i\pm\half}^{s-\half}, a^{+\half}_{\pm,i}\xi_{i\pm\half}^{s+\half} + a^{-\half}_{\pm,i}\xi_{i\pm\half}^{s-\half} \Big\rangle \displaybreak[2] \notag\\
&= q^{2s} \sum_{\pm,\, i} q^{-6i\mp2} \Big\langle q^{2s+1} a^{+\half}_{\pm,i}\xi_{i\pm\half}^{s+\half} + q^{2s-1}a^{-\half}_{\pm,i}\xi_{i\pm\half}^{s-\half}, a^{+\half}_{\pm,i}\xi_{i\pm\half}^{s+\half} + a^{-\half}_{\pm,i}\xi_{i\pm\half}^{s-\half} \Big\rangle \displaybreak[2]  \notag\\
&= q^{4s} \sum_{\pm,\, i} q^{-6i\mp2} \big(q (a^{+\half}_{\pm,i})^2 + q^{-1} (a^{-\half}_{\pm,i})^2 \big) \notag\\
&= q^{4s} \sum_{\substack{\pm \\ i\in\{-s,\ldots, s-1,s\}}} q^{-6i\mp2} \Big(qq^{(\mp s +i)}\frac{[s\pm i+1]_q}{[2s+1]_q} + q^{-1}q^{(\pm s +i\pm1)}\frac{[s\mp i]_q}{[2s+1]_q}\Big). \label{prop_estimates_generators_eq5}
\end{align}
It is a standard result that for a positive linear functional $\omega$ it holds that $\omega(a^*)=\overline{\omega(a)}$. Clearly $\varphi_s$ is positive, hence
\begin{equation}\label{prop_estimates_generators_eq2}
\varphi_s(ab) = \overline{\varphi_s(b^*a^*)}.
\end{equation}
Since $\tilde X_0$ is self-adjoint and \eqref{prop_estimates_generators_eq5} is real-valued we see
\begin{align*}
&d_\half d_s(\varphi_\half\ten\varphi_s)\big((\pi_\half\ten\pi_s)(\hat\Delta(\tilde X_0))^* (1_{\half}\ten\pi_s(\tilde X_0))\big)\\
&= q^{4s} \sum_{\substack{\pm \\ i\in\{-s,\ldots, s-1,s\}}} q^{-6i\mp2} \Big(qq^{(\mp s +i)}\frac{[s\pm i+1]_q}{[2s+1]_q} + q^{-1}q^{(\pm s +i\pm1)}\frac{[s\mp i]_q}{[2s+1]_q}\Big).
\end{align*}
Similarly using that $\tilde\lambda$ is central gives
\begin{align*}
&d_\half d_s(\varphi_\half\ten\varphi_s)\big((\pi_\half\ten\pi_s)(\hat\Delta(\tilde X_0))^*(\pi_\half\ten\pi_s)(\hat\Delta(\tilde X_0))\big)\\
&= (\Tr_\half\ten\Tr_s)\big((\pi_\half\ten\pi_s)(\hat\Delta(\tilde\lambda^{-2}k^4)(\rho^{-1}\ten\rho^{-1}))\big) \displaybreak[2]\\
&= (\Tr_\half\ten\Tr_s)\big((\pi_\half\ten\pi_s)(\hat\Delta(\tilde\lambda^{-2})(k^6\ten k^6)))\big) \displaybreak[2]\\
&= \sum_{\substack{\pm \\ i\in\{-s,\ldots, s-1,s\}}} \big\langle (\pi_\half\ten\pi_s)(\hat\Delta(\tilde\lambda^{-2})(k^6\ten k^6))\xi_{\pm\half}^\half\ten\xi_i^s,\xi_{\pm\half}^\half\ten\xi_i^s \big\rangle \displaybreak[2]\\
&= \sum_{\pm,\,i} q^{(-6i\mp 3)}\big\langle (\pi_\half\ten\pi_s)(\hat\Delta(\tilde\lambda^{-2}))\xi_{\pm\half}^\half\ten\xi_i^s,\xi_{\pm\half}^\half\ten\xi_i^s \big\rangle \displaybreak[2]\\
&= \sum_{\pm,\,i} q^{(-6i\mp 3)}\Big\langle a^{+\half}_{\pm,i}\, \pi_{s+\half}(\tilde\lambda^{-2})\xi_{i\pm\half}^{s+\half} + a^{-\half}_{\pm,i} \, \pi_{s-\half}(\tilde\lambda^{-2}) \xi_{i\pm\half}^{s-\half}, a^{+\half}_{\pm,i}\xi_{i\pm\half}^{s+\half} + a^{-\half}_{\pm,i}\xi_{i\pm\half}^{s-\half} \Big\rangle \displaybreak[2]\\
&= \sum_{\pm,\,i} q^{(-6i\mp 3)}\Big\langle q^{(4s+2)}a^{+\half}_{\pm,i}\xi_{i\pm\half}^{s+\half} + q^{(4s-2)}a^{-\half}_{\pm,i}\xi_{i\pm\half}^{s-\half}, a^{+\half}_{\pm,i}\xi_{i\pm\half}^{s+\half} + a^{-\half}_{\pm,i}\xi_{i\pm\half}^{s-\half} \Big\rangle\\
&= q^{4s}\sum_{\pm,\, i} q^{(-6i\mp 3)}\Big(q^2 q^{(\mp s +i)}\frac{[s\pm i+1]_q}{[2s+1]_q} + q^{-2}q^{(\pm s +i\pm1)}\frac{[s\mp i]_q}{[2s+1]_q}\Big).
\end{align*}
Last using that
\[
1=\langle \xi_{\pm\half}^{\half}\ten\xi_i^s, \xi_{\pm\half}^{\half}\ten\xi_i^s\rangle = \langle a^{+\half}_{\pm,i}\xi_{i\pm\half}^{s+\half} + a^{-\half}_{\pm,i}\xi_{i\pm\half}^{s-\half}, a^{+\half}_{\pm,i}\xi_{i\pm\half}^{s+\half} + a^{-\half}_{\pm,i}\xi_{i\pm\half}^{s-\half}\rangle
= (a^{+\half}_{\pm,i})^2 + (a^{-\half}_{\pm,i})^2
\]
gives
\begin{align*}
&d_\half d_s(\varphi_\half\ten\varphi_s)\big((1_\half\ten\pi_s(\tilde X_0))^*(1_\half\ten\pi_s(\tilde X_0))\big)\\
&=(\Tr_\half\ten\Tr_s)\big((\pi_\half\ten\pi_s)((1_\half\ten(\tilde\lambda^{-2}k^4))(\rho^{-1}\ten\rho^{-1}))\big) \displaybreak[2]\\
&= q^{4s} \sum_{\pm,\, i}\big\langle (\pi_\half(k^2)\ten\pi_s(k^6))\xi_{\pm\half}^\half\ten\xi_i^s,\xi_{\pm\half}^\half\ten\xi_i^s \big\rangle \displaybreak[2]\\
&= q^{4s} \sum_{\pm,\, i} q^{(-6i\mp1)} \\
&= q^{4s} \sum_{\pm,\, i} q^{(-6i\mp1)}\Big(q^{(\mp s +i)}\frac{[s\pm i+1]_q}{[2s+1]_q} + q^{(\pm s +i\pm1)}\frac{[s\mp i]_q}{[2s+1]_q}\Big) .
\end{align*}
Combining these four calculations yields
\begin{align*}
&(\varphi_\half\ten\varphi_s)\big((1_\half\ten\pi_s(\tilde X_0)-(\pi_\half\ten\pi_s)(\hat\Delta(\tilde X_0)))^*(1_\half\ten\pi_s(\tilde X_0)-(\pi_\half\ten\pi_s)(\hat\Delta(\tilde X_0)))\big)\\
&= (\varphi_\half\ten\varphi_s)\big((1_\half\ten\pi_s(\tilde X_0))^*(1_\half\ten\pi_s(\tilde X_0))\big)\\
&\qquad+ (\varphi_\half\ten\varphi_s)\big((\pi_\half\ten\pi_s)(\hat\Delta(\tilde X_0))^*(\pi_\half\ten\pi_s)(\hat\Delta(\tilde X_0))\big) \\
&\qquad -2  (\varphi_\half\ten\varphi_s)\big((1_{\half}\ten\pi_s(\tilde X_0))^*(\pi_\half\ten\pi_s)(\hat\Delta(\tilde X_0))\big) \displaybreak[2]\\
&= (d_\half d_s)^{-1}q^{4s}\Big( \sum_{\pm,\, i} q^{(-6i\mp1)}\Big(q^{(\mp s +i)}\frac{[s\pm i+1]_q}{[2s+1]_q} + q^{(\pm s +i\pm1)}\frac{[s\mp i]_q}{[2s+1]_q}\Big)\\
&\qquad + \sum_{\pm,\, i}q^{(-6i\mp 3)}\Big(q^2 q^{(\mp s +i)}\frac{[s\pm i+1]_q}{[2s+1]_q} + q^{-2}q^{(\pm s +i\pm1)}\frac{[s\mp i]_q}{[2s+1]_q}\Big)\\
&\qquad -2 \sum_{\pm,\, i} q^{(-6i\mp2)} \Big(qq^{(\mp s +i)}\frac{[s\pm i+1]_q}{[2s+1]_q} + q^{-1}q^{(\pm s +i\pm1)}\frac{[s\mp i]_q}{[2s+1]_q}\Big)\Big) \displaybreak[2]\\
&= (d_\half d_s)^{-1} [2s+1]_q^{-1}q^{4s}\Big( \sum_{\pm,\, i} q^{(-6i\mp1)}\Big( (1+q^{\mp2}q^2-2q^{\mp1}q)q^{(\mp s +i)}[s\pm i+1]_q \\
&\qquad + (1+q^{\mp2}q^{-2}-2q^{\mp1}q^{-1})q^{(\pm s +i\pm1)}[s\mp i]_q\Big)\Big) \displaybreak[2]\\
&= d_\half^{-1} d_s^{-2}q^{4s}\Big(\sum_{i=-s}^s q^{(-6i-1)}\Big( (1+1-2)q^{(i-s)}[s+i+1]_q
+ (1+q^{-4}-2q^{-2})q^{(s +i+1)}[s-i]_q\Big)\\
&\qquad + \sum_{i=-s}^s q^{(-6i+1)}\Big( (1+q^4-2q^2)q^{(s+i)}[s-i+1]_q
+ (1+1-2)q^{(-s +i-1)}[s+i]_q\Big)\Big) \displaybreak[2]\\
&= d_\half^{-1} d_s^{-2}q^{4s}\sum_{i=-s}^s q^{-6i}\Big(q^{-1}(1-q^{-2})^2q^{(s +i+1)}[s-i]_q\Big)+ \Big( q(1-q^2)^2q^{(s+i)}[s-i+1]_q\Big) \displaybreak[2]\\
&= d_\half^{-1} d_s^{-2}\frac{q^{4s}}{(q-q^{-1})}\sum_{i=-s}^s q^{-6i}\Big(q^{-3}(q-q^{-1})^2q^{(s +i+1)}(q^{s-i}-q^{-s+i})\Big)\\
&\qquad + \Big( q^3(q-q^{-1})^2q^{(s+i)}(q^{s-i+1}-q^{-s+i-1})\Big) \displaybreak[2]\\
&= d_\half^{-1} d_s^{-2}q^{4s}(q-q^{-1})\sum_{i=-s}^s q^{-6i}\big(q^{2s-2}-q^{2i-2}+ q^{2s+4}-q^{2i+2}\big) \displaybreak[2]\\
&= d_\half^{-1} d_s^{-2}q^{4s}(q-q^{-1})\Big((q^{2s-2}+q^{2s+4})\sum_{i=-s}^s q^{-6i} - (q^{-2}+q^2)\sum_{i=-s}^s q^{-4i}\Big) \displaybreak[2]\\
&= d_\half^{-1} d_s^{-2}q^{4s}(q-q^{-1})\big((q^{2s-2}+q^{2s+4})[2s+1]_{q^3} - (q^{-2}+q^2)[2s+1]_{q^2}\big) \displaybreak[2]\\
&= d_\half^{-1} d_s^{-2}q^{4s}(q-q^{-1})\big(O(q^{-4s}) +O(q^{-4s})\big)\\
&= d_\half^{-1} d_s^{-2} O(1),
\end{align*}
which proves the result for $\tilde X_0$. In fact this estimate is stronger than stated in the proposition, but we will not need that.

We deal with $\tilde{X}_1$ in an analogous way, however the estimates become slightly more involved, because the $\xi_i^s$ are no longer eigenvectors for $\pi_s(e)$ and $\pi_s(f)$. The calculations for $\tilde X_{-1}$ are similar and are omitted. First we calculate
\begin{align}
d_s\varphi_s\big(\pi_s(\tilde X_1^*\tilde X_1)\big) &=\Tr_s\big(\pi_s(\tilde\lambda^{-1}kf\tilde\lambda^{-1}ek)\pi_s(\rho^{-1})\big) \notag \displaybreak[2]\\
&= \Tr_s\big(\pi_s(\tilde\lambda^{-2}fek^4)\big) \label{prop_estimates_generators_eq3} \displaybreak[2]\\
&= q^{4s}\sum_{i=-s}^s\langle\pi_s(fek^4)\xi_i^s,\xi_i^s\rangle \notag \displaybreak[2]\\
&= q^{4s}\sum_{i=-s}^s q^{-4i}([s+i]_q[s-i+1]_q)^\half([s-(i-1)]_q[s+(i-1)+1]_q)^\half\langle\xi_i^s,\xi_i^s\rangle \notag \displaybreak[2]\\
&= q^{4s}\sum_{i=-s}^s q^{-4i}[s+i]_q[s-i+1]_q \label{prop_estimates_generators_eq4} \displaybreak[2]\\
&= \frac{q^{4s}}{(q-q^{-1})^2}\sum_{i=-s}^s q^{-4i}(q^{s+i}-q^{-s-i})(q^{s-i+1}-q^{-s+i-1}) \notag \displaybreak[2]\\
&= \frac{q^{4s}}{(q-q^{-1})^2}\sum_{i=-s}^s q^{-4i}(q^{2s+1}+q^{-2s-1}-q^{2i-1}-q^{-2i+1}) \notag \displaybreak[2]\\
&= \frac{q^{4s}}{(q-q^{-1})^2}\Big((q^{2s+1}+q^{-2s-1})\sum_{i=-s}^s q^{-4i}-q^{-1}\sum_{i=-s}^s q^{-2i} -q\sum_{i=-s}^s q^{-6i}\Big) \notag\\
&= \frac{q^{4s}}{(q-q^{-1})^2}\Big((q^{2s+1}+q^{-2s-1})[2s+1]_{q^2} - q^{-1}[2s+1]_q - q[2s+1]_{q^3}\Big). \notag
\end{align}
This gives
\begin{align*}
&d_\half d_s(\varphi_\half\ten\varphi_s)\big((1_\half\ten\pi_s(\tilde X_1))^*(1_\half\ten\pi_s(\tilde X_1))\big)\\
&= q^{4s}\frac{q+q^{-1}}{(q-q^{-1})^2}\Big((q^{2s+1}+q^{-2s-1})[2s+1]_{q^2} - q^{-1}[2s+1]_q - q[2s+1]_{q^3}\Big).
\end{align*}
Using \eqref{prop_estimates_generators_eq3} for $s\pm\half$ instead of $s$ we also get
\begin{align*}
&d_\half d_s(\varphi_\half\ten\varphi_s) \big((\pi_\half\ten\pi_s)(\hat\Delta(\tilde\lambda^{-1}ek)^*\hat\Delta(\tilde\lambda^{-1} ek))\big)\\
&= (\Tr_\half\ten\Tr_s) \big((\pi_\half\ten\pi_s)(\hat\Delta(kf\tilde\lambda^{-1}\tilde\lambda^{-1} ek)(k^2\ten k^2))\big) \displaybreak[2]\\
&= (\Tr_\half\ten\Tr_s) \big((\pi_\half\ten\pi_s)(\hat\Delta(\tilde\lambda^{-2} fek^4))\big) \displaybreak[2]\\
&= \Tr_{s-\half} \big(\pi_{s-\half}(\tilde\lambda^{-2} fek^4)\big)+ \Tr_{s+\half} \big(\pi_{s+\half}(\tilde\lambda^{-2} fek^4)\big) \displaybreak[2]\\
&= \frac{q^{4s-2}}{(q-q^{-1})^2}\Big((q^{2s}+q^{-2s})[2s]_{q^2} - q^{-1}[2s]_q - q[2s]_{q^3}\Big)\\
&\qquad + \frac{q^{4s+2}}{(q-q^{-1})^2}\Big((q^{2s+2}+q^{-2s-2})[2s+2]_{q^2} - q^{-1}[2s+2]_q - q[2s+2]_{q^3}\Big).
\end{align*}
Recall the abbreviations of the Clebsch--Gordan coefficients $a_{\pm,i}^{+\half}$ and $a_{\pm,i}^{-\half}$ introduced in \eqref{prop_estimates_generators_eq1} above. We obtain
\begin{align*}
&d_\half d_s(\varphi_\half\ten\varphi_s)\big((1_{\half}\ten\pi_s(\tilde X_1))^*(\pi_\half\ten\pi_s)(\hat\Delta(\tilde X_1))\big)\\
&= (\Tr_\half\ten\Tr_s)\big((1_{\half}\ten\pi_s(\tilde \lambda^{-1}ek))^*(\pi_\half\ten\pi_s)(\hat\Delta(\tilde \lambda^{-1}ek))(k^2\ten k^2)\big) \displaybreak[2]\\
&= (\Tr_\half\ten\Tr_s)\big(q^{2s}(\pi_\half\ten\pi_s)((1\ten kf)\hat\Delta(\tilde \lambda^{-1}e)(k^3\ten k^3))\big) \displaybreak[2]\\
&= q^{2s} \sum_{\substack{\pm\\ i\in\{-s,\ldots, s-1,s\}}} \big\langle (\pi_\half\ten\pi_s)((1\ten kf)\hat\Delta(\tilde \lambda^{-1}e)(k^3\ten k^3))\xi_{\pm\half}^\half\ten\xi_i^s, \xi_{\pm\half}^\half\ten\xi_i^s\big\rangle \displaybreak[2]\\
&= q^{2s} \sum_{\pm,\, i} q^{\mp\frac{3}{2}}q^{-3i} \big\langle (\pi_\half\ten\pi_s)(\hat\Delta(\tilde \lambda^{-1}e))\xi_{\pm\half}^\half\ten\xi_i^s, (\pi_\half\ten\pi_s)(1\ten ek) \xi_{\pm\half}^\half\ten\xi_i^s\big\rangle \displaybreak[2]\\
&= q^{2s} \sum_{\pm,\, i} q^{\mp\frac{3}{2}}q^{-3i} q^{-i} ([s+i]_q[s-i+1]_q)^\half \big\langle (\pi_\half\ten\pi_s)(\hat\Delta(\tilde \lambda^{-1}e))\xi_{\pm\half}^\half\ten\xi_i^s, \xi_{\pm\half}^\half\ten\xi_{i-1}^s\big\rangle \displaybreak[2]\\
&=q^{2s} \sum_{\pm,\, i} q^{\mp\frac{3}{2}}q^{-4i} ([s+i]_q[s-i+1]_q)^\half\\
&\quad \times\big\langle a^{+\half}_{\pm,i}\pi_{s+\half}(\tilde \lambda^{-1}e)\xi_{i\pm\half}^{s+\half} + a^{-\half}_{\pm,i}\pi_{s-\half}(\tilde \lambda^{-1}e)\xi_{i\pm\half}^{s-\half}, a^{+\half}_{\pm,i-1}\xi_{i-1\pm\half}^{s+\half} + a^{-\half}_{\pm,i-1}\xi_{i-1\pm\half}^{s-\half} \big\rangle \displaybreak[2]\\
&= q^{2s} \sum_{\pm,\, i} q^{\mp\frac{3}{2}}q^{-4i} ([s+i]_q[s-i+1]_q)^\half\\
&\quad \times\big\langle a^{+\half}_{\pm,i}q^{2s+1}\big([s+\half+i\pm\half]_q[s+\half-i\mp\half+1]_q\big)^\half\xi_{i\pm\half-1}^{s+\half} \\
&\qquad+a^{-\half}_{\pm,i}q^{2s-1}\big([s-\half+i\pm\half]_q[s-\half-i\mp\half+1]_q\big)^\half\xi_{i\pm\half-1}^{s-\half}, a^{+\half}_{\pm,i-1}\xi_{i-1\pm\half}^{s+\half} + a^{-\half}_{\pm,i-1}\xi_{i-1\pm\half}^{s-\half} \big\rangle \displaybreak[2]\\
&= q^{4s} \sum_{\pm,\, i} q^{\mp\frac{3}{2}}q^{-4i} \Big(q\big([s+i]_q[s-i+1]_q[s+i+\half\pm\half]_q[s-i+\frac{3}{2}\mp\half]_q\big)^\half a^{+\half}_{\pm,i} a^{+\half}_{\pm,i-1}\\
&\quad + q^{-1}\big([s+i]_q[s-i+1]_q[s+i-\half\pm\half]_q[s-i+\half\mp\half]_q\big)^\half a^{-\half}_{\pm,i}a^{-\half}_{\pm,i-1}\Big) \displaybreak[2]\\
&= q^{4s} \sum_{\pm,\, i} q^{\mp\frac{3}{2}}q^{-4i} \Big(q\big([s+i]_q[s-i+1]_q[s+i+\half\pm\half]_q[s-i+\frac{3}{2}\mp\half]_q\big)^\half\\
&\qquad \times q^{\half(\mp s+i)}\Big(\frac{[s\pm i+1]_q}{[2s+1]_q}\Big)^\half q^{\half(\mp s+i-1)}\Big(\frac{[s\pm i\mp1+1]_q}{[2s+1]_q}\Big)^\half \\
&\quad + q^{-1}\big([s+i]_q[s-i+1]_q[s+i-\half\pm\half]_q[s-i+\half\mp\half]_q\big)^\half \\
&\qquad \times q^{\half(\pm s+i\pm1)}\Big(\frac{[s\mp i]_q}{[2s+1]_q}\Big)^\half q^{\half(\pm s+i-1\pm1)}\Big(\frac{[s\mp i\pm1]_q}{[2s+1]_q}\Big)^\half\Big) \displaybreak[2]\\
&= \frac{q^{4s}}{[2s+1]_q} \sum_{\pm, i} q^{\mp\frac{3}{2}}q^{-4i} \Big(q\big([s+i]_q[s-i+1]_q [s+i+\half\pm\half]_q[s-i+\frac{3}{2}\mp\half]_q\\
&\qquad\times [s\pm i+1]_q[s\pm i\mp1+1]_q\big)^\half q^{\mp s+i-\half}\\
&\quad + q^{-1}\big([s+i]_q[s-i+1]_q[s+i-\half\pm\half]_q[s-i+\half\mp\half]_q\\
&\qquad\times[s\mp i]_q[s\mp i\pm1]_q\big)^\half q^{\pm s+i\pm1-\half}\Big) \displaybreak[2]\\
&= \frac{q^{4s}}{d_s} \sum_i q^{-\frac{3}{2}}q^{-4i} q\big([s+i]_q[s-i+1]_q [s+i+1]_q[s-i+1]_q [s+i+1]_q[s+i]_q\big)^\half q^{-s+i-\half}\\
&\quad + q^{\frac{3}{2}}q^{-4i} q\big([s+i]_q[s-i+1]_q [s+i]_q[s-i+2]_q [s-i+1]_q[s-i+2]_q\big)^\half q^{s+i-\half}\\
&\quad + q^{-\frac{3}{2}}q^{-4i}q^{-1}\big([s+i]_q[s-i+1]_q[s+i]_q[s-i]_q[s-i]_q[s-i+1]_q\big)^\half q^{s+i+\half}\\
&\quad + q^{\frac{3}{2}}q^{-4i}q^{-1}\big([s+i]_q[s-i+1]_q[s+i-1]_q[s-i+1]_q[s+i]_q[s+i-1]_q\big)^\half q^{-s+i-\frac{3}{2}} \displaybreak[2]\\
&= \frac{q^{4s}}{d_s} \sum_i q^{-4i} [s+i]_q[s-i+1]_q \Big(q^{-\frac{3}{2}} q[s+i+1]_q q^{-s+i-\half} + q^{\frac{3}{2}}q[s-i+2]_q q^{s+i-\half}\\
&\quad + q^{-\frac{3}{2}}q^{-1}[s-i]_q q^{s+i+\half} + q^{\frac{3}{2}}q^{-1}[s+i-1]_q q^{-s+i                                                                  -\frac{3}{2}}\Big) \displaybreak[2]\\
&= \frac{q^{4s}}{d_s} \sum_i q^{-4i} [s+i]_q[s-i+1]_q \Big([s+i+1]_q q^{-s+i-1} + [s-i+2]_q q^{s+i+2}
+[s-i]_q q^{s+i-2}\\
&\quad + [s+i-1]_q q^{-s+i-1}\Big) \displaybreak[2]\\
&=\frac{q^{4s}}{d_s} \sum_i q^{-4i} [s+i]_q[s-i+1]_q \Big(q[s+i+1+s-i+2]_q +q^{-1}[s-i+s+i-1]_q \Big) \displaybreak[2]\\
&=q^{4s}\,\frac{q[2s+3]_q+q^{-1}[2s-1]_q}{[2s+1]_q} \sum_{i=-s}^s q^{-4i} [s+i]_q[s-i+1]_q\\
&=q^{4s}\,\frac{q[2s+3]_q+q^{-1}[2s-1]_q}{[2s+1]_q(q-q^{-1})^2}\Big((q^{2s+1}+q^{-2s-1})[2s+1]_{q^2} - q^{-1}[2s+1]_q - q[2s+1]_{q^3}\Big),
\end{align*}
where we used in the third last line that $q^m[n]_q+q^{-n}[m]_q=[m+n]_q$ and in the last line \eqref{prop_estimates_generators_eq4}. Since this expression is real-valued, by \eqref{prop_estimates_generators_eq2} it follows that also
\begin{align*}
&d_\half d_s(\varphi_\half\ten\varphi_s)\big((\pi_\half\ten\pi_s)(\hat\Delta(\tilde X_1))^*(1_{\half}\ten\pi_s(\tilde X_1))\big)\\
&= q^{4s}\,\frac{q[2s+3]_q+q^{-1}[2s-1]_q}{[2s+1]_q(q-q^{-1})^2}\Big((q^{2s+1}+q^{-2s-1})[2s+1]_{q^2} - q^{-1}[2s+1]_q - q[2s+1]_{q^3}\Big).
\end{align*}
Adding these four terms shows that $\tilde X_1$ satisfies the required estimated. Indeed, using the asymptotic behaviour of the $q$-numbers (cf.\ Lemma \ref{Lem_q-numbers}) gives
\begin{align*}
&d_\half d_s(\varphi_\half\ten\varphi_s)\big((1_\half\ten\pi_s(\tilde X_1) - (\pi_\half\ten\pi_s)\hat\Delta(\tilde X_1))^*(1_\half\ten\pi_s(\tilde X_1) - (\pi_\half\ten\pi_s)\hat\Delta(\tilde X_1))\big) \displaybreak[1]\\
&= q^{4s}\,\frac{q+q^{-1}}{(q-q^{-1})^2}\Big((q^{2s+1}+q^{-2s-1})[2s+1]_{q^2} - q^{-1}[2s+1]_q - q[2s+1]_{q^3}\Big)\\
&\qquad + \frac{q^{4s-2}}{(q-q^{-1})^2}\Big((q^{2s}+q^{-2s})[2s]_{q^2} - q^{-1}[2s]_q - q[2s]_{q^3}\Big)\\
&\qquad + \frac{q^{4s+2}}{(q-q^{-1})^2}\Big((q^{2s+2}+q^{-2s-2})[2s+2]_{q^2} - q^{-1}[2s+2]_q - q[2s+2]_{q^3}\Big)\\
&\qquad - 2q^{4s}\,\frac{q[2s+3]_q+q^{-1}[2s-1]_q}{[2s+1]_q(q-q^{-1})^2}\Big((q^{2s+1}+q^{-2s-1})[2s+1]_{q^2} - q^{-1}[2s+1]_q - q[2s+1]_{q^3}\Big) \displaybreak[2]\\
&= \frac{q^{4s}}{(q-q^{-1})^2}\Big((q+q^{-1})\big(O(q^{-4s}) + q^{-2s-1}\,\frac{-q^{-4s-2}}{q^2-q^{-2}} - O(q^{-2s}) + O(q^{6s})- q\,\frac{-q^{-6s-3}}{q^3-q^{-3}}\big)\\
&\qquad + q^{-2}\big(O(q^{-4s}) + q^{-2s}\,\frac{-q^{-4s}}{q^2-q^{-2}} - O(q^{-4s}) - O(q^{6s}) - q\,\frac{-q^{-6s}}{q^3-q^{-3}}\big)\\
&\qquad + q^2\big(O(q^{-4s}) + q^{-2s-2}\,\frac{-q^{-4s-4}}{q^2-q^{-2}} - O(q^{-4s}) - O(q^{6s}) - q\,\frac{-q^{-6s-6}}{q^3-q^{-3}}\big)\\
&\qquad - 2\big(qq^{-2}(1+O(q^{4s})) + q^{-1}q^2(1+O(q^{4s}))\big)\\
&\qquad \quad \times\big(O(q^{-4s}) + q^{-2s-1}\,\frac{-q^{-4s-2}}{q^2-q^{-2}} + O(q^{-4s}) + O(q^{6s}) - q\,\frac{-q^{-6s-3}}{q^3-q^{-3}}\big)\Big) \displaybreak[2]\\
&= \frac{q^{4s}}{(q-q^{-1})^2}\Big((q+q^{-1})\big(\frac{-q^{-6s-3}}{q^2-q^{-2}} + \frac{q^{-6s-2}}{q^3-q^{-3}}\big) - \frac{q^{-6s-2}}{q^2-q^{-2}} + \frac{q^{-6s-1}}{q^3-q^{-3}} - \frac{q^{-6s-4}}{q^2-q^{-2}} \\
&\qquad + \frac{q^{-6s-3}}{q^3-q^{-3}} + 2(q+q^{-1})(1+O(q^{4s}))\big(\frac{q^{-6s-3}}{q^2-q^{-2}} - \frac{q^{-6s-2}}{q^3-q^{-3}}\big) + O(q^{-4s})\Big) \displaybreak[2]\\
&= \frac{q^{4s}}{(q-q^{-1})^2}\big(O(q^{-2s}) + O(q^{-4s})\big)\\
&= O(1),
\end{align*}
which proves the claim for $\tilde X_1$.
\end{pf}

\begin{Lem}\label{general_mu}
Let $C>0$, $x\in l^\infty(\widehat{\SUqt})$ and $\mu=\sum_{t\in\half\Int_+} \mu_t\delta_t$ be a probability measure on $\half\Int_+$. Denote $U=\bigoplus_{t\in\supp(\mu)} U_t$. Assume that the element $x$ satisfies the estimate
\begin{equation}\label{general_mu_eq1}
\|1_{\half}\ten\pi_s(x)-(\pi_{\half}\ten\pi_s)\hat\Delta(x)\|_{\varphi_{\half}\ten\varphi_s}^2 \leq C\, d_s^{-1}
\end{equation}
for every $s\in\half\Int_+$, then the inequality
\[
\|1_U\ten\pi_s(x)-(\pi_U\ten\pi_s)\hat\Delta(x)\|_{\varphi_\mu\ten\varphi_s}^2 \leq
C\,d_s^{-1}\frac{[2]_q}{(\sqrt{[2]_q}-\sqrt{2})^2}\sum_{t\in\half\Int_+} \mu_t\, \frac{(d_\half)^{2t}}{d_t}
\]
holds for every $s\in\half\Int_+$.
\end{Lem}
\begin{pf}
Let $x\in l^\infty(\widehat{\SUqt})$ satisfying \eqref{general_mu_eq1}. We start with the following estimate
\begin{align*}
&\|(1_\half\ten\cdots\ten1_\half\ten(\underbrace{\pi_\half\ten\cdots\ten\pi_\half}_k\ten\pi_s)\hat\Delta^k(x)) \\
&\qquad- (1_\half\ten\cdots\ten1_\half\ten1_\half\ten(\underbrace{\pi_{\half}\ten\cdots\ten\pi_{\half}}_{k-1}\ten\pi_s) \hat\Delta^{k-1}(x))\|_{\varphi_{\half}\ten\cdots\ten\varphi_{\half}\ten\varphi_s}^2\\
&= \|(\pi_{\half}\ten\cdots\ten\pi_{\half}\ten\pi_s)\big((\iota\ten\hat\Delta^{k-1})\hat\Delta(x) - (1_\half\ten\hat\Delta^{k-1}(x))\big)\|_{\varphi_{\half}\ten\cdots\ten\varphi_{\half}\ten\varphi_s}^2 \displaybreak[2] \\
&=(\varphi_{\half}\ten\cdots\ten\varphi_{\half}\ten\varphi_s) (\pi_{\half}\ten\cdots\ten\pi_{\half}\ten\pi_s)\\
&\qquad\big(((\iota\ten\hat\Delta^{k-1})\hat\Delta(x)-1_\half\ten\hat\Delta^{k-1}(x))^* ((\iota\ten\hat\Delta^{k-1})\hat\Delta(x)-1_\half\ten\hat\Delta^{k-1}(x))\big) \displaybreak[2]\\
&= \big((\varphi_{\half}\circ\pi_\half)\ten((\varphi_{\half}\ten \cdots\ten\varphi_{\half}\ten\varphi_s) (\pi_{\half}\ten\cdots\ten\pi_{\half}\ten\pi_s)\hat\Delta^{k-1})\big) \\
&\qquad \big((\hat\Delta(x)-1_\half\ten(x))^*(\hat\Delta(x)-1_\half\ten(x))\big) \displaybreak[2]\\
&= \sum_{r\in\half\Int_+} m_{\half^{\ten k-1}\ten s}^r \frac{d_r}{d_{\half}^{k-1}d_s} (\varphi_{\half}\ten\varphi_r)(\pi_{\half}\ten\pi_r)\big((\hat\Delta(x)-1_\half\ten(x))^* (\hat\Delta(x)-1_\half\ten(x))\big)\\
&\leq \sum_{r\in\half\Int_+} m_{\half^{\ten k-1}\ten s}^r \frac{d_r}{d_{\half}^{k-1}d_s}\, C\,d_r^{-1} \displaybreak[2]\\
&= d_{\half}^{-(k-1)}d_s^{-1}\, C \sum_{r\in\half\Int_+} m_{\half^{\ten k-1}\ten s}^r \displaybreak[2]\\
&\leq d_{\half}^{-(k-1)}d_s^{-1}\, C 2^{k-1}\\
&= C\, \Big(\frac{2}{d_{\half}}\Big)^{k-1}d_s^{-1}.
\end{align*}
Here we used \eqref{estimates_multiplicities_suqt_eq1} to obtain a bound on the sum of the multiplicities. Now let $t\in\half\Int_+$. Then the representation $\pi_t$ embeds into $(\pi_{\half})^{\ten 2t}$ with multiplicity 1. Therefore for any positive element $y$ it holds
\[
\frac{d_t}{(d_\half)^{2t}}\,\varphi_t(y)\leq \sum_{r\in\half\Int_+} m_{\half^{\ten 2t}}^r\,\frac{d_r}{(d_\half)^{2t}}\,\varphi_r(y) = (\varphi_\half)^{2t}(y).
\]
Combining this with the previous estimate we obtain
\begin{align}
&\|1_t\ten\pi_s(x)-(\pi_t\ten\pi_s)\hat\Delta(x)\|_{\varphi_t\ten\varphi_s}^2 \notag \\
&\leq \frac{(d_\half)^{2t}}{d_t} \|1_{\half}\ten \cdots\ten 1_{\half}\ten\pi_s(x)-(\pi_{\half}\ten\cdots\ten\pi_{\half}\ten\pi_s) (\hat\Delta^{2t-1}\ten\iota)\hat\Delta(x)\|_{\varphi_{\half}\ten \cdots\ten\varphi_{\half}\ten\varphi_s}^2 \notag \displaybreak[2] \\
&\leq \frac{(d_\half)^{2t}}{d_t}\Big(\|\underbrace{1_{\half}\ten \cdots\ten 1_{\half}}_{2t}\ten\pi_s(x)-\underbrace{1_{\half}\ten \cdots\ten 1_{\half}}_{2t-1}\ten((\pi_{\half}\ten\pi_s) \hat\Delta(x))\|_{\varphi_{\half}\ten\cdots\ten\varphi_{\half}\ten\varphi_s} + \ldots \notag\\
&\quad + \|1_{\half}\ten((\underbrace{\pi_{\half}\ten\cdots\ten\pi_{\half}}_{2t-1}\ten\pi_s)\hat\Delta^{2t-1}(x)) - (\underbrace{\pi_{\half}\ten\cdots\ten\pi_{\half}}_{2t}\ten\pi_s) \hat\Delta^{2t}(x)\|_{\varphi_{\half}\ten\cdots\ten\varphi_{\half}\ten\varphi_s}\Big)^2 \notag \displaybreak[2] \\
&\leq \frac{(d_\half)^{2t}}{d_t} \Big( (C\,d_s^{-1})^{1/2} + \Big(C\,\Big(\frac{2}{d_{\half}}\Big)d_s^{-1}\Big)^{1/2} + \ldots + \Big(C\, \Big(\frac{2}{d_{\half}}\Big)^{2t-1}d_s^{-1}\Big)^{1/2}\Big)^2 \notag \displaybreak[2]\\
&\leq C\,\frac{(d_\half)^{2t}}{d_td_s} \Big(\sum_{n=0}^\infty \Big(\frac{2}{[2]_q}\Big)^{n/2}\Big)^2 \notag \displaybreak[2]\\
&= C\,\frac{(d_\half)^{2t}}{d_td_s} \Big(\frac{1}{1-\sqrt{\frac{2}{[2]_q}}}\Big)^2 \notag\\
&= C\,\frac{(d_\half)^{2t}}{d_td_s} \frac{[2]_q}{(\sqrt{[2]_q}-\sqrt{2})^2}. \label{general_mu_eq2}
\end{align}
Since
\begin{align*}
\|1_U\ten\pi_s(x)-(\pi_U\ten\pi_s)\hat\Delta(x)\|_{\varphi_\mu\ten\varphi_s}^2 &= \sum_t \mu_t \|1_U\ten\pi_s(x)-(\pi_U\ten\pi_s)\hat\Delta(x)\|_{\varphi_t\ten\varphi_s}^2,
\end{align*}
the proof is complete.
\end{pf}

\begin{Cor}\label{cor_general_measure}
Let $C>0$, $x\in l^\infty(\widehat{\SUqt})$ and $\mu=\sum_{t\in\half\Int_+} \mu_t\delta_t$ be a probability measure on $\half\Int_+$. Assume that $\sum_t \mu_t(1+q^2)^{2t}<\infty$ and that the element $x$ satisfies estimate \eqref{general_mu_eq1} for every $s\in\half\Int_+$, then there exists a constant $C'$ independent of $n$ such that
\[
\|j_n(x)-j_{n+1}(x)\|^2_{\varphi_\mu^\infty}\leq C' \Big(\sum_{r\in\half\Int_+} \mu(r)\,\frac{\dim(U_r)}{d_r}\Big)^n.
\]
In particular if $\supp(\mu)\neq\{0\}$, then
\[
\sum_{n=0}^\infty \|j_n(x)-j_{n+1}(x)\|^2_{\varphi_\mu^\infty} <\infty.
\]
\end{Cor}
\begin{pf}
For $n\geq 1$ it holds that $\frac{q-q^{-1}}{q^2-1}\,q^n\leq \frac{1}{[n]_q} \leq q^n(q-q^{-1})$. Also $\sum_t \mu_t[2]_q^{2t}q^{2t} = \sum_t \mu_t(1+q^2)^{2t}$. So $\sum_t \mu_t\frac{(d_\half)^{2t}}{d_t}<\infty$ if and only if $\sum_t \mu_t (1+q^2)^{2t}<\infty$.
Use Lemma \ref{general_mu} and the decomposition of $\varphi_\mu^n$ using the constants $c_{n,s}(\mu)$ (see \eqref{constants_cnr_eq1}) to obtain the following chain of inequalities:
\begin{align*}
\|j_n(x)-j_{n+1}(x)\|^2_{\varphi_\mu^\infty} &= \|1\ten\hat\Delta^{n-1}(x)-\hat\Delta^{n}(x)\|^2_{\varphi_\mu^{\ten n+1}}\\
&= \|1\ten x-\hat\Delta(x)\|^2_{\varphi_\mu\ten(\varphi_\mu)^n} \displaybreak[2]\\
&= \sum_{s\in \half\Int_+} c_{n,s}(\mu) \|1\ten x-\hat\Delta(x)\|^2_{\varphi_\mu\ten\varphi_s} \displaybreak[2]\\
&\leq \sum_{s\in \half\Int_+} c_{n,s}(\mu) C\,d_s^{-1}\frac{[2]_q}{(\sqrt{[2]_q}-\sqrt{2})^2} \Big(\sum_{t\in\half\Int_+} \mu_t\, \frac{(d_\half)^{2t}}{d_t}\Big)\\
&\leq C\,\frac{[2]_q}{(\sqrt{[2]_q}-\sqrt{2})^2} \Big(\sum_{t\in\half\Int_+} \mu_t\, \frac{(d_\half)^{2t}}{d_t}\Big)\Big(\sum_{r\in\half\Int_+} \mu(r)\,\frac{\dim(U_r)}{d_r}\Big)^n.
\end{align*}
This proves the first part of the corollary. Write $d:= \sum_r \mu(r)\,\frac{\dim(U_r)}{d_r}$. If $\supp(\mu)\neq\{0\}$ the sum satisfies $0<d<1$ and thus $\sum_{n=0}^\infty d^n = \frac{1}{1-d}<\infty$.
\end{pf}

Now it is a matter of putting everything together to prove the regularity for random walks on $\widehat{\SUqt}$.

\medskip
\begin{pf}[ of the first part of Theorem \ref{Thm_convergence_SUq2}]
Proposition \ref{prop_estimates_generators} and Corollary \ref{cor_general_measure} imply that for elements $x=\tilde X_i$, $i=-1,0,1$ the following holds
\[
\sum_{n=0}^\infty (\cdots\ten\varphi_\mu\ten\varphi_\mu)\big((j_n(x)-j_{n+1}(x))^*(j_n(x)-j_{n+1}(x))\big) <\infty.
\]
Hence $\textrm{s-}\lim_n j_n(x)$ exists. As $X_0^*=X_0$ and $\tilde X_{-1}^*=q\tilde X_1$, we in fact have that $\sstar\lim_n j_n(x)$ exists. In other words, the elements $\tilde X_j$ are $\varphi_\mu$-regular. By Lemma \ref{Lem_span_generators} and Theorem \ref{Thm_Martin_SUqt} the Martin boundary $\tilde M(\widehat\SUqt,\mu)$  is generated as a \Cstar algebra by $c_0(\widehat{\SUqt})$ and $\tilde X_i$, $i=-1,0,1$. Proposition \ref{prop_reg_algebra} asserts that all elements in $c_0(\widehat\SUqt)$ are $\varphi_\mu$-regular and that the set of regular elements forms a \Cstar algebra. Hence all elements in $\tilde M(\widehat{\SUqt},\mu)$ are regular.
\end{pf}

\begin{Rem}
This proof of regularity of ``functions'' on the Martin compactification is fundamentally different from the proof in the classical case. There is no need to use some sort of ``noncommutative stopping times'' (it is even not clear what such objects should be). Due to the fact that $0<q<1$, there is exponentially fast convergence to the boundary if the measure $\mu$ is nice enough, this is much stronger convergence than classically. We emphasize that this proof does not work in the classical case, so not for $q=1$.
\end{Rem}

\section{Boundary convergence and monoidal equivalence}\label{sec_monoidal_equivalence}
In this section we will establish that the property of boundary convergence is stable under monoidal equivalence. In \cite{DeRijdtVanderVennet10} the authors showed how one can compute the Poisson and Martin boundary via monoidal equivalences. We will follow their approach. However, later we will tackle the same problem in a more general setting by defining boundary convergence for random walks on \Cstar tensor categories and showing compatibility with the definition for quantum groups (see Subsections \ref{subsec_categorical_boundary_convergence} and \ref{subsec_categorical_correspondence} below). For this reason we will be very brief in this section and only outline the main arguments.

\begin{Def}[\protect{\cite[Def.\ 3.1]{BichonDeRijdtVaes06}}]
Two compact quantum groups $G_i:=(C(G_i),\Delta_i)$ $(i=1,2)$ are said to be {\it monoidally equivalent} if the representation categories $\Rep(G_1)$ and $\Rep(G_2)$ are unitarily monoidally equivalent as \Cstar tensor categories.
\end{Def}

The following theorem is proven in \cite{BichonDeRijdtVaes06} the formulation is from \cite[Thm.\ 4.2]{DeRijdtVanderVennet10}. Since this result plays a fundamental role and we use different conventions, we recall it here for convenience.

\begin{Thm}[\protect{\cite{BichonDeRijdtVaes06}}]\label{Thm_link_algebra}
Given two monoidally equivalent compact quantum groups $G_1$ and $G_2$ together with a unitary monoidal equivalence $\kappa\colon \Rep(G_1)\ra \Rep(G_2)$, the following holds:
\begin{enumeraterm}
\item there exists a (up to isomorphism) unique unital $*$-algebra $\Bcal$ equipped with a faithful state $\omega$ and unitary elements $X^s\in \Bcal\ten B(\Hcal_s,\Hcal_{\kappa(s)})$ satisfying
\begin{itemize}
\item $X^s_{12}X^t_{13}(\iota_\Bcal\ten S) = (\iota_\Bcal\ten\kappa(S))X^r$, whenever $S\in\Hom(r,s\ten t)$;
\item the matrix coefficients of $\{X^s\}_{s\in\Irr(G_1)}$ form a basis in $\Bcal$ as a vector space;
\item $(\omega\ten\iota)(X^s)=0$ if $s\neq 0$;
\end{itemize}
\item there exist unique commuting ergodic actions $\delta_1\colon\Bcal\ra\Bcal\ten_{\alg}\Com[G_1]$ and $\delta_2\colon\Bcal\ra \Com[G_2]\ten_{\alg}\Bcal$ satisfying for every $s\in\Irr(G_1)$
\begin{align*}
&(\delta_1\ten\iota)(X^s) = (X^s_{13})(U_s)_{23} \in \Bcal\ten \Com[G_1]\ten B(\Hcal_s,\Hcal_{\kappa(s)});\\
&(\delta_2\ten\iota)(X^s) = (U_{\kappa(s)})_{13}(X^s_{23}) \in \Com[G_2]\ten\Bcal\ten B(\Hcal_s,\Hcal_{\kappa(s)});
\end{align*}
\item the state $\omega$ is invariant under $\delta_1$ and $\delta_2$ and is given by $\omega(b)1_B=(\iota\ten h)\delta_1(b)$.
\end{enumeraterm}
\end{Thm}

This algebra $\Bcal$ is called the {\it link algebra} of $G_1$ and $G_2$ under the monoidal equivalence $\kappa$. The algebra $\Bcal$ can be constructed explicitly, however such a definition is irrelevant for our purposes. We will only work with the properties described above. Using the faithful state $\omega$ this can all be extended to the von Neumann framework. Let $B$ be the von Neumann algebra generated by $\Bcal$ in the GNS representation of $\omega$. The actions $\delta_1$ and $\delta_2$ have unitary implementations and can therefore be extended to $\delta_1\colon B\ra B\bar\ten L^\infty(G_1)$ and $\delta_2\colon B\ra L^\infty(G_2)\bar\ten B$ (see \cite[Rem.\ 3.12]{DeRijdtVanderVennet10}).

\begin{Not}\label{Not_maps_kn}
In this section $X$ denotes the following element\footnote{Later the symbol $X$ will again be used to indicate objects in categories.}
\[
X:= \prod_{s\in\Irr(G_1)} X^s\in \bds_{s\in\Irr(G_1)} B\bar\ten B(\Hcal_s,\Hcal_{\kappa(s)}).
\]
Use this to define the collection of maps
\begin{align*}
k_n\colon\bigotimes_{i=-n}^{-1} l^\infty(\hat G_2) &\ra B\bar\ten\Big(\bigotimes_{i=-n}^{-1} l^\infty(\hat G_1)\Big), \qquad (n\geq1),\\
x &\mapsto (X_{1,n+1}^*\cdots X_{1,2}^*)(1_B\ten x)(X_{1,2}\cdots X_{1,n+1});\\
k_0\colon\Com\ra B,\qquad z&\mapsto z1_B.
\end{align*}
\end{Not}

Observe that if $W^{(i)}$ is the multiplicative unitary of $G_i$, then $(\delta_1\ten\iota)(X) = X_{13}W_{23}^{(1)}$ and $(\iota\ten\delta_2)(X)= W_{13}^{(2)}X_{23}$ when viewed in the multiplier algebra.

\bigskip
For the remainder of this section we assume that $G_1$ and $G_2$ are two monoidally equivalent compact quantum groups with unitary monoidal equivalence $\kappa\colon \Rep(G_1)\ra \Rep(G_2)$.

\begin{Not}\label{Not_box_product}
If $\alpha\colon N\ra L^\infty(G_1)\bar\ten N$ is a left action of the compact quantum group $G_1$ on a von Neumann algebra $N$, denote
\[
B\boxtimes^\alpha N :=\{x\in B\bar\ten N \,:\, (\delta_1\ten\iota)(x)=(\iota\ten\alpha)(x)\}.
\]
In the algebraic setting, if $\Ncal$ is a unital $*$-algebra and $\alpha\colon\Ncal\ra \Com[G_1]\ten_{\alg}\Ncal$ an action, let
\[
\Bcal\boxtimes_{\alg}^\alpha\Ncal:= \{x\in \Bcal\ten_{\alg} \Ncal\,:\, (\delta_1\ten\iota)(x)=(\iota\ten\alpha)(x)\}.
\]
The \Cstar algebraic case is more involved due to technicalities with the range, see also \cite[\textsection 8]{Voigt11}. We formulate it in general, but we will only need it if the \Cstar algebra under consideration equals the Martin boundary or Martin compactification.
Suppose $D$ is a \Cstar algebra with an action $\alpha\colon D\ra C(G_1)\ten D$. Define for $s\in\Irr(G_1)$ the spectral subspaces
\[
K_s := \{x\in D\ten\bar\Hcal_s\,:\, (\alpha\ten\iota)(x)=x_{13}U_{23}^s\}.
\]
Denote
\[
\Dcal_s:=\Span\{x(1\ten\xi)\,:\, x\in K_s,\, \xi\in\Hcal_s\}, \qquad \qquad \Dcal:=\Span\{\Dcal_s\,:\, s\in\Irr(G_1)\}.
\]
It follows that $\alpha\colon\Dcal_s\ra \Com[G_1]_s\ten_{\alg}\Dcal_s$ is a Hopf algebra coaction. An action $\alpha$ is called {\it reduced} if the conditional expectation onto the fixed point algebra $(\iota\ten h)\alpha\colon D\ra D^{\alpha}$ is faithful. For a reduced action $\alpha\colon D\ra C(G_1)\ten D$ define $B\boxtimes_{\textrm{red}}^\alpha D$ to be the norm closure of $\Bcal\boxtimes_{\alg}^\alpha \Dcal$ in $B\ten D$.
\end{Not}

Proposition \ref{Prop_Properties_link_map} below will play a crucial role to transport the convergence properties of one quantum group to a monoidal equivalent one. The result shows that all the algebraic operations can be transferred from $G_1$ to $G_2$ by means of the link algebra $B$. The construction is motivated by and can be proved in a spirit similar to \cite[\textsection~7 and \textsection~8]{DeRijdtVanderVennet10}, it extends their results to higher tensor powers. When defining random walks De Rijdt and Vander Vennet work with states of the form $\psi_s :=d_s^{-1}\Tr(\pi_s(\cdot\rho))$ and slice in the right leg, while we work with $\varphi_s:=d_s^{-1}\Tr(\pi_s(\cdot\rho^{-1}))$ slicing in the left leg. Moreover they interchange the roles of $G_1$ and $G_2$ in the monoidal equivalence, while we do not. Therefore our maps defining the isomorphisms have a slightly different form than theirs.

\begin{Prop}\label{Prop_Properties_link_map}
Let $\alpha_l$ be the adjoint action as defined by \eqref{extended_adjoint_action_eq1}. Using the notation introduced above, the following holds for $0\leq m\leq n-1$:
\begin{enumeraterm}
\item $k_n\colon \bigotimes_{-n}^{-1} l^\infty(\hat{G}_2)\ra B\boxtimes^{\alpha_l}\big(\bigotimes_{-n}^{-1} l^\infty(\hat{G}_1)\big)$ is a $*$-isomorphism;
\item $k_{n+1}\circ(\iota^{\ten m}\ten\hat{\Delta}_2\ten\iota^{\ten n-m-1}) = (\iota_B\ten\iota^{\ten m}\ten\hat{\Delta}_1\ten\iota^{\ten n-m-1})\circ k_n$;
\item $k_{n-1}\circ(\iota^{\ten m}\ten\varphi_{\kappa(s)}\ten\iota^{\ten n-m-1}) = (\iota_B\ten\iota^{\ten m}\ten\varphi_{s}\ten\iota^{\ten n-m-1})\circ k_n$;
\item $k_{n+1}(1\ten x) = (k_n(x))_1\ten1_{l^\infty(G_1)}\ten (k_n(x))_{2,\ldots,n} = (\iota_B\ten 1_{l^\infty(\hat G_1)}\ten\iota^{\ten n})(k_n(x))$, if $x\in \bigotimes_{-n}^{-1} l^\infty(\hat{G}_2)$;
\end{enumeraterm}
\end{Prop}

Suppose that $\mu$ is a probability measure on $\Irr(G_1)$. Recall that the pushforward measure $\kappa_*(\mu)$ is defined by $\kappa_*(\mu)(s):=\mu(\kappa^{-1}(s))$. It is a probability measure on $\Irr(G_2)$ and satisfies $\kappa(\bar s)=\overline{\kappa(s)}$ and thus $\kappa_*(\bar\mu)=\overline{\kappa_*(\mu)}$.

\begin{Thm}\label{Thm_boundary_conv_monoidal_equiv}
Suppose that $G_1$ and $G_2$ are two monoidally equivalent quantum groups with a unitary monoidal equivalence $\kappa\colon \Rep(G_1)\ra \Rep(G_2)$ and that $\mu$ is a probability measure on $\Irr(G_1)$. If the random walk on $\hat G_1$ defined by $\mu$ converges to the boundary, then so does the random walk defined by $\kappa_*(\mu)$ on $\hat G_2$.
\end{Thm}
\begin{sketch}
Suppose $x\in l^\infty(\hat G_2)$. From Proposition \ref{Prop_Properties_link_map} it is easy to compute that
\begin{align*}
&(\cdots\ten\varphi_{\kappa_*(\mu)}\ten\varphi_{\kappa_*(\mu)}) \big((j_n(x)-j_m(x))^*(j_n(x)-j_m(x))\big) \notag\\
&= (\omega\ten(\cdots\ten\varphi_\mu\ten\varphi_\mu)) \big(((\iota_B\ten j_n)k_1(x)-(\iota_B\ten j_m)k_1(x))^*((\iota_B\ten j_n)k_1(x)-(\iota_B\ten j_m)k_1(x))\big).
\end{align*}
If $y\in \Bcal\boxtimes_{\textrm{\alg}}^{\alpha_l}\tilde \Mcal(\hat G_1,\mu)$, then $y$ is of the form $y=\sum_{i=1}^n b_i\ten y_i$ for some finite sum and elements $y_i\in \tilde\Mcal(\hat G_1,\mu)$ and $b_i\in\Bcal$. Using convergence to the boundary for $\hat G_1$ one can show that $\sstar\lim_n (\iota_B\ten j_n)(y)$ exists.  Using density this extends to all elements $y\in B\boxtimes_{\textrm{reg}}^{\alpha_l}\tilde M(\hat G_1,\mu)$. Since $k_1\colon \tilde M(\hat G_2,\kappa_*(\mu)) \ra B\boxtimes^{\alpha_l}_{\textrm{red}} \tilde M(\hat G_1,\mu)$ is a $*$-isomorphism (cf.\ \cite[Thm.\ 10.1]{DeRijdtVanderVennet10}) $\sstar\lim_n j_n(x)$ exists for all $x\in\tilde M(\hat G_2,\kappa_*(\mu))$.

To prove the representation of harmonic elements first show that $k_1\colon H^\infty(\hat G_2,\kappa_*(\mu)) \ra B \boxtimes^{\alpha_l} H^\infty(\hat G_1,\mu)$ is a $*$-isomorphism. Then using the identities of Proposition \ref{Prop_Properties_link_map} and the spectral subspaces of $H^\infty(\hat G_2,\kappa(\mu))$ it is not hard to show that the following limits exists and equalities hold
\begin{align*}
\nu_{G_2}(K_{\kappa_*(\bar\mu)}(x)h) &= \lim_n(\varphi_{\kappa_*(\mu)}\ten \ldots\ten \varphi_{\kappa_*(\mu)})\big(\hat\Delta^n(K_{\kappa_*(\bar\mu)}(x)h)\big) \\
&= \lim_n(\omega\ten\varphi_\mu\ten \ldots\ten \varphi_\mu) \big((\iota_B\ten\hat\Delta^n)(k_1(K_{\kappa_*(\bar\mu)}(x)h))\big) \displaybreak[2] \\
&= (\omega\ten\nu_{G_1})\big(k_1(K_{\kappa_*(\bar\mu)}(x)h)\big)\\
&= (\omega\ten\hat\psi_{G_1})\big(k_1(x)k_1(h)\big) = \hat\psi_{G_2}(xh),
\end{align*}
which is the required identity.
\end{sketch}

\begin{Exam}
In \cite{VanDaeleWang96} the free orthogonal quantum groups were introduced. Let $n \geq 2$ and $F\in \textrm{GL}_n(\Com)$ be a matrix with $F\bar F=\pm1$, where $\bar F$ is the matrix $(\bar F)_{ij} = \bar F_{ij}$. Define $A_0(F)$ to be the universal \Cstar algebra generated by elements $u_{ij}$, $1\leq i,j,\leq n$ satisfying:
\[
U=(u_{ij})_{i,j} \textrm{ is unitary}, \qquad U=FU^cF^{-1},
\]
where $(U^c)_{i,j}=u_{ij}^*$. The comultiplication is given by $\Delta(u_{ij}):=\sum_{k=1}^n u_{ik}\ten u_{kj}$. Note that for
\[
F_q=\begin{pmatrix} 0 & -\sqrt q \\ {\sqrt{q}}^{-1} &0\end{pmatrix}
\]
it holds $A_0(F_q)=\SUqt$. Also $(A_0(F_1),\Delta_1)$ is monoidally equivalent to $(A_0(F_2),\Delta_2)$ if and only if $\sgn(F_1\bar F_1) = \sgn(F_2\bar F_2)$ and $\Tr(F_1^* F_1)=\Tr(F_2^*F_2)$ (see \cite[Thm.\ 5.3, Cor.\ 5.4]{BichonDeRijdtVaes06}). It follows that every free orthogonal quantum group $A_0(F)$ with $F\in \textrm{GL}_n(\Com)$ and $n\geq 2$ is monoidally equivalent to $\SUqt$ for some $q\in[-1,1]\setminus\{0\}$.

Let $A_0(F)$ be a free orthogonal quantum group monoidally equivalent to $\SUqt$ for some $q\in(0,1)$ (thus $F\bar F=-1$ and $\Tr(F^*F)>2$), so $q$ satisfies $q+q^{-1}=\Tr(F^*F)$. Suppose we have an explicit identification $\kappa\colon\Irr(A_0(F))\ra \half\Int_+\cong \Irr(\SUqt)$ (see also \cite[Thm.\ 1]{Banica96}) and assume that a generating probability measure $\mu$ on $\Irr(A_0(F))$ satisfies the condition
\[
\sum_{s\in\Irr(A_0(F))} \mu(s) (1+q^2)^{2\kappa(s)}<\infty,
\]
then by Theorem \ref{Thm_convergence_SUq2} the random walk defined by $\kappa_*(\mu)$ on $\widehat{\SUqt}$ converges to the boundary. Hence Theorem \ref{Thm_boundary_conv_monoidal_equiv} shows that the random walk on $\widehat{A_0(F)}$ defined by $\mu$ converges to the boundary.
\end{Exam}

\section{A categorical approach to the Martin boundary}\label{sec_categorical_Martin_boundary}
The results of \cite{DeRijdtVanderVennet10} and Section \ref{sec_monoidal_equivalence} combined with the paper \cite{NeshveyevYamashita14c} indicate that the Martin boundary and compactification and also the notion of convergence to the boundary can be defined on the level of \Cstar tensor categories. In this section we show that this is indeed the case.

\subsection{Random walks on \texorpdfstring{\Cstar}{C*-}tensor categories}\label{sec_random_walks on_categories}
Neshveyev and Yamashita found a way to define random walks on \Cstar tensor categories \cite{NeshveyevYamashita14a}. For such random walks they defined a Poisson boundary and used it to prove their characterisation of quantum groups with the same representation theory as $\textrm{SU}(N)$. Motivated by their work we will construct a Martin boundary and formulate convergence to the boundary for such random walks. This subsection will form the starting point of this theory. We introduce the analogue of the algebra of functions on the path space and the necessary functors on these spaces. It is motivated by \cite[\textsection~3.1]{NeshveyevYamashita14a}.

\begin{Not}\label{Not_cat_path_spaces}
Let $\Ccal$ be a strict \Cstar tensor category with simple unit $\unit$. For $n\geq 1$ and an object $U\in\Ob(\Ccal)$ define the $n$-ary functor
\begin{align*}
(\iota^{\ten n}\ten U)\colon \Ccal\times\cdots\times\Ccal&\ra \Ccal,\\
(\iota^{\ten n}\ten U)(X_1,\ldots,X_n)&:= (X_1\ten\cdots\ten X_n\ten U), \\
(\iota^{\ten n}\ten U)(f_1,\ldots,f_n)&:=f_1\ten\cdots\ten f_n\ten\iota_U.
\end{align*}
Clearly $\iota^{\ten n}\ten U$ is unitary. For two objects $U,V\in\Ob(\Ccal)$ consider the space of natural transformations $\nat(\iota^{\ten n}\ten U, \iota^{\ten n}\ten V)$. Since every object can be decomposed into simple ones and $\iota^{\ten n}\ten U$ respects direct sums, a natural transformation $\eta\colon \iota^{\ten n}\ten U\ra \iota^{\ten n}\ten V$ is completely determined by its action on simple objects $U_s$. Thus
\begin{equation}\label{Not_cat_path_spaces_eq3}
\nat(\iota^{\ten n}\ten U, \iota^{\ten n}\ten V)\cong \prod_{s_1,\ldots,s_n\in\Irr(\Ccal)} \Hom_{\Ccal}(U_{s_1}\ten\cdots\ten U_{s_n}\ten U, U_{s_1}\ten\cdots\ten U_{s_n}\ten V).
\end{equation}
A natural transformation $\eta\in\nat(\iota^{\ten n}\ten U, \iota^{\ten n}\ten V)$ is {\it (uniformly) bounded} if
\[
\|\eta\|_\infty :=\sup\{\|\eta_{U_{s_1},\ldots,U_{s_n}}\|\,:\,s_1,\ldots,s_n\in\Irr(\Ccal)\}<\infty.
\]
Denote by $\nat_b(\iota^{\ten n}\ten U,\iota^{\ten n}\ten V)$ the uniformly bounded natural transformations. Define $\Ccal_{-n}$ as the \Cstar category which is the subobject completion of the category with objects $\Ob(\Ccal)$ and morphisms
\[
\Hom_{\Ccal_{-n}}(U,V) := \nat_b(\iota^{\ten n}\ten U,\iota^{\ten n}\ten V), \qquad (U, V\in\Ob(\Ccal)).
\]
So a morphism $\eta=(\eta_{X_1,\ldots, X_n})_{X_1,\ldots,X_n}\in \Hom_{\Ccal_{-n}}(U,V) $ consists of a bounded collection of morphisms
\[
\eta_{X_1,\ldots,X_n}\colon X_1\ten\cdots\ten X_n\ten U\ra X_1\ten\cdots\ten X_n\ten V
\]
natural in $X_1,\ldots, X_n\in\Ob(\Ccal)$. In $\Ccal_{-n}$ the multiplication of morphisms is given by composition and the involution is coming from the $*$ on $\Ccal$. The \Cstar norm on $\Hom_{\Ccal_{-n}}(U,V)$ is given by $\|\eta\|_\infty:= \sup \{\|\eta_{s_1,\ldots,s_n}\|\,:\, s_1,\ldots, s_n\in\Irr(\Ccal)\}$ for $\eta\in\nat_b(\iota^{\ten n}\ten U, \iota^{\ten n}\ten V)$.\\
Observe that there exists a canonical functor $\Ecal\colon\Ccal\ra\Ccal_{-n}$ given on objects by $U\mapsto U$ and on morphisms
\begin{equation}\label{Not_cat_path_spaces_eq4}
\Ecal(T):=\iota^{\ten n}\ten T, \qquad \textrm{where } (\iota^{\ten n}\ten T)_{X_1,\ldots, X_n}:=\iota_{X_1}\ten\cdots\ten\iota_{X_n}\ten T.
\end{equation}
Obviously $\Ecal$ is a unitary functor.
\end{Not}

The category $\Ccal_{-n}$ should be thought of as the space of functions on paths of length $n$ (see also Corollary \ref{Cor_equiv_path_spaces}). In particular for $n=1$ we obtain the functions on the space.

\begin{Not}
For $n=1$ there is more structure present. Define the tensor product of objects in $\Ccal_{-1}$ to be inherited from $\Ccal$. For natural transformations define
\begin{equation}\label{Def_tensor_product_nat_transf_eq0}
(\nu\ten\eta):= (\nu\ten\iota_X)(\iota_U\ten\eta), \qquad (\nu\in\Hom_{\Ccal_{-1}}(U,V),\, \eta\in\Hom_{\Ccal_{-1}}(W,X)),
\end{equation}
where
\begin{align}
(\nu\ten\iota_X)_Y&:=\nu_Y\ten\iota_X\colon Y\ten U\ten X \ra Y\ten V\ten X;  \label{Def_tensor_product_nat_transf_eq1}\\
(\iota_U\ten\eta)_Y&:=\eta_{Y\ten U}\colon Y\ten U\ten W \ra Y\ten U\ten X. \label{Def_tensor_product_nat_transf_eq2}
\end{align}
Extend this to the subobject completion. With this tensor structure $\Ccal_{-1}$ becomes a \Cstar tensor category. For $\eta\in\nat(\iota\ten U,\iota\ten V)$ write $\supp(\eta):=\{s\in\Irr(\Ccal)\,:\, \eta_s\neq 0\}$ for the {\it support} of $\eta$. A natural transformation $\eta\in \nat(\iota\ten V, \iota\ten W)$ is called a
\begin{enumerate}
\item {\it natural transformation vanishing at infinity} if $\{\|\eta_s\|\}_{s\in\Irr(\Ccal)}\in c_0(\Irr(\Ccal))$;
\item {\it compactly} or {\it finitely supported natural transformation} if $|\supp(\eta)|<\infty$.
\end{enumerate}
These classes will be denoted by $\nat_0(\iota\ten V, \iota\ten W)$ and respectively $\nat_{00}(\iota\ten V, \iota\ten W)$. By the same argument as for \eqref{Not_cat_path_spaces_eq3} we get
\begin{align*}
&\nat_0(\iota\ten V, \iota\ten W) \cong \vds_{s\in\Irr(\Ccal)} \Hom_{\Ccal}(U_s\ten V, U_s\ten W); \displaybreak[2]\\
&\nat_{00}(\iota\ten V, \iota\ten W) \cong \ds_{s\in\Irr(\Ccal)} \Hom_{\Ccal}(U_s\ten V, U_s\ten W).
\end{align*}
Observe that $M(\nat_0(\iota\ten V, \iota\ten V))=\nat_b(\iota\ten V, \iota\ten V)$, where $M$ indicates the multiplier algebra.
\end{Not}

\begin{Not}\label{Def_cat_maps}
For $0\leq m<n$ define the functors $(\iota\ten\cdot)\colon\Ccal_{-n}\ra \Ccal_{-(n+1)}$ and $\iota^{\ten m}\ten\hat\Delta\ten\iota^{\ten n-m-1}\colon\Ccal_{-n}\ra \Ccal_{-(n+1)}$.  On objects $U\in\Ob(\Ccal)\subset\Ob(\Ccal_{-n})$ they are given by the identity while on morphisms $\eta\in\Hom_{\Ccal_{-n}}(U,V)$ they are defined by
\begin{align}
(\iota\ten\eta)_{X_1,\ldots,X_{n+1}}:=&\iota_{X_1}\ten\eta_{X_2,\ldots,X_{n+1}}; \label{Def_cat_maps_eq1}\\
((\iota^{\ten m}\ten\hat\Delta\ten\iota^{\ten n-m-1})(\eta))_{X_1,\ldots,X_{n+1}}:= &\eta_{X_1,\ldots,X_m,X_{m+1}\ten X_{m+2}, X_{m+3},\ldots, X_{n+1}}. \label{Def_cat_maps_eq2}
\end{align}
These functors extend uniquely to the subobject completions. Given objects $U,\,V,\,X\in\Ob(\Ccal)$ pick a standard solution $(R_X,\bar R_X)$ of the conjugate equations for X. Define the normalised partial trace
\begin{align*}
&\,\tr_X\ten\iota_V\colon\Hom_\Ccal(X\ten U,X\ten V)\ra\Hom_\Ccal(U,V) \\
&(\tr_X\ten\iota_V)(T):= d_X^{-1}(R_X^*\ten\iota_V)(\iota_{\bar X}\ten T)(R_X\ten\iota_U).
\end{align*}
Notice that $\Hom_\Ccal(X\ten U,X\ten V)$ is not a tensor product, so this notation might seem a bit misleading at first. Use these partial traces to define for each $X\in\Ob(\Ccal)$ and $n\geq 2$ the collection of completely positive linear maps
\begin{align}
\tr_X\ten\iota^{\ten n-1}&\colon\Hom_{\Ccal_{-n}}(U,V)\ra\Hom_{\Ccal_{-(n-1)}}(U,V); \notag\\
((\tr_X\ten\iota^{\ten n-1})(\eta))_{X_1,\ldots,X_{n-1}}:=& (\tr_X\ten\iota^{\ten n-1}\ten\iota_V)(\eta_{X,X_1,\ldots,X_{n-1}}) \label{Def_cat_maps_eq3}\\
=& d_X^{-1}(R_X^*\ten\iota^{n-1}\ten\iota_V)(\iota_{\bar X}\ten\eta_{X,X_1,\ldots,X_{n-1}})(R_X\ten\iota^{n-1}\ten\iota_U),\notag
\end{align}
where $(R_X, \bar R_X)$ is a standard solution of the conjugate equations for $X$. Note that $(\tr_X\ten\iota^{\ten n-1})$ does not define a functor, because it does not preserve composition of morphisms. Similarly if $n=1$ and $\eta\in\Ccal_{-1}(U,V)$ denote
\begin{equation}\label{Def_cat_maps_eq4}
\tr_X(\eta):=(\tr_X\ten\iota_V)(\eta_X)\in\Hom_\Ccal(U,V)
\end{equation}
Also the maps $(\tr_X\ten\iota^{\ten n-1})$ can be extended to the subobject completion. The lemma below shows that \eqref{Def_cat_maps_eq1}--\eqref{Def_cat_maps_eq4} are well-defined. \\
Since $\Ccal$ is assumed to be strict, the functor $\hat\Delta$ is coassociative. Define inductively $\hat\Delta^n:=(\hat\Delta\ten\iota^{\ten n-1})\hat\Delta^{n-1}$. By linearity define for a probability measure $\mu$ on $\Irr(\Ccal)$
\[
\tr_\mu\ten\iota^{\ten n-1}:= \sum_{s\in\Irr(\Ccal)} \mu(s) (\tr_{U_s}\ten\iota^{\ten n-1}).
\]
Note that simplicity of the unit $\unit$ ensures that the traces $\tr_\mu$ take values in $\Com\cong\End_\Ccal(\unit)$.
\end{Not}

\begin{Lem}\label{Lem_functor_well-defined}
If $\eta\in\Hom_{\Ccal_{-n}}(U,V)$, then $(\iota\ten\eta)$ and $\iota^{\ten m}\ten\hat\Delta\ten\iota^{\ten n-m-1}(\eta)$ are in $\Hom_{\Ccal_{-(n+1)}}(U,V)$. Moreover if $\mu$ is a probability measure on $\Irr(\Ccal)$, the series
\[
\sum_{s\in\Irr(\Ccal)} \mu(s) (\tr_{U_s}\ten\iota^{\ten n-1}) (\eta_{X_1,\ldots, X_{n-1}})
\]
converges in norm in $\Hom_\Ccal(X_1\ten\cdots\ten X_{n-1}\ten U, X_1\ten\cdots\ten X_{n-1}\ten V)$ and $(\tr_\mu\ten\iota^{\ten n-1}(\eta))\in\Hom \Ccal_{-(n-1)}(U,V)$.
\end{Lem}
\begin{pf}
The proof is straightforward. We will only show that $\hat\Delta$ defines a natural transformation. The other two are similar. Let $f_i\colon X_i\ra Y_i$, for some objects $X_i,\,Y_i\in\Ob(\Ccal)$ and write $k:=n-m-1$. Obviously $f_{m+1}\ten f_{m+2}\colon X_{m+1}\ten X_{m+2}\ra Y_{m+1}\ten Y_{m+2}$, so that
\begin{align*}
((\iota^{\ten m}&\ten \hat\Delta\ten\iota^{\ten k})(\eta))_{Y_1,\ldots,Y_{n+1}}(f_1\ten\cdots\ten f_{n+1}\ten\iota_U)\\
&= \eta_{Y_1,\ldots, Y_{m+1}\ten Y_{m+2},\ldots, Y_{n+1}} (f_1\ten \cdots\ten (f_{m+1}\ten f_{m+2})\ten\cdots\ten f_{n+1}\ten\iota_U) \displaybreak[2]\\
&=(f_1\ten \cdots\ten (f_{m+1}\ten f_{m+2})\ten\cdots\ten f_{n+1}\ten\iota_V)\eta_{X_1,\ldots, X_{m+1}\ten X_{m+2},\ldots, X_{n+1}}\\
&= (f_1\ten\cdots\ten f_{n+1}\ten\iota_V)((\iota^{\ten m}\ten \hat\Delta\ten\iota^{\ten k})(\eta))_{X_1,\ldots,X_{n+1}}.
\end{align*}
Clearly $(\iota\ten\eta)$ and $\iota^m\ten\hat\Delta\ten\iota^{n-m-1}(\eta)$ are bounded by $\|\eta\|_\infty$. While for $\tr_\mu$ we get
\[
\Big\|\sum_s \mu(s) (\tr_{U_s}\ten\iota^{\ten n-1}) (\eta_{X_1,\ldots, X_{n-1}})\Big\| \leq \sum_s \mu(s) d_s^{-1} \|R_s\|^2 \|\iota_{\bar U_s}\ten\eta_{U_s, X_1,\ldots,X_{n-1}}\| \leq \|\eta\|_\infty,
\]
so the series converges and defines a natural bounded transformation.
\end{pf}

\begin{Def}[\protect{\cite{NeshveyevYamashita14a}}]\label{Def_categorical_random_walk}
Suppose that $\mu$ is a probability measure on $\Irr(\Ccal)$. Define the collection of positive linear maps
\[
P_\mu\colon \Hom_{\Ccal_{-1}}(U,V)\ra \Hom_{\Ccal_{-1}}(U,V), \qquad P_\mu:=(\tr_\mu\ten\iota)\hat\Delta.
\]
$P_\mu$ is called the {\it Markov operator} of $\mu$. $P_\mu$ is well-defined due to Lemma \ref{Lem_functor_well-defined} above. Explicitly we have
\[
P_\mu(\eta)_X = \sum_{s\in\Irr(\Ccal)}\mu(s)(\tr_{U_s}\ten\iota_X\ten\iota_V)(\eta_{U_s\ten X}),\quad (\eta\in\nat_b(\iota\ten U,\iota\ten V)\subset\Hom_{\Ccal_{-1}}(U,V)).
\]
\end{Def}

As in the quantum case we need a classical discrete Markov chain on the set of irreducible objects. To define this random walk we need a couple of central natural transformations. At the moment only $V=\unit$ is relevant, but later we will also need general objects $V$.

\begin{Def}
For $V\in\Ob(\Ccal)$ and $t\in\Irr(\Ccal)$ put
\begin{align*}
\kappa^{t,V}\in \Hom_{\Ccal_{-1}}(V,V) \cong \bds_{s\in\Irr(\Ccal)} \End_\Ccal(U_s\ten V), \qquad \kappa^{t,V}_s:=\begin{cases} \iota_{U_t\ten V} &\textrm{if } s=t;\\
0 &\textrm{if } s\neq t.
\end{cases}
\end{align*}
From \eqref{Not_cat_path_spaces_eq3} we see $\nat_b(\iota\ten\unit,\iota\ten\unit)\cong \bds_s\Hom_{\Ccal}(U_s,U_s) = \bds_s\Com\iota_s$. Therefore the matrix $\{p_\mu(s,t)\}_{s,t\in\Irr(\Ccal)}$ determined by the identity
\begin{equation}\label{Def_RW_on_center_eq1}
P_\mu(\kappa^{t,\unit})\kappa^{s,\unit} = p_\mu(s,t)\kappa^{s,\unit}
\end{equation}
is well-defined and has entries in $\Com$. It describes a Markov kernel on $\Irr(\Ccal)$. Indeed, if $\mu$ is a probability measure, then $\|P_\mu\|\leq 1$. Clearly $P_\mu$ is positive, thus $p_\mu(s,t)\in[0,1]$ for all $s,t$. Furthermore if we denote by $I$ the natural transformation which is the identity on all objects, then $\sum_{t}P_r(\kappa^{t,\unit})_s = P_r(I)_s=\iota_s$ for all $r,s$. Therefore $\sum_t p_\mu(s,t) \iota_s=\sum_t P_\mu(\kappa^{t,\unit})_s = \iota_s$ and hence $\sum_t p_\mu(s,t)=1$.

Define the map $\vee\colon P_s\mapsto \check P_s:=P_{\bar s}$ and extend it to operators $P_U$ and $P_\mu$ by antilinearity.
\end{Def}

\begin{Lem}\label{on_center}
For the operators $P_U$, $P_\mu$ defined above the following properties hold:
\begin{enumeraterm}
\item $P_U\circ P_V = P_{V\ten U}$;
\item $\check P_\mu= P_{\bar{\mu}}$ and $\check P_U = P_{\bar{U}}$;
\item $\check P_U \circ \check P_V = \check P_{U\ten V}$;
\item $p_\mu(s,t) = \sum_r \mu(r) m_{r,s}^t \,\frac{d_t}{d_rd_s}$;
\item $p^n_{\bar\mu}(s,t) = \big(\frac{d_t}{d_s}\big)^2 p^n_\mu(t,s)$;
\item $P^n_\mu(\kappa^{t,\unit})\kappa^{s,\unit}=p_\mu^n(s,t)\kappa^{s,\unit}$.
\end{enumeraterm}
\end{Lem}

This Lemma can be proved by the same methods as used in the proof of \cite[Lem.\ 2.4]{NeshveyevTuset04}. We leave it as an exercise for the reader.

\subsection{The categorical Martin boundary}\label{subsec_categorical_Martin_boundary}
Even though this subsection is called ``The categorical Martin boundary'', we will start by presenting the Poisson boundary, as that one is easier to define and we will need it later in this paper. After stating the main properties we move our attention to the Martin boundary for random walks on \Cstar tensor categories. We need to prove some additional results before we can give its definition.

\begin{Def}[\cite{NeshveyevYamashita14a}]
Let $\mu$ be a probability measure on $\Irr(\Ccal)$. A natural transformation $\eta\in \nat_b(\iota\ten U,\iota\ten V)$ is called {\it $\mu$-harmonic} if $P_\mu(\eta)=\eta$.
\end{Def}

Given $T\in\Hom_\Ccal(X,Y)$, consider $\Ecal(T)=(\iota_U\ten T)_U\in\Hom_{\Ccal_{-1}}(X,Y)$ (see \eqref{Not_cat_path_spaces_eq4}). The tensor product of natural transformations \eqref{Def_tensor_product_nat_transf_eq0} gives $\kappa^{t,V}\ten\Ecal(T) = \kappa^{t,V}\ten T$ which acts as
\[
(\kappa^{t,V}\ten \Ecal(T))_U:=\kappa^{t,V}_U\ten T\in\Hom(U\ten V\ten X, U\ten V\ten Y).
\]
Note that in $\Hom_\Ccal(U\ten V, U\ten W)$
\begin{align*}
P_s(\kappa^{t,\unit}\ten T)_U = (\tr_s\ten\iota_U\ten\iota_W)((\kappa^{t,\unit}\ten T)_{U_s\ten U}) = (\tr_s\ten\iota_U)(\kappa^{t,\unit})_{U_s\ten U}\ten T = P_s(\kappa^{t,\unit})_U\ten T
\end{align*}
and thus
\begin{equation}\label{P_kappa_eq1}
P_\mu(\kappa^{t,\unit}\ten T) = P_\mu(\kappa^{t,\unit})\ten T.
\end{equation}
In particular as $\kappa^{t,V} = \kappa^{t,\unit}\ten\iota_V$, it follows that
\begin{equation}\label{P_kappa_eq2}
P_\mu(\kappa^{t,V}) = P_\mu(\kappa^{t,\unit})\ten\iota_V.
\end{equation}

\begin{Def}[\cite{NeshveyevYamashita14a}]\label{Def_cat_Poisson_boundary}
The {\it categorical Poisson boundary} of $\Ccal$ with respect to a probability measure $\mu$ on $\Irr(\Ccal)$ consists of a pair $(\Pcal(\Ccal,\mu),\Ecal)$, where $\Pcal(\Ccal,\mu)$ is the \Cstar tensor category which is the subobject completion of $\Ccal$ with morphism sets given by
\[
\Hom_{\Pcal(\Ccal,\mu)}(U,V):=\{\eta\in\nat_b(\iota\ten U,\iota\ten V)\,:\, \eta \textrm{ is $\mu$-harmonic}\}.
\]
The composition of morphisms in $\Pcal(\Ccal,\mu)$ is given by $(\nu\cdot\eta)_X := \lim_n P_\mu^n(\nu\eta)_X$. For objects $U,V\in\Ob(\Ccal)\subset\Ob(\Pcal(\Ccal,\mu))$ the tensor product is the same as in $\Ccal$, while on morphisms $\nu,\eta$ it is defined by
\[
(\nu\ten\eta)_Y:= (\nu\ten\iota_X)\cdot (\iota_U\ten\eta),
\]
where $\nu\ten\iota_X$ and $\iota_U\ten\eta$ are as in \eqref{Def_tensor_product_nat_transf_eq1} and \eqref{Def_tensor_product_nat_transf_eq2}. $\Ecal$ is the tensor functor $\Ecal\colon\Ccal\ra\Pcal(\Ccal,\mu)$ defined as $\Ecal(U)=U$ on objects and $\Ecal(T)_U:= (\iota_U\ten T)_U$ on morphisms.
\end{Def}

This Poisson boundary is well-defined. Indeed, in \cite{NeshveyevYamashita14a} Neshveyev and Yamashita showed that if $\eta$ and $\nu$ are bounded $\mu$-harmonic natural transformations, then $\eta\cdot\nu$ is again a bounded $\mu$-harmonic natural transformation. Moreover this product is associative. In addition for natural bounded transformations it holds that
\begin{align}
(P_\mu(\eta)\ten\iota_Y)_W = P_\mu(\eta\ten\iota_Y)_W, &&
(\iota_U\ten P_\mu(\nu))_W = P_\mu(\iota_U\ten\nu)_W.  \label{Def_cat_Poisson_boundary_eq1}
\end{align}
So the tensor product is well-defined. Also $(\iota_U\ten T)_U$ is harmonic for every $T\in\Hom_{\Ccal}(X,Y)$ and thus the functor $\Ecal$ is well-defined. Note that in general the unit object $\unit$ need no longer be simple in $\Pcal(\Ccal,\mu)$.

\begin{Def}\label{Def_cat_generating_transient}
If for every $s\in\Irr(\Ccal)$ there exists an $n\in\Nat$ such that $s\in\supp(P_\mu^n(\kappa^{0,\unit}))$, then $\mu$ is called {\it generating}. Consider the classical random walk $(\Irr(\Ccal),p_\mu)$ if this Markov chain is transient we call $\mu$ or $P_\mu$ {\it transient}. If $P_\mu$ is transient, we define the {\it Green kernel} $G_\mu:=\sum_{n=0}^\infty P_\mu^n$.
\end{Def}

\begin{Rem}\label{conjugate_measure}
If $\Ccal=\Rep(G)$ for a compact quantum group $G$. Then $\kappa^{s,\unit}$ corresponds to $I_s$. So a probability measure $\mu$ on $\Irr(G)$ is generating respectively transient in the sense of quantum groups (Definition \ref{Def_quantum_generating_transient}) if and only if $\mu$ is generating respectively transient in the sense of categories (Definition \ref{Def_cat_generating_transient}).

From statement (v) of Lemma \ref{on_center} it is obvious that $\mu$ is generating (respectively transient) if and only if $\bar\mu$ is generating (respectively transient).
\end{Rem}

\begin{Lem}\label{Lem_Green_well-defined}
If $\eta\in\nat_{00}(\iota\ten V,\iota\ten V)\subset\Hom_{\Ccal_{-1}}(V,V)$ and $\mu$ is transient, then $\|G_\mu(\eta)\|_\infty<\infty$ and thus $G_\mu(\eta)_U$ is well-defined for all $U\in\Ob(\Ccal)$. If in addition $\mu$ is generating, then $(G_\mu(\kappa^{0,V}))_U$ is invertible for all $U,V\in\Ob(\Ccal)$.
\end{Lem}
\begin{pf}
Put $C:=\|\eta\|_\infty$. Then $|\eta_s|\leq C\sum_{t\in\supp(\eta)} \kappa^{t,V}_s$, for all $s\in\Irr(\Ccal)$. The Markov kernel $P_\mu$ is completely positive, so the Green kernel $G_\mu$ being a sum of Markov kernels is positive as well. It follows that
\[
\|G_\mu(\eta)\|_\infty \leq \|G_\mu(|\eta|)\|_\infty \leq C \sum_{s\in\supp(\eta)} \sup_{t\in\Irr(\Ccal)}\|(G_\mu(\kappa^{s,V}))_t\| = C \sum_{s\in\supp(\eta)}\big(\sup_{t\in\Irr(\Ccal)}g_\mu^n(t,s)\big).
\]
By the maximum principle (see for instance \cite[\textsection~2.1]{Revuz84}) each of the summands $\sup_t g_\mu^n(t,s)$ is finite. Since $\eta$ is finitely supported, the quantity $\|G_\mu(\eta)\|_\infty$ is finite.

The support $\supp(\kappa^{0,V})=\{0\}$, so $G_\mu(\kappa^{0,V})_U$ is well-defined for every object $U$. Use the generating property of $\mu$ to find for each $t\in\Irr(\Ccal)$ an $n_t\in\Nat$, such that $t\in\supp(P_\mu^{n_t}(\kappa^{0,\unit}))$. From \eqref{P_kappa_eq2} we get $P_\mu^{n_t}(\kappa^{0,V}) = P_\mu^{n_t}(\kappa^{0,\unit})\ten\iota_V$. The latter is nonzero by choice of $n_t$. Since $P_\mu$ is positive and $\kappa^{0,V}$ is a positive natural transformation, we get that $P_\mu^{n_t}(\kappa^{0,V})_t$ is strictly positive (meaning that there exists $c>0$ such that $P_\mu^{n_t}(\kappa^{0,V})_t>c\iota_{U_t\ten V}$). Let $U\in\Ob(\Ccal)$ and decompose $U=\bigoplus_t m_U^t U_t$. Denote $I:= \{t\in\Irr(\Ccal)\,:\, m_U^t\neq0\}$ and $N:=\max\{n_t\,:\, t\in I\}$. Then for each $t\in I$ the inequality $\sum_{n=0}^N P_\mu^n(\kappa^{0,V})_t \geq P_\mu^{n_t}(\kappa^{0,V})_t$ implies that $\sum_{n=0}^N P_\mu^n(\kappa^{0,V})_t$ is strictly positive. Now
\[
G_\mu(\kappa^{0,V})_U \geq \bigoplus_{t\in I} m_V^t \sum_{n=0}^N P_\mu(\kappa^{0,V})_t
\]
shows that $G_\mu(\kappa^{0,V})_U$ is strictly positive, thus invertible.
\end{pf}

\begin{Def}
Suppose that $\mu$ is a generating and transient probability measure on $\Irr(\Ccal)$. The {\it Martin kernel} for $P_\mu$ is the operator given by
\[
K_{\bar\mu}\colon \nat_{00}(\iota\ten V,\iota\ten W)\ra\Hom_{\Ccal_{-1}}(V,W), \qquad K_{\bar\mu}(\eta)_U:=G_{\bar\mu}(\eta)_U\big(G_{\bar\mu}(\kappa^{0,V})_U\big)^{-1}.
\]
By Remark \ref{conjugate_measure} $\bar\mu$ is generating and transient, so Lemma \ref{Lem_Green_well-defined} implies that the Martin kernel $K_{\bar\mu}$ makes sense.
\end{Def}

Note that \eqref{P_kappa_eq2} implies by linearity that $G_\mu(\kappa^{t,V}) = G_\mu(\kappa^{t,\unit})\ten\iota_V$. As $P_s$ and thus $G_{\bar\mu}$ preserve the center, for any $\nu\in\nat_b(\iota\ten V,\iota\ten W)$ it holds that
\[
G_{\bar\mu}(\kappa^{0,W})_U\,\nu_U = (G_{\bar\mu}(\kappa^{0,\unit})\ten\iota_W)_U\,\nu_U = \nu_U\, (G_{\bar\mu}(\kappa^{0,\unit})\ten\iota_V)_U = \nu_U\, G_{\bar\mu}(\kappa^{0,V})_U.
\]
Therefore we obtain
\[
G_{\bar\mu}(\eta)_U\big(G_{\bar\mu}(\kappa^{0,V})_U\big)^{-1} = \big(G_{\bar\mu}(\kappa^{0,W})_U\big)^{-1}G_{\bar\mu}(\eta)_U.
\]
This means that the appropriate inverse can be placed on either side of $G_{\bar\mu}(\eta)$ when defining the Martin kernel.

\begin{Def}\label{Def_Martin_compactification}
Given a strict \Cstar tensor category $\Ccal$ with simple unit $\unit$ and a generating and transient probability measure $\mu$ on $\Irr(\Ccal)$, let $\tilde\Mcal'$ be the smallest \Cstar subcategory of $\Ccal_{-1}$ containing the objects $\Ob(\Ccal)$ and morphism sets $\nat_0(\iota\ten U,\iota\ten V)$ and $K_{\bar\mu}(\nat_{00}(\iota\ten U,\iota\ten V))$. The composition of morphisms in $\tilde\Mcal'$ is given by the composition of natural transformations.

The {\it Martin compactification} of $\Ccal$ with respect to $\mu$ consists of the pair $(\tilde\Mcal(\Ccal,\mu), \Ecal)$, where $\tilde\Mcal(\Ccal,\mu)$ is direct sum and subobject completion of $\tilde\Mcal'$ and $\Ecal\colon\Ccal\ra\tilde\Mcal(\Ccal,\mu)$ is the restriction of the canonical functor $\Ecal\colon\Ccal\ra \Ccal_{-1}$ introduced in Notation \ref{Not_cat_path_spaces} above.
\end{Def}

The functor $\Ecal$ is well-defined. Indeed, by linearity Identity \eqref{P_kappa_eq1} implies that $G_\mu(\kappa^{0,\unit}\ten T) = G_\mu(\kappa^{0,\unit})\ten T$ for $T\in\Hom_\Ccal(V,W)$.
Therefore
\[
K_\mu(\kappa^{0,\unit}\ten T) = G_\mu(\kappa^{0,\unit}\ten T)G_\mu(\kappa^{0,V})^{-1} = (G_\mu(\kappa^{0,\unit})\ten T)(G_\mu(\kappa^{0,\unit})\ten\iota_V)^{-1} = \iota\ten T.
\]
Again note that the unit object $\unit$ of $\tilde\Mcal(\Ccal,\mu)$ need no longer be simple.

\begin{Lem}\label{Martin_compactification_tensor_product}
The category $\tilde\Mcal(\Ccal,\mu)$ forms a \Cstar tensor category.
\end{Lem}
\begin{pf}
Since $\Ccal_{-1}$ is a \Cstar tensor category we only need to show that $\tilde\Mcal(\Ccal,\mu)$ is closed under tensor products. Recall the tensor product in $\Ccal_{-1}$ \eqref{Def_tensor_product_nat_transf_eq0}. It suffices to show that
\begin{enumerate}[label=(\alph*)]
\item $(\eta\ten\iota_Y)\in\nat_0(\iota\ten U\ten Y,\iota\ten V\ten Y)$ if $\eta\in\nat_0(\iota\ten U,\iota\ten V)$;
\item $(\iota_U\ten\nu)\in\nat_0(\iota\ten U\ten X,\iota\ten U\ten Y)$ if $\nu\in\nat_0(\iota\ten X,\iota\ten Y)$;
\item $K_{\bar\mu}(\eta)\ten\iota_Y = K_{\bar\mu}(\eta\ten\iota_Y)$ if $\eta\in\nat_{00}(\iota\ten U,\iota\ten V)$;
\item $\iota_U\ten K_{\bar\mu}(\nu) \in \tilde\Mcal(\Ccal,\mu)$ if $\nu\in\nat_{00}(\iota\ten X,\iota\ten Y)$.
\end{enumerate}
(a) is trivial, because $(\eta\ten\iota_Y)_U = \eta_U\ten\iota_Y$, thus $\supp(\eta\ten\iota_Y)=\supp(\eta)$ and $\|(\eta\ten\iota_Y)_W\| = \|\eta_W\|$.\\
For (b), suppose that $\eta$ has finite support. By Frobenius reciprocity
\begin{align*}
|\supp(\iota_U\ten \nu)| &= |\{s\in \Irr(\Ccal)\,:\, \eta_{U_s\ten U}\neq 0\}|\\
&= |\{s\in \Irr(\Ccal)\,:\, \exists t\in\supp(\nu),\, m_{s,U}^t\neq 0\}| \\
&= |\{s\in \Irr(\Ccal)\,:\, \exists t\in\supp(\nu),\, m_{U, \bar t}^{\bar s}\neq 0\}|\\
&\leq \sum_{t\in\supp(\nu)} |\{s\in\Irr(\Ccal)\,:\, m_{U, \bar t}^{\bar s}\neq 0\}| <\infty.
\end{align*}
Hence $\iota_U\ten\nu \in \nat_{00}(\iota\ten(U\ten X),\iota\ten(U\ten Y))$. Since $\nat_0(\iota\ten(U\ten X),\iota\ten(U\ten Y))$ is the norm-closure of $\nat_{00}(\iota\ten(U\ten X),\iota\ten(U\ten Y))$ the claim follows. \\
(c) Use \eqref{Def_cat_Poisson_boundary_eq1} to obtain
\begin{align*}
K_{\bar\mu}(\eta)\ten\iota_Y &= (G_{\bar\mu}(\eta)G_{\bar\mu}(\kappa^{0,U})^{-1})\ten\iota_Y = (G_{\bar\mu}(\eta)\ten\iota_Y)(G_{\bar\mu}(\kappa^{0,U})\ten\iota_Y)^{-1}\\
&= (G_{\bar\mu}(\eta\ten\iota_Y))(G_{\bar\mu}(\kappa^{0,U}\ten\iota_Y))^{-1} = K_{\bar\mu}(\eta\ten\iota_Y).
\end{align*}
Part (d) is slightly more complicated. Again by \eqref{Def_cat_Poisson_boundary_eq1} we obtain
\begin{align*}
\iota_U\ten K_{\bar\mu}(\nu) &= \iota_U\ten(G_{\bar\mu}(\nu)G_{\bar\mu}(\kappa^{0,X})^{-1}) = (\iota_U\ten G_{\bar\mu}(\nu))(\iota_U\ten(G_{\bar\mu}(\kappa^{0,X}))^{-1}) \\
&= G_{\bar\mu}(\iota_U\ten\nu))(G_{\bar\mu}(\iota_U\ten\kappa^{0,X}))^{-1}.
\end{align*}
Since $\mu$ is generating by Lemma \ref{Lem_Green_well-defined} the morphism $G_{\bar\mu}(\kappa^{0,U\ten X})$ is invertible. It follows that $\iota_U\ten K_{\bar\mu}(\nu)$ can be written as
\begin{align*}
\iota_U\ten K_{\bar\mu}(\nu) &= G_{\bar\mu}(\iota_U\ten\nu))(G_{\bar\mu}(\kappa^{0,U\ten X}))^{-1} G_{\bar\mu}(\kappa^{0,U\ten X})(G_{\bar\mu}(\iota_U\ten\kappa^{0,X}))^{-1}\\
&= K_{\bar\mu}(\iota_U\ten\nu) K_{\bar\mu}(\iota_U\ten\kappa^{0,X})^{-1},
\end{align*}
which is an element of $\tilde\Mcal(\Ccal,\mu)$.
\end{pf}

Let $\tilde \Mcal'$ be as in Definition \ref{Def_Martin_compactification}. Define $\Mcal'$ to be the category with $\Ob(\Mcal'):=\Ob(\tilde\Mcal')=\Ob(\Ccal)$ and
\[
\Hom_{\Mcal'}(V,W):=\{[\eta]\,:\, \eta\in\Hom_{\tilde\Mcal'}(V,W)\},
\]
where $[\eta]$ denotes the equivalence class of $\eta$ in $\nat_b(\iota\ten V,\iota\ten W)/\nat_0(\iota\ten V,\iota\ten W)$. Let $\Mcal(\Ccal,\mu)$ be the subobject and direct sum completion of $\Mcal'$. Then we can make the following observations:
\begin{enumeraterm}
\item The equivalence relation is well-defined and turns $\Mcal'$ into a category;
\item $\lbrack\,\cdot\,\rbrack\colon\tilde\Mcal'\ra\Mcal'$ is a functor and defines a tensor and $*$-structure on $\Mcal'$ by $[\nu]\ten [\eta]:=[\nu\ten\eta]$ and $[\eta]^*:=[\eta^*]$;
\item $\lbrack\,\cdot\,\rbrack$ extends to a full unitary tensor functor $\tilde\Mcal(\Ccal,\mu)\ra\Mcal(\Ccal,\mu)$.
\end{enumeraterm}

\begin{Def}\label{Def_Martin_boundary}
The {\it Martin boundary} of $\Ccal$ with respect to $\mu$ is the pair $(\Mcal(\Ccal,\mu),\Ecal)$ where $\Mcal(\Ccal,\mu)$ is the \Cstar tensor category introduced above and $\Ecal\colon \Ccal\ra \Mcal(\Ccal,\mu)$ is the unitary tensor functor which is the composition of the functor $\Ecal\colon\Ccal\ra\tilde\Mcal(\Ccal,\mu)$ with $\lbrack\,\cdot\,\rbrack$.
\end{Def}

It is often easier to work with $\nat_b(\iota\ten V,\iota\ten V)$ than with $\nat_b(\iota\ten V,\iota\ten W)$ for $V\neq W$. Fortunately most of the times one can reduce to the former case by taking the direct sum of $V$ and $W$.

\begin{Not}
Let $p_V\in\Hom_\Ccal(V\oplus W,V)$ and $p_W\in\Hom_\Ccal(V\oplus W,W)$ be the projections onto $V$ and respectively $W$. To be precise $p_Vp_V^*=\iota_V$ and $p_V^*p_V=q_V$ where $q_V$ is a projection in $\End_\Ccal(V\oplus W)$, similarly for $p_W$. Then $q_V+q_W=\iota_{V\oplus W}$. Suppose that $\eta\in\nat_b(\iota\ten V,\iota\ten W)$, the {\it extension} of $\eta$ is defined as $\eta^e:=(\eta_U^e)_U$, where $\eta_U^e$ equals the composition
\[
\xymatrix{U\ten (V\oplus W) \ar[r]^-{\iota\ten p_V} & U\ten V \ar[r]^{\eta_U} & U\ten W \ar[r]^-{\iota\ten p_W^*} & U\ten (V\oplus W)}.
\]
Similarly for $\rho\in\nat_b(\iota\ten(V\oplus W), \iota\ten(V\oplus W))$, the {\it restriction} of $\rho$ is defined as $\rho^r:=(\rho^r_U)_U$, where $\rho^r_U$ equals the composition
\[
\xymatrix{U\ten V \ar[r]^-{\iota\ten p_V^*} & U\ten (V\oplus W) \ar[r]^{\rho_U} & U\ten(V\oplus W) \ar[r]^-{\iota\ten p_W} & U\ten W}.
\]
\end{Not}

The extension and restriction again yield natural transformations. Indeed, suppose that $T\in\Hom_\Ccal(U,X)$ and $\eta$ is a natural transformation. Since $T\ten\iota$ acts in different legs than the projections $\iota\ten p_V$ and $\iota\ten p_W$, it is immediate that
\begin{align*}
(T\ten\iota_{V\oplus W})(\iota\ten p_W^*)\eta_U(\iota\ten p_V) &= (\iota\ten p_W^*)(T\ten\iota_W)\eta_U(\iota\ten p_V) = (\iota\ten p_W^*)\eta_U(T\ten\iota_V)(\iota\ten p_V)\\
&= (\iota\ten p_W^*)\eta_U(\iota\ten p_V)(T\ten\iota_{V\oplus W}).
\end{align*}
Thus $\eta^e$ is a natural transformation. Since $\|\eta^e_U\| = \|(\iota\ten p_W^*)\eta_U (\iota\ten p_V)\|\leq \|\eta_U\|$, the natural transformation $\eta^e$ is uniformly bounded whenever $\eta$ is and it is vanishing at infinity if $\eta$ is vanishing at infinity. A similar argument works for the restriction. Moreover observe that since $p_Vp_V^*=\iota_V$ and $p_Wp_W^*=\iota_W$ it holds that
\[
(\eta^e)^r_U = (p_W\ten\iota)(p_W^*\ten\iota)\eta_U(p_V\ten\iota)(p_V^*\ten\iota)=\eta_U,
\]
so $(\eta^e)^r=\eta$.

\begin{Lem}\label{equivalent_Martin_kernels}
Suppose that $\mu$ is a generating and transient probability measure on $\Irr(\Ccal)$, then \begin{align*}
K_\mu(\eta^e) &= K_\mu(\eta)^e &&\textrm{for }\eta\in\nat_{00}(\iota\ten V,\iota\ten W); \\
K_\mu(\rho^r) &= K_\mu(\rho)^r &&\textrm{for }\rho\in\nat_{00}(\iota\ten (V\oplus W),\iota\ten (V\oplus W)).
\end{align*}
\end{Lem}
\begin{pf}
First suppose that $\eta\in\nat_b(\iota\ten V,\iota\ten W)$ and $t\in\Irr(\Ccal)$, then
\begin{align*}
P_t(\eta^e)_U &= (\tr_{U_t}\ten\iota_{U\ten (V\oplus W)})\big((\iota_{U_t} \ten \iota_U\ten p_W^*)\eta_{U_t\ten U}(\iota_{U_t}\ten \iota_U\ten p_V)\big)\\
&= (\iota_U\ten p_W^*) \big((\tr_{U_t}\ten\iota_{U\ten W})(\eta_{U_t\ten U})\big)(\iota_U\ten p_V) = (P_t(\eta))_U^e.
\end{align*}
Thus by linearity in $\mu$ we get $P_\mu(\eta^e)= P_\mu(\eta)^e$. Taking sums gives $G_\mu(\eta^e)=G_\mu(\eta^e)$ for finitely supported $\eta$.
We know from \eqref{P_kappa_eq2} that $G_\mu(\kappa^{t,V}) = G_\mu(\kappa^{t,\unit})\ten\iota_V$ and thus
\begin{align*}
K_\mu(\eta)^e &= (\iota\ten p_W^*) G_\mu(\eta)(G_\mu(\kappa^{0,\unit}\ten\iota_V))^{-1} (\iota\ten p_V)\\
&= (\iota\ten p_W^*) G_\mu(\eta)(\iota\ten p_V)(G_\mu(\kappa^{0,\unit})^{-1}\ten\iota_{V\oplus W}) \\
&= G_\mu(\eta)^e(G_\mu(\kappa^{0,\unit})^{-1}\ten\iota_{V\oplus W})= K_\mu(\eta^e).
\end{align*}
The second identity can be proved similarly.
\end{pf}

\begin{Rem}
Observe that $\End_{\tilde\Mcal(\Ccal,\mu)}(U)$ equals the \Cstar algebra generated by the morphisms $K_{\bar\mu}(\nat_{00}(\iota\ten U,\iota\ten U))$ and $\nat_0(\iota\ten U,\iota\ten U)$. Indeed, it suffices to show that if $V\in\Ob(\Ccal)$ is another object and $\eta$ is an element of the \Cstar algebra generated by $\nat_0(\iota\ten (U\oplus V),\iota\ten (U\oplus V))$ and $K_{\bar\mu}(\nat_{00}(\iota\ten (U\oplus V),\iota\ten (U\oplus V)))$, then $(\iota\ten p_U)\eta(\iota\ten p_U^*)$ is an element of the \Cstar algebra generated by $\nat_0(\iota\ten U,\iota\ten U)$ and $K_{\bar\mu}(\nat_{00}(\iota\ten U,\iota\ten U))$. If $\eta$ is a generator, so an element in $\nat_0(\iota\ten (U\oplus V),\iota\ten (U\oplus V)) \cup K_{\bar\mu}(\nat_{00}(\iota\ten (U\oplus V),\iota\ten (U\oplus V)))$, this can directly be verified by a proof similar to Lemma \ref{equivalent_Martin_kernels}. But then it immediately holds for all $\eta$.

Since $\tilde\Mcal(\Ccal,\mu)$ is a \Cstar category we also have that
\[
\Hom_{\tilde\Mcal'}(V,W) = (\End_{\tilde\Mcal'}(V\oplus W))^r,
\]
where $r$ denotes the restriction.
\end{Rem}

\subsection{Categorical convergence to the boundary}\label{subsec_categorical_boundary_convergence}
It is also possible to put convergence to the boundary in a categorical framework. In this section we work towards the definition.

\begin{Not}
Let $\mu$ be a probability measure on $\Irr(\Ccal)$. For $n>m$ define recursively the maps $\Hom_{\Ccal_{-n}}(U,V)\ra\Hom_{\Ccal_{-m}}(U,V)$ by
\[
(\underbrace{\tr_\mu\ten\cdots\ten\tr_\mu}_{n-m}\ten\iota^{\ten m})(\eta):= (\tr_\mu\ten\iota^{\ten m}) ((\underbrace{\tr_\mu\ten\cdots\ten\tr_\mu}_{n-m-1}\ten\iota^{\ten m+1})(\eta)), \qquad (\eta\in\Hom_{\Ccal_{-n}}(U,V)),
\]
where the base case $m=n-1$ is defined by \eqref{Def_cat_maps_eq3}. Denote
\[
(\tr_\mu\ten\cdots\ten\tr_\mu\ten\tr_U)(\eta) := \tr_U((\tr_\mu\ten\cdots\ten\tr_\mu)(\eta)), \qquad(\eta\in\Hom_{\Ccal_{-n}}(U,U)).
\]
Write
\begin{equation}\label{Not_norm_path_space_eq1}
\|\eta\|_{\mu^{\ten n}}:= \big((\underbrace{\tr_\mu\ten\cdots\ten\tr_\mu}_n\ten\tr_U)(\eta^*\eta)\big)^\half.
\end{equation}
Since $\tr_\mu$ is a positive linear functional, $\|\cdot\|_{\mu^{\ten n}}$ defines a seminorm on $\Hom_{\Ccal_{-n}}(U,V)$. As $\eta^*\in\Hom_{\Ccal_{-n}}(V,U)$ Equation \eqref{Not_norm_path_space_eq1} also defines
\[
\|\eta^*\|_{\mu^{\ten n}}= \big((\underbrace{\tr_\mu\ten\cdots\ten\tr_\mu}_n\ten\tr_V)(\eta\eta^*)\big)^\half.
\]
\end{Not}

\begin{Def}\label{Def_cat_regular}
A natural transformation $\eta\in\Hom_{\Ccal_{-1}}(U,V)$ is called {\it $\mu$-regular} if the following condition holds: for every $\eps>0$ there exists $N\in\Nat$ such that for all $n>m\geq N$
\begin{align*}
\big\|\hat\Delta^{n-1}(\eta)-\iota^{\ten n-m}\ten\hat\Delta^{m-1}(\eta)\big\|_{\mu^{\ten n}}^2<\eps, \qquad
\big\|(\hat\Delta^{n-1}(\eta)-\iota^{\ten n-m}\ten\hat\Delta^{m-1}(\eta))^*\big\|_{\mu^{\ten n}}^2<\eps.
\end{align*}
Denote the set of $\mu$-regular elements in $\Hom_{\Ccal_{-1}}(U,V)$ by $R_{\Ccal,\mu}(U,V)$ or simply $R_\mu(U,V)$. These morphisms can be collected into a category. Denote $R(\Ccal,\mu)$ for the subcategory of $\Ccal_{-1}$ with $\Hom_{R(\Ccal,\mu)}(U,V):=R_{\Ccal,\mu}(U,V)$, this is indeed a subcategory, see Proposition \ref{Prop_regular_category}.
\end{Def}

For $U=V$ regularity can be formulated in an equivalent way. First, we define the path spaces of paths of infinite length (cf.\ \cite[\textsection~3.1]{NeshveyevYamashita14a}). For this consider the von Neumann algebras $\Hom_{\Ccal_{-n}}(U,U)$. Fix $\mu$ and define a conditional expectation
\begin{equation}\label{Def_path_space_eq1}
E_{n+1,n}\colon\Hom_{\Ccal_{-(n+1)}}(U,U)\ra\Hom_{\Ccal_{-n}}(U,U), \qquad E_{n+1,n}(\eta):= (\tr_\mu\ten\iota^{\ten n})(\eta).
\end{equation}
Clearly $E_{n+1,n}(\iota\ten\eta)=\eta$, so the embedding $\Hom_{\Ccal_{-n}}(U,U) \hookrightarrow \Hom_{\Ccal_{-(n+1)}}(U,U)$ is preserved by $E_{n+1,n}$. For $n>m$ define $E_{n,m}:= E_{m+1,m}\circ\cdots\circ E_{n,n-1}$. Define a state $\varphi^{(n)}_U:= \tr_U\circ E_{n,0}$ on $\Hom_{\Ccal_{-n}}(U,U)$. This collection $(\varphi_U^{(n)})_n$ gives a state $\varphi^{(\infty)}_U$ on the union $\bigcup_n \Hom_{\Ccal_{-n}}(U,U)$. Define the von Neumann algebra $\Hom_{\Ccal_{-\infty}}(U,U)$ as the completion of $\bigcup_n \Hom_{\Ccal_{-n}}(U,U)$ in the GNS representation defined by the state $\varphi^{(\infty)}_U$.
Denote the composition
\begin{align*}
j_n\colon\Hom&_{\Ccal_{-1}}(U,U)\ra \Hom_{\Ccal_{-n}}(U,U) \hookrightarrow \Hom_{\Ccal_{-\infty}}(U,U),\\
&\eta\mapsto \hat\Delta^{n-1}(\eta)\mapsto \cdots\ten\iota\ten\iota\ten \hat\Delta^{n-1}(\eta).
\end{align*}
It follows that a morphism $\eta\in \Hom_{\Ccal_{-1}}(U,U)$ is $\mu$-regular if and only if $\sstar\lim_n j_n(\eta)$ exists in the von Neumann algebra $\Hom_{\Ccal_{-\infty}}(U,U)$. Again often it is easier to work with $U=V$.

\begin{Lem}\label{Lem_extension_restriction}
Given $\eta\in\Hom_{\Ccal_{-1}}(U,V)$ and $\nu\in\Hom_{\Ccal{-1}}(U\oplus V, U\oplus V)$. The following holds:
\begin{enumeraterm}
\item if $\eta\in R_\mu(U,V)$, then $\eta^e\in R_\mu(U\oplus V, U\oplus V)$;
\item if $\nu\in R_\mu(U\oplus V, U\oplus V)$, then $\nu^r\in R_\mu(U,V)$.
\end{enumeraterm}
\end{Lem}
\begin{pf}
Straightforward computations show that for $\eta\in\Hom_{\Ccal_{-n}}(U,V)$ and $\nu\in\Hom_{\Ccal{-n}}(U\oplus V, U\oplus V)$ the restriction and extension operations satisfy:
\begin{align*}
(\iota^{\ten m}\ten\hat\Delta\ten\iota^{\ten n-m-1})(\eta^e) &= ((\iota^{\ten m}\ten\hat\Delta\ten\iota^{\ten n-m-1})(\eta))^e;\\
(\iota^{\ten m}\ten\hat\Delta\ten\iota^{\ten n-m-1})(\nu^r) &= ((\iota^{\ten m}\ten\hat\Delta\ten\iota^{\ten n-m-1})(\nu))^r;\\
\iota\ten\eta^e &= (\iota\ten\eta)^e;\\
\iota\ten\nu^r &= (\iota\ten\nu)^r.
\end{align*}
Similarly we have
\begin{align*}
\|\eta^e\|_{\mu^{\ten n}} = \frac{d_U}{d_U+d_V}\|\eta\|_{\mu^{\ten n}}; &&
\|\nu^r\|_{\mu^{\ten n}} \leq \frac{d_U+d_V}{d_U}\|\nu\|_{\mu^{\ten n}}.
\end{align*}
To prove these two we will restrict ourselves to $n=1$, the case $n>1$ is only more complicated in notation. Note that
\[
\tr_{U\oplus V} = \frac{1}{d_{U\oplus V}}\Tr_{U\oplus V} = \frac{1}{d_U+d_V}\Tr_U\oplus\Tr_V = \frac{d_U}{d_U+d_V}\tr_U\oplus\,\frac{d_V}{d_U+d_V}\tr_V.
\]
As $\iota_{U\oplus V}\geq p_V^*p_V$, it follows that $\nu^*\nu = \nu^*(\iota\ten\iota_{U\oplus V})\nu\geq \nu^*(\iota\ten p_V^*)(\iota\ten p_V)\nu$ and thus
\begin{align*}
\|\eta^e\|_\mu^2&= (\tr_\mu\ten\tr_{U\oplus V})((\eta^e)^*\eta^e)\\
&= \big(\tr_\mu\ten\big(\frac{d_U}{d_U+d_V}\tr_U\oplus\,\frac{d_V}{d_U+d_V}\tr_V\big)\big) \big((\iota\ten p_U^*)\eta^*(\iota\ten p_V)(\iota\ten p_V^*)\eta(\iota\ten p_U)\big)\displaybreak[2]\\
&=\big(\tr_\mu\ten\big(\frac{d_U}{d_U+d_V}\tr_U\oplus\,\frac{d_V}{d_U+d_V}\tr_V\big)\big) \big((\iota\ten p_U^*)\eta^*\eta(\iota\ten p_U)\big) \displaybreak[2]\\
&= \frac{d_U}{d_U+d_V}(\tr_\mu\ten\tr_U)(\eta^*\eta) = \frac{d_U}{d_U+d_V}\|\eta\|_\mu^2; \\
\|\nu^r\|_\mu^2 &= (\tr_\mu\ten\tr_U)((\nu^r)^*\nu^r) = (\tr_\mu\ten\tr_U)\big((\iota\ten p_U)\nu^*(\iota\ten p_V^*)(\iota\ten p_V)\nu (\iota\ten p_U^*)\big)\\
&\leq (\tr_\mu\ten\tr_U)\big((\iota\ten p_U)\nu^*\nu (\iota\ten p_U^*)\big) \displaybreak[2]\\
&\leq \frac{d_U+d_V}{d_U}(\tr_\mu\ten(\tr_{U\oplus V}))\big((\iota\ten p_U)\nu^*\nu (\iota\ten p_U^*)\oplus (\iota\ten p_V)\nu^*\nu (\iota\ten p_V^*)\big)\\
&= \frac{d_U+d_V}{d_U}(\tr_\mu\ten\tr_{U\oplus V})(\nu^*\nu) = \frac{d_U+d_V}{d_U}\|\nu\|_\mu^2.
\end{align*}
The lemma can now be proved by combining the previous identities. Indeed,
\begin{align*}
\big\|\hat\Delta^{n-1}(\eta^e)-\iota^{\ten n-m}\ten\hat\Delta^{m-1}(\eta^e)\big\|_{\mu^{\ten n}}^2 &= \big\|\big(\hat\Delta^{n-1}(\eta)-\iota^{\ten n-m}\ten\hat\Delta^{m-1}(\eta)\big)^e\big\|_{\mu^{\ten n}}^2\\
&= \frac{d_U}{d_U+d_V}\big\|\hat\Delta^{n-1}(\eta)-\iota^{\ten n-m}\ten\hat\Delta^{m-1}(\eta)\big\|_{\mu^{\ten n}}^2; \displaybreak[2]\\
\big\|\hat\Delta^{n-1}(\nu^r)-\iota^{\ten n-m}\ten\hat\Delta^{m-1}(\nu^r)\big\|_{\mu^{\ten n}}^2 &= \big\|\big(\hat\Delta^{n-1}(\nu)-\iota^{\ten n-m}\ten\hat\Delta^{m-1}(\nu)\big)^r\big\|_{\mu^{\ten n}}^2\\
&\leq \frac{d_U+d_V}{d_U}\big\|\hat\Delta^{n-1}(\nu)-\iota^{\ten n-m}\ten\hat\Delta^{m-1}(\nu)\big\|_{\mu^{\ten n}}^2
\end{align*}
and similarly for the adjoint.
\end{pf}

Since $(\eta^e)^r=\eta$ this lemma implies that a morphism $\eta\in\Hom_{\Ccal_{-1}}(U,V)$ is $\mu$-regular if and only if $\sstar\lim_n j_n(\eta^e)$ exists in the von Neumann algebra $\Hom_{\Ccal_{-\infty}}(U\oplus V, U\oplus V)$.

\begin{Prop}\label{Prop_regular_category}
The following holds:
\begin{enumeraterm}
\item $R(\Ccal,\mu)$ forms a \Cstar tensor subcategory of $\Ccal_{-1}$;
\item if $\eta$ is $\mu$-harmonic, then $\eta$ is $\mu$-regular;
\item if $\eta\in \Hom_{R(\Ccal,\mu)}(U,V)$, then $\lim_n (\tr_\mu\ten\cdots\ten\tr_\mu)(\hat\Delta^{n-1}(\eta))\in \Hom_\Ccal(U,V)$ exists in norm.
\end{enumeraterm}
\end{Prop}
We write
\begin{equation}\label{Prop_regular_category_eq2}
\tr_\mu^\infty(\eta):=\lim_n(\tr_\mu\ten\cdots\ten\tr_\mu)(\hat\Delta^{n-1}(\eta))\in\Hom_\Ccal(U,V).
\end{equation}
Note that we cannot say that the Poisson boundary $\Pcal(\Ccal,\mu)$ is a subcategory of $R(\Ccal,\mu)$ since the product in $\Pcal(\Ccal,\mu)$ is different.

\medskip
\begin{pf}[ of Proposition \ref{Prop_regular_category}]
(i) From the definition of $\mu$-regularity it is immediate that $\eta^*\in R_\mu(V,U)$ whenever $\eta\in R_\mu(U,V)$. Moreover it is clear that $R_\mu(U,V)$ is a linear space. $R_\mu(U,V)$ has norm $\|\cdot\|_\infty$ which is inherited from $\Ccal$ and thus has the properties of a norm in a \Cstar category. Observe that $\tr_U(x^*x)\leq \tr_U(\|x\|^21)=\|x\|^2$ for any $x\in\Hom_\Ccal(U,V)$, this gives that $\|\eta\|_{\mu^{\ten n}}\leq \|\eta\|_\infty$ whenever $\eta\in\Hom_{\Ccal_{-n}}(U,V)$. Now a $3\eps$-argument can be used to show that the homomorphism sets are closed in norm, we omit the details. \\
For multiplicativity we use the same estimate as in the proof of Lemma \ref{Lem_extension_restriction}, namely that $x^*x\leq \|x\|^21$, so that $y^*x^*xy\leq \|x\|^2y^*y$ and thus
\[
\|xy\|_U^2=\tr_U((xy)^*(xy))\leq\|x\|^2\tr_U(y^*y)=\|x\|^2\|y\|_U^2.
\]
Let $\nu\in R_\mu(U,V)$, $\eta\in R_\mu(V,W)$ and $\eps>0$. Then for $n$ and $m$ large enough
\begin{align*}
\|\hat\Delta^{n-1}&(\eta\nu)-\iota^{\ten n-m}\ten\hat\Delta^{m-1}(\eta\nu)\|_{\mu^{\ten n}}\\
&\leq \|\hat\Delta^{n-1}(\eta)(\hat\Delta^{n-1}(\nu)-\iota^{\ten n-m}\ten\hat\Delta^{m-1}(\nu))\|_{\mu^{\ten n}}\\
&\qquad + \|(\hat\Delta^{n-1}(\eta)-\iota^{\ten n-m}\ten\hat\Delta^{m-1}(\eta))(\iota^{\ten n-m}\ten\hat\Delta^{m-1}(\nu))\|_{\mu^{\ten n}} \displaybreak[2]\\
&\leq \|\hat\Delta^{n-1}(\eta)\|_\infty \|\hat\Delta^{n-1}(\nu)-\iota^{\ten n-m}\ten\hat\Delta^{m-1}(\nu)\|_{\mu^{\ten n}}\\
&\qquad + \|\iota^{\ten n-m}\ten\hat\Delta^{m-1}(\nu)\|_\infty \|\hat\Delta^{n-1}(\eta)-\iota^{\ten n-m}\ten\hat\Delta^{m-1}(\eta)\|_{\mu^{\ten n}}\\
&\leq \|\eta\|_\infty \eps + \|\nu\|_\infty \eps.
\end{align*}
A similar argument works for the adjoint. Thus $\eta\nu\in R_\mu(U,W)$ and $R(\Ccal,\mu)$ is a \Cstar category.

To show that $R_\mu(\Ccal,\mu)$ admits a tensor structure we must verify that
\[
(\eta\ten\nu)\in \Hom_{R(\Ccal,\mu)}(U\ten X,V\ten Y), \qquad \textrm{if }\eta\in \Hom_{R(\Ccal,\mu)}(U,V),\;\nu\in \Hom_{R(\Ccal,\mu)}(X,Y).
\]
Here $(\eta\ten\nu)$ is defined by Identity \eqref{Def_tensor_product_nat_transf_eq0}. We have already shown that $R_\mu$ is closed under composition, so it suffices to show that the natural transformations $\eta\ten\iota_Y$ and $\iota_U\ten\nu$ are $\mu$-regular. We compute
\begin{align*}
&(\tr_\mu\ten\cdots\ten\tr_\mu\ten\tr_{U\ten Y})\\
&\qquad \big((\hat\Delta^{n-1}(\eta\ten\iota_Y)-\iota^{\ten n-m}\ten\hat\Delta^{m-1}(\eta\ten\iota_Y))(\hat\Delta^{n-1}(\eta\ten\iota_Y)-\iota^{\ten n-m}\ten\hat\Delta^{m-1}(\eta\ten\iota_Y))^*\big)\\
&=(\tr_\mu\ten\cdots\ten\tr_\mu\ten\tr_U\ten \tr_Y)\\
&\qquad \big(((\hat\Delta^{n-1}(\eta)-\iota^{\ten n-m}\ten\hat\Delta^{m-1}(\eta))\ten\iota_Y)((\hat\Delta^{n-1}(\eta)-\iota^{\ten n-m}\ten\hat\Delta^{m-1}(\eta))^*\ten\iota_Y^*)\big)\\
&=\tr_Y(\iota_Y\iota_Y^*)(\tr_\mu\ten\cdots\ten\tr_\mu\ten \tr_U)\\
&\qquad\big((\hat\Delta^{n-1}(\eta)-(\iota^{\ten n-m}\ten\hat\Delta^{m-1}(\eta)))(\hat\Delta^{n-1}(\eta)-(\iota^{\ten n-m}\ten\hat\Delta^{m-1}(\eta)))^*\big),
\end{align*}
which goes to zero as $n,\,m\ra\infty$ because $\eta$ is $\mu$-regular. Similarly it can be shown that
\[
\|\hat\Delta^{n-1}(\eta\ten\iota_Y)-\iota^{\ten n-m}\ten\hat\Delta^{m-1}(\eta\ten\iota_Y)\|_{\mu^{\ten n}}^2\ra 0, \qquad \textrm{as } n,m\ra\infty,
\]
so $\eta\ten\iota_Y\in R_\mu(U\ten Y,V\ten Y)$. For $\iota_U\ten\nu$, assume that $U$ is simple, say $U=U_t$. Suppose that $t\in\supp(\mu^{*k})$ for some $k\geq1$, then
\begin{align*}
(\tr_\mu\ten\cdots\ten\tr_\mu\ten\tr_Y)\hat\Delta^{k-1} &= \sum_{s_1,\ldots,s_k,r} m_{s_1,\ldots,s_k}^r \frac{d_r}{d_{s_1}\cdots d_{s_k}} \mu(s_1)\cdots\mu(s_k)(\tr_r\ten\tr_Y)\\
&\geq m_{s_1,\ldots,s_k}^t d_t\,\frac{\mu(s_1)\cdots\mu(s_k)}{d_{s_1}\cdots d_{s_k}}\tr_t\ten\tr_Y.
\end{align*}
Hence if $s_1,\ldots,s_k\in\supp(\mu)$ are such that $m_{s_1,\ldots,s_k}^t\geq1$, then
\[
\tr_{U_t}\ten\tr_Y\leq d_t^{-1}\,\frac{d_{s_1}\cdots d_{s_k}}{\mu(s_1)\cdots\mu(s_k)}(\tr_\mu\ten\cdots\ten\tr_\mu\ten\tr_Y)\hat\Delta^{k-1}.
\]
Observe that $(\iota_U\ten\nu)_X = \nu_{X\ten U} = \hat\Delta(\nu)_{X,U}$. So $\hat\Delta^{n-1}(\iota_U\ten\nu) = \hat\Delta^n(\nu)_{\,\cdot\,,U}$, where on the spot $\cdot$ one has to put an $n$-tuple of objects. We still assume that $U=U_t$ for some simple object $U_t$. We obtain
\begin{align*}
&(\tr_\mu\ten\cdots\ten\tr_\mu\ten\tr_{U_t\ten Y})\\
&\quad \big((\hat\Delta^{n-1}(\iota_{U_t}\ten\eta)-\iota^{\ten n-m}\ten\hat\Delta^{m-1}(\iota_{U_t}\ten\eta))(\hat\Delta^{n-1}(\iota_{U_t}\ten\eta)-\iota^{\ten n-m}\ten\hat\Delta^{m-1}(\iota_{U_t}\ten\eta))^*\big)\\
&=(\tr_\mu\ten\cdots\ten\tr_\mu\ten\tr_{U_t\ten Y}) \\
&\quad\big((\hat\Delta^n(\eta)_{\,\cdot\,,U_t}-\iota^{\ten n-m}\ten\hat\Delta^m(\eta)_{\,\cdot\,,U_t})(\hat\Delta^n(\eta)_{\,\cdot\,,U_t}-\iota^{\ten n-m}\ten\hat\Delta^m(\eta)_{\,\cdot\,,U_t})^*\big) \displaybreak[2]\\
&\leq d_t^{-1}\frac{d_{s_1}\cdots d_{s_k}}{\mu(s_1)\cdots\mu(s_k)}(\tr_\mu\ten\cdots\ten\tr_\mu\ten\tr_Y)\\
&\quad \big((\hat\Delta^{n+k-1}(\eta)-\iota^{\ten n-m}\ten\hat\Delta^{m+k-1}(\eta))(\hat\Delta^{n+k-1}(\eta)-\iota^{\ten n-m}\ten\hat\Delta^{m+k-1}(\eta))^*\big)\\
&=d_t^{-1}\frac{d_{s_1}\cdots d_{s_k}}{\mu(s_1)\cdots\mu(s_k)} \|\hat\Delta^{n+k-1}(\eta)-\iota^{\ten n-m}\ten\hat\Delta^{m+k-1}(\eta)\|^2_{\mu^{\ten n+k}}\ra0, \qquad\textrm{as } n,m\ra\infty.
\end{align*}
Similarly it can be shown that
\[
\|(\hat\Delta^{n-1}(\iota_{U_t}\ten\eta)-\iota^{\ten n-m}\ten\hat\Delta^{m-1}(\iota_{U_t}\ten\eta))^*\|^2_{\mu^{\ten n}}\ra 0, \qquad\textrm{as } m,n\ra\infty.
\]
If $U$ is not simple, then decomposing $U$ into simple objects and taking direct sums gives the result.\\
(ii) This is similar to \cite[\textsection~3.1]{NeshveyevYamashita14a}. Let $\eta\in\Hom_{\Ccal_{-1}}(U,V)$ be $\mu$-harmonic. By (the results in the proofs of) Lemmas \ref{equivalent_Martin_kernels} and \ref{Lem_extension_restriction} we may assume $U=V$, since if $U\neq V$ we can consider the extension $\eta^e$ to $U\oplus V$.
Again we apply the noncommutative martingale convergence theorem. Consider the sequence $(j_n(\eta))_{n=1}^\infty\subset\Hom_{\Ccal_{-\infty}}(U,U)$. Recall the conditional expectations $E_{n,m}$ defined in \eqref{Def_path_space_eq1}. It holds that for $n\geq m$
\begin{align}
E_{n,m} (\hat\Delta^{n-1}(\eta)) &= (\tr_\mu\ten\cdots\ten\tr_\mu\ten\iota^{\ten m})\hat\Delta^{n-1}(\eta) \notag\\
&= (\tr_\mu\ten\cdots\ten\tr_\mu\ten\iota^{\ten m})\big((\iota^{\ten n-m}\ten\hat\Delta^{m-1})\hat\Delta^{n-m}(\eta)\big) \notag\\
&=\hat\Delta^{m-1}(P_\mu^{n-m}(\eta)) = \hat\Delta^{m-1}(\eta). \label{Prop_regular_category_eq3}
\end{align}
Write $\tilde E_{n,m}$ for the composition of $E_{n,m}$ with the embedding $\Hom_{\Ccal_{-m}}(U,U)\hookrightarrow \Hom_{\Ccal_{-\infty}}(U,U)$. Consider $\Hom_{\Ccal_{-m}}(U,U)$ as a subalgebra of $\Hom_{\Ccal_{-\infty}}(U,U)$ and denote
\[
\tilde E_m\colon\Hom_{\Ccal_{-\infty}}(U,U)\ra\Hom_{\Ccal_{-m}}(U,U)\subset \Hom_{\Ccal_{-\infty}}(U,U)
\]
for the conditional expectation defined by $\{\tilde E_{n,m}\}_{n\geq m}$. It follows from \eqref{Prop_regular_category_eq3} that $\tilde E_m(j_n(\eta)) = j_m(\eta)$ for $n\geq m$. As we are dealing with von Neumann algebras the noncommutative martingale convergence theorem (cf.\ Lemma \ref{lem_martingale_convergence}) shows that there exists an $\tilde\eta\in\Hom_{\Ccal_{-\infty}}(U,U)$ such that $j_n(\eta)= \tilde E_n(\tilde\eta)$ and $(j_n(\eta))_{n=1}^\infty$ converges in strong$^*$ topology to $\tilde\eta$. So $\eta$ is $\mu$-regular. \\
(iii) Let $\eta\in R_\mu(U,V)$. Using the identities in the proof of Lemma \ref{Lem_extension_restriction} we get
\begin{align}
\lim_n(\tr_\mu\ten\cdots\ten\tr_\mu)(\hat\Delta^{n-1}(\eta)) &= \lim_n(\tr_\mu\ten\cdots\ten\tr_\mu)\big((\hat\Delta^{n-1}(\eta^e))^r\big) \notag\\
&= \lim_n(\tr_\mu\ten\cdots\ten\tr_\mu)\big((\iota^{\ten n}\ten p_V)(\hat\Delta^{n-1}(\eta^e))(\iota^{\ten n}\ten p_U^*)\big) \notag \displaybreak[2]\\
&= \lim_n p_V(\tr_\mu\ten\cdots\ten\tr_\mu)\big(\hat\Delta^{n-1}(\eta^e)\big) p_U^* \notag \displaybreak[2]\\
&= p_V \big(\lim_n (\tr_\mu\ten\cdots\ten\tr_\mu)\big(\hat\Delta^{n-1}(\eta^e)\big)\big) p_U^* \notag\\
&= p_V \big(\lim_n (\cdots\ten\tr_\mu\ten\tr_\mu)(j_{n-1}(\eta^e))\big) p_U^*. \label{Prop_regular_category_eq1}
\end{align}
From the observation following the proof of Lemma \ref{Lem_extension_restriction} it follows that $\sstar\lim_n j_n(\eta^e)$ exists in $\Hom_{\Ccal_{-\infty}}(U\oplus V, U\oplus V)$. As $\Hom_\Ccal(U\oplus V, U\oplus V)$ is a finite dimensional \Cstar algebra limit \eqref{Prop_regular_category_eq1} exists in norm.
\end{pf}

\begin{Def}\label{Def_cat_convergence_boundary}
Let $\Ccal$ be a strict \Cstar tensor category with simple unit $\unit$ and $\mu$ be a probability measure on the set of irreducible objects $\Irr(\Ccal)$. The random walk on $\Ccal$ with Markov kernel $P_\mu$ {\it converges to the boundary} if the following two conditions hold:
\begin{enumeraterm}
\item $K_{\bar\mu}(\eta)\in R_\mu(U,V)$ for every $\eta\in \nat_{00}(\iota\ten U,\iota\ten V)\subset\Hom_{\Ccal_{-1}}(U,V)$;
\item for every $\eta\in\nat_{00}(\iota\ten U,\iota\ten V)$ and $\nu\in\nat_b(\iota\ten X,\iota\ten Y)$ with $P_\mu(\nu)=\nu$ it holds
\begin{equation}\label{Def_cat_convergence_boundary_eq1}
\sum_{s\in\Irr(\Ccal)} d_s^2 \tr_s((\eta\ten\nu)_{U_s}) = \tr_\mu^\infty\big(K_{\bar\mu}(\eta)\ten\nu\big).
\end{equation}
\end{enumeraterm}
Here the tensor product of natural transformations is defined by \eqref{Def_tensor_product_nat_transf_eq0} and the functionals $\tr_s$ and $\tr_\mu^\infty$ are defined by \eqref{Def_cat_maps_eq4} and respectively \eqref{Prop_regular_category_eq2}. Note that both sides of \eqref{Def_cat_convergence_boundary_eq1} are in $\Hom_\Ccal(U\ten X,V\ten Y)$.
\end{Def}

Observe that requirement (ii) of this definition makes sense due to Proposition \ref{Prop_regular_category}.

\section{Correspondence with quantum groups}
To make sure that the definition of a Martin boundary and the definition of convergence to the boundary for random walks on \Cstar tensor categories are sensible, we need to check that they correspond in some way to theory that already exists for random walks on discrete quantum groups. So for a compact quantum group $G$ one should be able to reconstruct the Martin boundary $M(\hat G,\mu)$ from $\Mcal(\Rep(G),\mu)$ and vice versa. Similarly a random walk $(\hat G,\mu)$ should converge to the boundary if and only if $(\Rep(G),\mu)$ converges to the boundary.

\subsection{Duality between \texorpdfstring{$G$-\Cstar}{G-C*}algebras and categories} \label{subsec_duality_algebras_categories}
This subsection again contains preliminary material. We review some of the results of \cite{DeCommerYamashita13}, \cite{Neshveyev13} and \cite{NeshveyevYamashita14c} that we need to prove a correspondence between the categorical picture and the quantum group picture of random walks in the next subsection.

\bigskip
Let $B$ be a $G$-\Cstar algebra, with left action $\alpha\colon B\ra C(G)\ten B$. The {\it regular subalgebra} of $B$ is denoted by
\[
\Bcal:=\{x\in B\,:\,\alpha(x)\in \Com[G]\ten_{\alg}B\}.
\]
It is a dense $*$-subalgebra of $B$ (\cite[Lem.\ 4.3]{DeCommerYamashita13}) and $\alpha$ restricts to a Hopf algebra coaction $\alpha\colon \Bcal \ra \Com[G]\ten_\alg\Bcal$.

\begin{Def}
If $\Dcal$ is a \Cstar category, define the category $\End(\Dcal)$ with objects given by \Cstar functors $\Dcal\ra\Dcal$ and morphisms $\Hom_{\End(\Dcal)}(F,G):=\nat_b(F,G)$. $\Dcal$ is called a {\it $\Rep(G)$-module category} if $\Dcal$ comes equipped with a unitary tensor functor $\Rep(G)\ra\End(\Dcal)$. If $U\in\Ob(\Rep(G))$ the induced functor in $\End(\Dcal)$ is denoted $X\mapsto X\times U$. An object $X\in\Ob(\Dcal)$ is {\it generating} if for any object $Y\in\Ob(\Dcal)$ there exists $U\in\Rep(G)$ such that $Y$ is a subobject of $X\times U$.
\end{Def}

\begin{Thm}[\protect{\cite[Thm.\ 6.4]{DeCommerYamashita13}, \cite[Thm.\ 3.3]{Neshveyev13}}]\label{thm_duality_G-Cstar}
Let $G$ be a reduced compact quantum group. The following two categories are equivalent:
\begin{enumeraterm}
\item the category of unital $G$-\Cstar algebras with unital $G$-equivariant $*$-homomorphisms as morphisms;
\item the category of pairs $(\Dcal, X)$, where $\Dcal$ is a $\Rep(G)$-module \Cstar category and $X$ is a generating object in $\Dcal$ with morphisms given by equivalence classes of unitary $\Rep(G)$-module functors respecting the designated generating objects.
\end{enumeraterm}
\end{Thm}

The formulation of this result is taken from \cite[Thm.\ 1.1]{NeshveyevYamashita14c}. For later use we explicitly describe the correspondence.
Given a $G$-\Cstar algebra $B$. Let $\Dcal'$ be the \Cstar category with objects $\Ob(\Rep(G))$ but morphism sets
\[
\Hom_{\Dcal'}(U,V):=\{T\in \Bcal\ten B(\Hcal_U,\Hcal_V)\,:\, V_{13}^*(\alpha\ten \iota)(T)U_{13} = 1\ten T\}.
\]
Let $\Dcal_B$ be the subobject completion of $\Dcal'$. Then $\Dcal_B$ is the category corresponding to $B$. The generating object is given by unit object $\unit\in\Ob(\Rep(G))\subset\Ob(\Dcal_B)$. If $U\in\Ob(\Rep(G))$, then $\iota\ten U$ defines a functor on $\Rep(G)$ as in Notation \ref{Not_cat_path_spaces}. This functor $\iota\ten U$ can be extended to the completion $\Dcal_B$ the extension is again denoted by $\iota\ten U$. The unitary tensor functor $\Rep(G)\ra\End(\Dcal_B)$ is given by $U\mapsto \iota\ten U$. \\
Conversely, let $\Dcal$ be a $\Rep(G)$-module category with generating object $X$. We may assume that $\Dcal$ is equivalent to an idempotent completion of $\Rep(G)$ with some larger morphism sets and that $X=\unit\in\Rep(G)$. Indeed, let $\Dcal'$ be the category with $\Hom_{\Dcal'}(U,V):=\Hom_\Dcal(X\times U, X\times V)$ and take the idempotent completion. Define
\begin{align*}
\Bcal&:=\bigoplus_{s\in\Irr(G)}(\bar\Hcal_s\ten\Hom_\Dcal(\unit,U_s)); &\tilde\Bcal &:= \bigoplus_{U\in\Rep(G)}(\bar\Hcal_U\ten\Hom_\Dcal(\unit,U)).
\end{align*}
For every $U\in\Rep(G)$ fix isometries $w_i\colon \Hcal_{s_i}\ra\Hcal_U$ to decompose $U$ into irreducibles. Define
\[
\pi\colon\tilde\Bcal\ra\Bcal,\qquad \pi(\bar\xi\ten T):=\sum_i\overline{w_i^*\xi}\ten w_i^*T,\qquad\textrm{for }\bar\xi\ten T\in\bar\Hcal_U\ten\Hom_\Dcal(\unit,U).
\]
Then $\pi$ is independent of the choices of the isometries $w_i$. The space $\tilde\Bcal$ becomes an associative algebra with product $\cdot$ given by
\begin{equation}\label{thm_duality_G-Cstar_eq1}
(\bar\xi\ten T)\cdot(\bar\zeta\ten S):=\overline{(\xi\ten\zeta)}\ten((T\ten\iota) S),
\end{equation}
for $\bar\xi\ten T\in\bar\Hcal_U\ten\Hom_\Dcal(\unit,U)$ and $\bar\eta\ten S\in\bar\Hcal_V\ten\Hom_\Dcal(\unit,V)$. We get a product on $\Bcal$ by $\pi(x)\pi(y):=\pi(x\cdot y)$ for $x,\,y\in\tilde\Bcal$. Define an antilinear map $\bullet$ on $\tilde\Bcal$
\begin{equation}\label{thm_duality_G-Cstar_eq2}
(\bar\xi\ten T)^\bullet:= \overline{\overline{\rho^{-1/2}\xi}}\ten(T^*\ten\iota)\bar R_U,\qquad \textrm{for }\bar\xi\ten T\in\bar\Hcal_U\ten\Hom_\Dcal(\unit,U).
\end{equation}
This map is not an involution on $\tilde\Bcal$, but defines an involution on $\Bcal$ by $\pi(x)^*:=\pi(x^\bullet)$ whenever $x\in\tilde\Bcal$. Note that $\overline{\rho^{-1/2}\xi} = (\iota\ten\bar\xi)R_U(1)$. There exists a left coaction of $\Com[G]$ on $\Bcal$ defined in the following way. Let $\{\xi_i\}_i$ be an orthonormal basis of $\Hcal_U$ and $u_{ij}$ the matrix coefficients of the representation $U$ with respect to this basis. Define $\alpha(\pi(\bar\xi_i\ten T)):=\sum_j u_{ij}\ten\pi(\bar\xi_j\ten T)$.
It can be shown that there exists a unique completion of $\Bcal$ turning it into a \Cstar algebra $B$ such that $\alpha$ extends to a left $G$-action on $B$, see \cite[\textsection~4]{DeCommerYamashita13}. This is the $G$-\Cstar algebra $B$ corresponding to the category $\Dcal$. \\
Suppose we start with a unital $G$-\Cstar algebra $B$. We form $\Dcal_B$ and let $B'$ be the algebra corresponding to $\Dcal_B$. Then there is a $*$-isomorphism \cite[\textsection~2.5]{NeshveyevYamashita14c}
\begin{equation}\label{thm_duality_G-Cstar_eq3}
\lambda\colon B'\ra B, \qquad \pi(\bar\xi\ten T)\mapsto (\iota\ten\bar\xi)T.
\end{equation}

In case there is more structure present, there is again an equivalence of categories.

\begin{Def}\label{Def_bcyd_alg}
Let $B$ be a unital \Cstar algebra. Assume that there exists a continuous actions $\alpha\colon B\ra C(G)\ten B$ and $\beta\colon B\ra M(B\ten c_0(\hat G))$. Define the left $\Com[G]$-module algebra structure
\[
\vtr\colon\Com[G]\ten B\ra B, \qquad x\vtr a:=(\iota\ten x)\beta(a), \qquad \textrm{for } x\in\Com[G] \textrm{ and } a\in B.
\]
Here $c_0(\hat G)$ is identified with a subalgebra of $\Com[G]^*$ as described by the isomorphism \eqref{def_dual_eq1}. Consider the regular subalgebra $\Bcal\subset B$. Let $S$ be the antipode of the Hopf algebra $(\Com[G],\Delta)$. The algebra $B$ is called a {\it Yetter--Drinfeld $G$-\Cstar algebra} if
\begin{equation}\label{Def_bcyd_alg_eq1}
\alpha(x\vtr a) = x_{(1)}a_{(1)}S(x_{(3)})\ten (x_{(2)}\vtr a_{(2)}), \qquad \textrm{for all } x\in\Com[G] \textrm{ and } a\in\Bcal.
\end{equation}
Here Sweedler's sumless notation is used, so $\Delta(x)=x_{(1)}\ten x_{(2)}$ and $\alpha(a)=a_{(1)}\ten a_{(2)}$. A Yetter--Drinfeld $G$-\Cstar algebra $B$ is called {\it braided-commutative} whenever
\begin{equation}\label{Def_bcyd_alg_eq2}
ab=b_{(2)}(S^{-1}(b_{(1)})\vtr a), \qquad\textrm{for all } a,b\in\Bcal.
\end{equation}
\end{Def}

\begin{Thm}[\protect{\cite[Thm.\ 2.1]{NeshveyevYamashita14c}}]\label{Thm_duality_YD}
Let $G$ be a reduced compact quantum group. The following two categories are equivalent:
\begin{enumeraterm}
\item the category of unital braided-commutative Yetter--Drinfeld $G$-\Cstar algebras with unit\-al $G$- and $\hat G$-equivariant $*$-homomorphisms as morphisms;
\item the category of pairs $(\Ccal, \Ecal)$, where $\Ccal$ is a \Cstar tensor category and $\Ecal\colon\Rep(G)\ra\Ccal$ is a unitary tensor functor such that $\Ccal$ is generated by the image of $\Ecal$. The morphisms $(\Ccal,\Ecal)\ra(\Ccal',\Ecal')$ in this category are given by the set of equivalence classes of pairs $(\Fcal,\eta)$ where $\Fcal\colon\Ccal\ra\Ccal'$ is a unitary tensor functor and $\eta\colon \Fcal\Ecal\ra\Ecal'$ is a natural unitary monoidal isomorphism.
\end{enumeraterm}
\end{Thm}

The correspondence between these two categories is the same as the correspondence given by Theorem \ref{thm_duality_G-Cstar}, but one needs to account for the extra structure present. We will only describe how the tensor product on the category can be reconstructed. Let $B$ be a unital braided-commutative Yetter--Drinfeld $G$-\Cstar algebra. Using the construction following Theorem \ref{thm_duality_G-Cstar}, we obtain a $\Rep(G)$-module \Cstar category $\Ccal_B$. If $U,V\in\Ob(\Rep(G))\subset\Ob(\Ccal_B)$, the tensor product $U\ten V$ in $\Ccal_B$ is defined as in $\Rep(G)$. For morphisms $S\in \Hom_\Ccal(U,V)$ and $T=\sum_i b_i\ten T_i\in\Hom_\Ccal(W,X)$, define
\begin{align*}
\iota_U\ten T &:= \sum_{i,k,l} (u_{kl}\vtr b_i)\ten m_{kl}\ten T_i\in \Hom_\Ccal(U\times W, U\times X);\\
S\ten \iota_X &:= S\ten\iota_{\Hcal_X} \in \Hom_\Ccal(U\times X, V\times X).
\end{align*}
Then $S\ten T:= (S\ten\iota_X)(\iota_U\ten T)$, which equals $(\iota_V\ten T)(S\ten\iota_W)$. The functor $\Ecal_B\colon\Rep(G)\ra\Ccal_B$ is given as the identity map on objects and $T\mapsto 1_B\ten T$ whenever $T\in\Hom_{\Rep(G)}(U,V)$.

\subsection{The correspondence with discrete quantum groups}\label{subsec_categorical_correspondence}
The Martin boundary and Martin compactification of a random walk on a discrete quantum group define braided-commutative Yetter--Drinfeld $G$-\Cstar algebras. Therefore by the duality described in Subsection \ref{subsec_duality_algebras_categories} these algebras define \Cstar tensor categories. These categories are shown to be unitarily monoidally equivalent to the previously defined categorical Martin boundary and compactification of $\Rep(G)$. The results in Section \ref{sec_monoidal_equivalence} above suggest that convergence to the boundary for random walks on discrete quantum groups is a property from the underlying category $\Rep(G)$ and not from the actual realization via a fiber functor $\Rep(G)\ra\textrm{Hilb}_{\textrm{f}}$. We prove that this is indeed the case.

\begin{Not}
Recall the left adjoint action $\alpha_l$ defined by \eqref{extended_adjoint_action_eq1}. Define the space
\[
\Big(\bigotimes_{-n}^{-1} l^\infty(\hat G)\Big)_{\alg}:=\Big\{x\in \bigotimes_{-n}^{-1} l^\infty(\hat G)\,:\, \alpha_l(x)\in \Com[G]\ten_{\alg} \bigotimes_{-n}^{-1} l^\infty(\hat G)\Big\}
\]
and denote $C_{-n}(\hat G)$ for the norm-closure of $(\bigotimes_{-n}^{-1} l^\infty(\hat G))_{\alg}$ in $\bigotimes_{-n}^{-1} l^\infty(\hat G)$.
The restriction $\alpha_l\colon C_{-n}(\hat G)\ra C(G)\ten C_{-n}(\hat G)$ defines an action of \Cstar algebras.
Indeed, since $\alpha_l$ is a Hopf-algebra coaction we obtain that $\alpha_l(x)\in \Com[G]\ten_{\alg} \big(\bigotimes_{-n}^{-1} l^\infty(\hat G)\big)_{\alg}$ if $x\in \big(\bigotimes_{-n}^{-1} l^\infty(\hat G)\big)_{\alg}$. Since $\alpha_l$ is continuous in norm, it extends to a map of the norm closures, again denoted by $\alpha_l\colon C_{-n}(\hat G)\ra C(G)\ten C_{-n}(\hat G)$. As $(\eps\ten\iota)\alpha_l=\iota$ on $\bigotimes_{-n}^{-1}l^\infty(\hat G)$ this identity also holds on $\big(\bigotimes_{-n}^{-1}l^\infty(\hat G)\big)_{\alg}$. It follows from \cite[Cor.\ 1.4]{NeshveyevTuset04} that $\alpha_l$ is a left $G$-action on the \Cstar algebra $C_{-n}(\hat G)$.
\end{Not}

Recall the categorical path spaces of Notation \ref{Not_cat_path_spaces}. We specialise to $\Ccal=\Rep(G)$.

\begin{Lem}\label{Lem_categorical_path_space}
For $\Ccal=\Rep(G)$ and $U,V$ finite dimensional unitary representations of $G$ it holds that
\[
\Hom_{\Rep(G)_{-n}}(U,V)\cong \Big\{T\in\Big(\bigotimes_{-n}^{-1}l^\infty(\hat G)\Big)\ten B(\Hcal_U,\Hcal_V) \,:\, V_{13}^*(\alpha_l\ten\iota_V)(T)U_{13}=1\ten T\Big\}.
\]
Consider $\Hcal_{X_i} =\bigoplus_s m_{X_i}^s \Hcal_s$. The isomorphism is given explicitly by $\eta(T)\mapsfrom T$, where
\begin{equation}\label{Lem_categorical_path_space_eq2}
\eta(T)_{X_1,\ldots, X_n}:= T|_{\Hcal_{X_1}\ten\cdots\ten\Hcal_{X_n}\ten \Hcal_U}.
\end{equation}
\end{Lem}
\begin{pf}
Note that $T\in\big(\bigotimes_{-n}^{-1}l^\infty(\hat G)\big)\ten B(\Hcal_U,\Hcal_V) $ satisfies
\begin{equation}\label{Lem_categorical_path_space_eq1}
V_{13}^*(\alpha_l\ten\iota_V)(T)U_{13}=1\ten T
\end{equation}
if and only if
\[
(1\ten T)(W\times\cdots\times W\times U) = (W\times\cdots\times W\times V)(1\ten T),
\]
so if and only if $T$ intertwines the representations $(W^{\times n}\times U)$ and $(W^{\times n}\times V)$.
Since any finite dimensional representation embeds in the regular representation coming from the multiplicative unitary $W$ the above holds if and only if
\[
(1\ten T|_{\Hcal_{X_1}\ten\cdots\ten\Hcal_{X_n}\ten \Hcal_U})(X_1\times\cdots\times X_n\times U) = (X_1\times\cdots\times X_n\times V)(1\ten T|_{\Hcal_{X_1}\ten\cdots\ten\Hcal_{X_n}\ten \Hcal_U})
\]
for all finite dimensional representations $X_1,\ldots, X_n$. Which states exactly that
\[
T|_{\Hcal_{X_1}\ten\cdots\ten\Hcal_{X_n}\ten\Hcal_U}\in\Hom_{\Rep(G)}(X_1\times\cdots\times X_n\times U, X_1\times\cdots\times X_n\times V),
\]
for all $X_1,\ldots,X_n\in\Ob(\Rep(G))$.
Now if all $X_i$ and $Y_i$ are irreducible representations and $f_i\colon X_i\ra Y_i$ is a morphism in $\Rep(G)$. Then $f_i$ is a multiple of the identity if $X_i\cong Y_i$ or zero if $X_i\not\cong Y_i$. So clearly the following diagram commutes
\[
\xymatrix{
X_1\times \cdots \times X_n\times U \ar[rr]^--{f_1\ten \cdots\ten f_n\ten\iota_U} \ar[d]_{(\eta_T)_{X_1,\ldots, X_n}} &&Y_1\times \cdots\times Y_n\times U \ar[d]^{(\eta_T)_{Y_1, \ldots, Y_n}} \\
X_1\times \cdots \times X_n\times V \ar[rr]_--{f_1\ten \cdots\ten f_n\ten\iota_V} &&Y_1\times \cdots \times Y_n\times V}.
\]
If $X_i$ and $Y_i$ are not necessarily simple then by decomposing the representations into irreducible ones it follows that the above diagram again commutes. Therefore we conclude that $T$ satisfies \eqref{Lem_categorical_path_space_eq1} if and only if $\eta(T)\in \nat_b(\iota^{\ten n}\ten U,\iota^{\ten n}\ten V)$. Thus $\eta(T)\in\Hom_{\Rep(G)_{-n}}(U,V)$. \\
Clearly $T\mapsto \eta(T)$ is an injective $*$-homomorphism. It thus remains to show that it is surjective. For this suppose $\eta\in\Hom_{\Rep(G)_{-n}}(U,V)$. Then $\eta_{X_1,\ldots, X_n}\in\Hom_{\Rep(G)}(X_1\ten\cdots\ten X_n\ten U, X_1\ten\cdots\ten X_n\ten V)$. If we write $T=\bigoplus_{s_1,\ldots,s_n}\eta_{U_{s_1},\ldots,U_{s_n}}$, then it is immediate that $\eta=\eta(T)$ and the previous computations imply that $T$ satisfies \eqref{Lem_categorical_path_space_eq1}.
\end{pf}

Let $U$ and $V$ be unitary representations. Decompose $V$ as $V=\sum_{i,j=1}^{\dim V} v_{ij}\ten m_{ij}^V\in \Com[G]\ten B(\Hcal_V)$ and similarly decompose $U$. Observe that \eqref{Lem_categorical_path_space_eq1} can be written as
\begin{equation}\label{Lem_categorical_path_space_eq3}
(\alpha_l\ten\iota)(T) = V_{13}(1\ten T)U_{13}^* = \sum_{i,j=1}^{\dim V}\sum_{k,l=1}^{\dim U} v_{ij}u_{kl}^*\ten((1\ten m_{ij}^V) T (1\ten(m_{kl}^U)^*)),
\end{equation}
which lies in $\Com[G]\ten_{\alg} \big(\bigotimes_{-n}^{-1}l^\infty(\hat G)\big)\ten B(\Hcal_U,\Hcal_V)$. It follows that if $T$ satisfies \eqref{Lem_categorical_path_space_eq1}, then $T\in (\bigotimes_{-n}^{-1} l^\infty(\hat G))_{\alg}\ten B(\Hcal_U,\Hcal_V)$. \\
It is known (see \cite[\textsection~4.1]{NeshveyevYamashita14c}) that $C_{-1}(\hat G)$ contains more structure. It is a unital braided-commutative Yetter--Drinfeld $G$-\Cstar algebra with the right action of $\hat G$ given by the comultiplication $\hat\Delta$, so the left $\Com[G]$-module structure is defined by
\begin{equation}\label{Exam_Discrete_dual_bcydalgebra_eq1}
u\vtr x:=(\iota\ten u)\hat\Delta(x), \qquad (x\in C_{-1}(\hat G),\ u\in\Com[G]).
\end{equation}
It corresponds to the \Cstar tensor category $\Rep(G)_{-1}$.
From the discussion following Theorem \ref{thm_duality_G-Cstar} and we can immediately conclude the following.

\begin{Cor}\label{Cor_equiv_path_spaces}
For a reduced compact quantum group $G$ the unital $G$-\Cstar algebra corresponding to the \Cstar category $\Rep(G)_{-n}$ via the equivalence in Theorem \ref{thm_duality_G-Cstar} is isomorphic to $C_{-n}(\hat G)$. Or equivalently the \Cstar category corresponding to the unital G-\Cstar algebra $C_{-n}(\hat G)$ is unitarily equivalent to the category $\Rep(G)_{-n}$ as a right $\Rep(G)$-module category. For $n=1$ we have unitarily monoidal equivalence of \Cstar tensor categories.
\end{Cor}

In a similar spirit as \eqref{Lem_categorical_path_space_eq3} we can write the tensor product in $\Rep(G)_{-1}$ as follows. Suppose $\eta=\eta(T)$ and $\nu=\eta(S)$ (see \eqref{Lem_categorical_path_space_eq2}) for some $T=\sum_i x_i\ten T_i \in l^\infty_{\alg}(\hat G)\ten B(\Hcal_{U_1},\Hcal_{U_2})$ and $S=\sum_j y_j\ten S_j \in l^\infty_{\alg}(\hat G)\ten B(\Hcal_{V_1},\Hcal_{V_2})$. Then \eqref{Def_tensor_product_nat_transf_eq0} translates to
\begin{equation}\label{Lem_categorical_path_space_eq4}
\eta\ten\nu = \sum_{i,j} x_iy_j^{(1)}\ten T_i\pi_{U_1}(y_j^{(2)})\ten S_j.
\end{equation}

\begin{Lem}\label{Lem_properties_isom_path_space}
For $\Ccal=\Rep(G)$, the isomorphism \eqref{Lem_categorical_path_space_eq2} satisfies
\begin{enumeraterm}
\item $(\iota\ten\eta(T)) = \eta(1\ten T)$;
\item $(\iota^{\ten m}\ten \hat\Delta\ten\iota^{\ten n-m-1})(\eta(T)) = \eta((\iota^{\ten m}\ten \hat\Delta\ten\iota^{\ten n-m-1}\ten\iota_V)(T))$;
\item $(\tr_\mu\ten\iota^{\ten n-1}(\eta(T)) = \eta((\varphi_\mu\ten\iota^{\ten n-1}\ten\iota_V)(T))$.
\end{enumeraterm}
If $n=1$, then in addition we have
\begin{enumeraterm}[resume]
\item $P_\mu(\eta(T)) = \eta((P_\mu\ten\iota)(T))$;
\item $K_\mu(\eta(T)) = \eta((K_\mu\ten\iota)(T))$.
\end{enumeraterm}
Here the first $P_\mu$ is the Markov operator on $\Rep(G)$, while the one on the right is the Markov operator on $\hat G$ and similarly for $K_\mu$.
\end{Lem}
\begin{pf}
Assertions (i) and (iii) are trivial. To prove the remaining ones, assume that $T\in\big(\bigotimes_{-n}^{-1}l^\infty(\hat G)\big)\ten B(\Hcal_U,\Hcal_V)$ satisfies \eqref{Lem_categorical_path_space_eq1}. Let $\tilde X$ be a finite dimensional representation of $G$ and $\tilde a\in\Hom(\tilde X, X_{m+1}\ten X_{m+2})$. Write $a:=\iota^{\ten m}\ten \tilde a\ten \iota^{\ten n-m-2}$ and $k:=n-m-1$. Then by definition of $\hat\Delta$ and naturality of $\eta(T)$
\begin{align*}
(\iota^{\ten m}\ten\hat\Delta\ten\iota^{\ten k})&(\eta(T))_{X_1,\ldots,X_{n+1}}\circ (a\ten\iota_U)\\
&= (\pi_{X_1}\ten\cdots\ten \pi_{X_{m+1}\ten X_{m+2}}\ten\cdots\ten \pi_{X_{n+1}}\ten\pi_U)(T)\circ (a\ten\iota_U) \displaybreak[2]\\
&= (a\ten\iota_V)\circ (\pi_{X_1}\ten\cdots\ten \pi_{\tilde X}\ten\cdots\ten\pi_{X_{n+1}}\ten\pi_U)(T) \displaybreak[2]\\
&= (\iota^{\ten m}\ten\hat\Delta\ten\iota^{\ten k}\ten\iota_V)\big((\pi_{X_1}\ten\cdots\ten\pi_{X_{n+1}}\ten\pi_U)(T)\big)\circ (a\ten\iota_U)\\
&= (\eta((\iota^{\ten m}\ten\hat\Delta\ten\iota^{\ten k}\ten\iota_V)(T)))_{X_1,\ldots,X_{n+1}}\circ (a\ten\iota_U).
\end{align*}
As this holds for all $\tilde a$, we get the second identity.
The identity with $P_\mu$ is obvious from (ii) and (iii). To prove the last one, observe that from (iv) it follows that $\eta((G_\mu\ten\iota)(T)) = G_\mu(\eta(T))$. Moreover $\kappa^{0,U}=\eta(I_0\ten \iota_U)$. Therefore
\begin{align*}
K_\mu(\eta(T)) &= G_\mu(\eta(T))(G_\mu(\eta(I_0\ten\iota_U)))^{-1} \\
&= \eta((G_\mu\ten\iota)(T))\,\eta((G_\mu\ten\iota)(I_0\ten\iota_U))^{-1} \displaybreak[2]\\
&= \eta\big(((G_\mu\ten\iota)(T))(G_\mu(I_0)^{-1}\ten\iota_U)\big) = \eta((K_\mu\ten\iota)(T)),
\end{align*}
as desired.
\end{pf}

Define the algebras
\begin{align*}
\Bcal_{-n}(\hat G)&:=\bigoplus_{s\in\Irr(G)} \bar\Hcal_s\ten\Hom_{\Rep(G)_{-n}}(\unit,U_s);\\
\tilde\Bcal_{-n}(\hat G)&:=\bigoplus_{U\in\Ob(G)} \bar\Hcal_U\ten\Hom_{\Rep(G)_{-n}}(\unit,U)
\end{align*}
and denote
\begin{equation}\label{Def_equivalent_algebras_eq1}
\lambda_n\colon \Bcal_{-n}(\hat G)\ra C_{-n}(\hat G), \qquad \pi(\bar\xi\ten \eta)\mapsto (\iota\ten\bar\xi)\eta.
\end{equation}
Let $B_{-n}(\hat G)$ be the unique completion of $\Bcal_{-n}(\hat G)$ to a $G$-\Cstar algebra such that $\lambda_n$ extends to a $G$-equivariant isomorphism $\lambda_n\colon B_{-n}(\hat G)\ra C_{-n}(\hat G)$ (cf.\ \eqref{thm_duality_G-Cstar_eq3}).

\bigskip
Let us consider the case $n=0$. Put $B=\Com$, this is a unital braided-commutative Yetter--Drinfeld $G$-\Cstar algebra with trivial left and right actions. The associated \Cstar tensor category $\Ccal_B$ equals $\Rep(G)$. Indeed, for objects $U,V\in\Rep(G)\subset\Ob(\Ccal_B)$
\begin{align*}
\Hom_{\Ccal_B}(U,V)&=\{T\in B\ten B(\Hcal_U,\Hcal_V)\,:\, V_{13}^*(\alpha\ten\iota)(T)U_{13}=1\ten T\}\\
&\cong\{S\in B(\Hcal_U,\Hcal_V)\,:\, (\iota\ten S)U= V(\iota\ten S)\} \\
&=\Hom_{\Rep(G)}(U,V).
\end{align*}
So put $\Rep(G)_0:=\Rep(G)$ and $\Bcal_0(\hat G):=\bigoplus_{s\in\Irr(G)} \bar\Hcal_s\ten\Hom_{\Ccal_B}(\unit,U_s) \cong \End_{\Rep(G)}(\unit)$. Denote the isomorphism $\lambda_0\colon \Bcal_0(\hat G)\ra \Com$, $\pi(\bar\xi\ten T)\mapsto (\iota\ten\bar\xi)T$.
The algebra $\Bcal_0(\hat G)$ is already complete, but to be consistent in notation we write $B_0(\hat G)=\Bcal_0(\hat G)$.

\begin{Lem}\label{Lem_lambda_isomorphism_subalgebras}
The map $\lambda_1$ restricts to $G$-equivariant $*$-isomorphisms
\begin{align}
\lambda_1&\colon \bigoplus_s \bar\Hcal_s\ten\nat_{00}(\iota\ten\unit,\iota\ten U_s)\ra c_{00}(\hat G); \label{Lem_lambda_isomorphism_subalgebras_eq1}\\
\lambda_1&\colon \overline{\bigoplus_s \bar\Hcal_s\ten\nat_0(\iota\ten\unit,\iota\ten U_s)}^{\|\cdot\|}\ra c_0(\hat G), \label{Lem_lambda_isomorphism_subalgebras_eq2}
\end{align}
The closure indicates the unique completion such that the space becomes a \Cstar algebra with left action $\alpha_l$, see \textsection~\ref{subsec_duality_algebras_categories}.
\end{Lem}
\begin{pf}
$B(\Hcal_s)$ is a unital $G$-\Cstar algebra with left action $\alpha_{l,s}$ given by the restriction of the left adjoint action $\alpha_l$, so $\alpha_{l,s}(x):= U_s^*(1\ten x)U_s$. Form the $\Rep(G)$-module category $\Ccal_{B(\Hcal_s)}$ as described by Theorem \ref{thm_duality_G-Cstar}. From this theorem it follows that
\[
\lambda_s\colon \bigoplus_t \bar\Hcal_t\ten\Hom_{\Ccal_{B(\Hcal_s)}}(\unit,U_t) \ra B(\Hcal_s), \qquad \pi(\bar\xi\ten\eta) \mapsto (\iota\ten\bar\xi)\eta.
\]
is a  $G$-equivariant $*$-isomorphism. Since
\begin{align*}
\Hom_{\Ccal_{B(\Hcal_s)}}(U,V) &= \{T\in B(\Hcal_s)\ten B(\Hcal_U,\Hcal_V)\,:\, U_{13}^*(\alpha_{l,s}\ten\iota)(T)V_{13} = 1\ten T\}\\
&\cong\{\eta\in\nat_b(\iota\ten U,\iota\ten V)\,:\,\supp(\eta)\subset\{s\}\},
\end{align*}
we obtain that the function $\lambda:=\ds_s \lambda_s$ acting as
\[
\lambda\colon \bigoplus_s \bigoplus_t \bar\Hcal_t\ten \{\eta\in\nat_b(\iota\ten U,\iota\ten V)\,:\,\supp(\eta)\subset\{s\}\} \ra \bigoplus_s B(\Hcal_s)
\]
is an isomorphism. But this map equals exactly $\lambda_1$ of \eqref{Lem_lambda_isomorphism_subalgebras_eq1}. Taking closures gives the second isomorphism.
\end{pf}

If $\Fcal\colon\Rep(G)_{-n}\ra\Rep(G)_{-k}$ is a functor of module categories, we obtain maps
\[
\bar\Hcal_s\ten\Hom_{\Rep(G)_{-n}}(\unit,U_s)\ra \bar\Hcal_s\ten\Hom_{\Rep(G)_{-k}}(\unit,U_s), \qquad
\bar\xi\ten\eta\mapsto \bar\xi\ten\Fcal(\eta).
\]
Taking direct sums and passing to the completion gives a $*$-morphism of \Cstar algebras (see \cite[Prop.\ 4.5]{DeCommerYamashita13}) we denote it by
\begin{equation}\label{def_functor_to_star-morphism}
(\iota_{\bar\Hcal}\ten\Fcal)\colon B_{-n}(\hat G)\ra B_{-k}(\hat G).
\end{equation}
Similarly for $\tr_\mu$ we obtain positive maps
\[
\bar\Hcal_s\ten\Hom_{\Rep(G)_{-n}}(\unit,U_s)\ra \bar\Hcal_s\ten\Hom_{\Rep(G)_{-(n-1)}}(\unit,U_s), \qquad
\bar\xi\ten\eta\mapsto \bar\xi\ten(\tr_\mu\ten\iota^{\ten n-1})(\eta).
\]
By taking direct sums and completion these maps extended to positive maps
\begin{equation}\label{def_trace_to_cond_expect}
(\iota_{\bar\Hcal}\ten\tr_\mu\ten\iota^{\ten n-1})\colon B_{-n}(G)\ra B_{-(n-1)}(G).
\end{equation}
Using the functor $(\iota\ten\cdot)$ one can embed $\bar\Hcal_s\ten\Hom_{\Rep(G)_{-(n-1)}}(\unit,U_s)\hookrightarrow \bar\Hcal_s\ten\Hom_{\Rep(G)_{-n}}(\unit,U_s)$ and thus one obtains an embedding $B_{-(n-1)}(G)\hookrightarrow B_{-n}(G)$. The map $(\iota_{\bar\Hcal}\ten\tr_\mu\ten\iota^{\ten n-1})$ defines a conditional expectation $B_{-n}(G)\ra B_{-(n-1)}(G)$.

\begin{Cor}\label{Cor_compatibility_lambda}
The following identities hold on $B_{-n}(\hat G)$ for $n\geq1$:
\begin{enumeraterm}
\item $\lambda_{n+1}\circ(\iota_{\bar\Hcal}\ten\iota\ten\cdot) = (1\ten\cdot)\circ\lambda_n$;
\item $\lambda_{n+1}\circ(\iota_{\bar\Hcal}\ten\hat\Delta\ten\iota^{\ten n-1}) = (\hat\Delta\ten\iota^{\ten n-1})\circ\lambda_n$;
\item $\lambda_{n-1}\circ(\iota_{\bar\Hcal}\ten\tr_\mu\ten\iota^{\ten n-1}) = (\varphi_\mu\ten\iota^{\ten n-1})\circ\lambda_n$,
\end{enumeraterm}
where the $*$-morphisms $(\iota_{\bar\Hcal}\ten\iota\ten\cdot)$, $(\iota_{\bar\Hcal}\ten\hat\Delta\ten\iota^{\ten n-1})$ are defined by \eqref{def_functor_to_star-morphism} and $(\iota_{\bar\Hcal}\ten\tr_\mu\ten\iota^{\ten n-1})$ by \eqref{def_trace_to_cond_expect}. In particular it holds that
\begin{enumeraterm}[resume]
\item $\lambda_1\circ(\iota_{\bar\Hcal}\ten P_\mu) = P_\mu\circ\lambda_1$;
\item $\lambda_1\circ(\iota_{\bar\Hcal}\ten K_\mu) = K_\mu\circ\lambda_1$.
\end{enumeraterm}
\end{Cor}
\begin{pf}
This is more a matter of notation than actually something new. The key part is Lemma \ref{Lem_properties_isom_path_space}.
Let $\bar\xi\ten\eta\in\bar\Hcal_s\ten\Hom_{\Rep(G)_{-n}}(\unit,U_s)\subset \Bcal_{-n}(\hat G)\subset B_{-n}(\hat G)$. By Lemma \ref{Lem_categorical_path_space} we may assume that $\eta$ is of the form $\eta(T)$ for some $T\in (\bigotimes_{-n}^{-1}l^\infty(\hat G))\ten B(\Com,\Hcal_s)$. Using Lemma \ref{Lem_properties_isom_path_space} we conclude
\begin{align*}
\lambda_{n+1}((\iota_{\bar\Hcal}\ten\hat\Delta\ten\iota^{\ten n-1})&(\bar\xi\ten\eta(T)))|_{s_1\ten\cdots\ten s_{n+1}}\\
&= (\iota^{\ten n+1}\ten\bar\xi)\big((\hat\Delta\ten\iota^{\ten n-1})(\eta(T))\big)|_{s_1\ten\cdots\ten s_{n+1}} \displaybreak[2]\\
&=(\iota^{\ten n+1}\ten\bar\xi)\big(\eta((\hat\Delta\ten\iota^{\ten n-1}\ten\iota_{\Hcal_s})(T))\big)|_{s_1\ten\cdots\ten s_{n+1}} \displaybreak[2]\\
&= (\hat\Delta\ten\iota^{\ten n-1}\ten\bar\xi)(T)|_{s_1\ten\cdots\ten s_{n+1}} \displaybreak[2]\\
&= (\hat\Delta\ten\iota^{\ten n-1})\big((\iota^{\ten n+1}\ten\bar\xi)(\eta(T))\big)|_{s_1\ten\cdots\ten s_{n+1}}\\
&=(\hat\Delta\ten\iota^{\ten n-1})\lambda_n(\bar\xi\ten\eta(T))|_{s_1\ten\cdots\ten s_{n+1}}.
\end{align*}
Identities (ii) and (iii) can be verified in an analogous way and (iv) follows immediately from the definition of $P_\mu$ (see Definition \ref{Def_categorical_random_walk}).\\
(v) Consider $\bar 1\ten\kappa^{0,\unit}$ as an element of
\[
\bar\Com\ten \Hom_{\Rep(G)}(\unit,\unit)\subset \bigoplus_s \bar\Hcal_s\ten \Hom_{\Rep(G)_{-1}}(\unit,U_s)=\Bcal_{-1}(\hat G).
\]
Then $\lambda_1(\bar 1\ten\kappa^{0,\unit})=I_0\in B(\Hcal_0)$. By linearity in $\mu$ we see from (iv) that $\lambda_1\circ(\iota_{\bar\Hcal}\ten G_\mu) = G_\mu\circ\lambda_1$. Now take $\bar\xi\ten\eta\in \bar\Hcal_s\ten\nat_{00}(\iota\ten\unit,\iota\ten U_s)$ for some $s$. We get
\begin{align*}
K_\mu\circ\lambda_1(\bar\xi\ten\eta)&= G_\mu(\lambda_1(\bar\xi\ten\eta)) G_\mu(\lambda_1(\bar 1\ten \kappa^{0,\unit}))^{-1} \\
&= \lambda_1\big((\bar\xi\ten G_\mu(\eta))\cdot(\bar 1\ten G_\mu(\kappa^{0,\unit})^{-1})\big) \displaybreak[2]\\
&= \lambda_1\big(\bar\xi\ten G_\mu(\eta)G_\mu(\kappa^{0,\unit})^{-1}\big)\\
&= \lambda_1(\bar\xi\ten K_\mu(\eta)),
\end{align*}
as desired.
\end{pf}

Suppose that $\mu$ is a generating and transient probability measure on $\Irr(G)$. The action $\alpha_l$ defines adjoint actions on the Martin compactification $\tilde M(\hat G,\mu)$ and Martin boundary $M(\hat G,\mu)$ (cf.\ \cite[Thm.\ 3.5]{NeshveyevTuset04}). Denote $M(\hat G,\mu)_{\alg}:= \tilde M(\hat G,\mu) \cap l^\infty_{\alg}(\hat G)$. This is a $*$-algebra which is norm-dense in $\tilde M(\hat G,\mu)$ (see Subsection \ref{subsec_duality_algebras_categories}). For the Martin boundary we consider
\[
M(\hat G,\mu)_{\alg}:= \tilde M(\hat G,\mu)_{\alg} / (c_0(\hat G)\cap l^\infty_{\alg}(\hat G)) = \{x\in M(\hat G,\mu)\,:\, \alpha_l(x)\in \Com[G]\ten_{\alg} M(\hat G,\mu)\}.
\]
Then again $M(\hat G,\mu)_{\alg}$ is norm-dense in $M(\hat G,\mu)$.

\begin{Lem}\label{Martin_YD-algebra}
The \Cstar algebras $\tilde M(\hat G,\mu)$ and $M(\hat G,\mu)$ are unital braided-commutative Yetter--Drinfeld $G$-\Cstar algebras. The left action of $G$ on $\tilde M(\hat G,\mu)$ is given by the restriction of $\alpha_l$ to $\tilde M(\hat G,\mu)$. The left $\Com[G]$-module structure $\vtr$ is defined by the restriction of \eqref{Exam_Discrete_dual_bcydalgebra_eq1} to $\tilde M(\hat G,\mu)$. Both actions factor through $M(\hat G,\mu)$.
\end{Lem}
\begin{pf}
By definition $\tilde M(\hat G,\mu)$ and $M(\hat G,\mu)$ are \Cstar algebras. They are unital, because $K_{\bar{\mu}}(I_0) = G_{\bar{\mu}}(I_0)G_{\bar{\mu}}(I_0)^{-1}=1$.
The mappings $\alpha_l$ and $\hat\Delta$ define a left $G$-action and respectively a right $\hat G$-action on both $\tilde M(\hat G,\mu)$ and $M(\hat G,\mu)$ (\cite[Thm.\ 3.5]{NeshveyevTuset04}). Thus $\tilde M(\hat G,\mu)$ and $M(\hat G,\mu)$ are closed under $\alpha_l$ and $\vtr$. As $C_{-1}(\hat G)$ is a unital braided-commutative Yetter--Drinfeld $G$-\Cstar algebra (\cite[\textsection~4.1]{NeshveyevYamashita14c}) and $\tilde M(\hat G,\mu)\subset C_{-1}(\hat G)$ it follows that the defining identities \eqref{Def_bcyd_alg_eq1} and \eqref{Def_bcyd_alg_eq2} also hold $\tilde M(\hat G,\mu)$ on $M(\hat G,\mu)$.
\end{pf}

\begin{Thm}\label{Thm_equivalence_Martin_boundary}
Let $G$ be a reduced compact quantum group and $\mu$ a generating and transient probability measure on $\Irr(G)$. Denote by $\tilde B(\hat G,\mu)$ and $B(\hat G,\mu)$ the unital braided-commutative Yetter--Drinfeld algebras associated to the categorical Martin compactification $\tilde \Mcal(\Rep(G),\mu)$ and respectively the categorical Martin boundary $\Mcal(\Rep(G),\mu)$ as in Theorem \ref{Thm_duality_YD}. Then $\lambda_1\colon\tilde B(\hat G,\mu)\ra \tilde M(\hat G,\mu)$ is a $*$-isomorphism preserving the left $G$-action and right $\hat G$-action. Moreover $\lambda_1$ factors through the Martin boundary $\lambda_1\colon B(\hat G,\mu)\ra M(\hat G,\mu)$.

Equivalently, the \Cstar tensor categories associated to the unital braided-commutative Yetter--Drinfeld $G$-\Cstar algebras $\tilde M(\hat G,\mu)$ and $M(\hat G,\mu)$ are unitarily monoidally equivalent to the categorical Martin compactification $\tilde\Mcal(\Rep(G),\mu)$ and respectively Martin boundary $\Mcal(\Rep(G),\mu)$ of the random walk defined by $\mu$ on the \Cstar tensor category $\Rep(G)$. These monoidal equivalences preserve the functors of $\Rep(G)$ into the respective categories.
\end{Thm}
\begin{pf}
From the construction of the categorical Martin compactification $\tilde\Mcal(\Rep(G),\mu)$ we see that $\Hom_{\tilde\Mcal(\Rep(G),\mu)}(U,V)\subset\Hom_{\Rep(G)_{-1}}(U,V)$ for any objects $U,V\in\Ob(\Rep(G))$, thus $\tilde B(\hat G,\mu)\subset B_{-1}(\hat G)$. Hence by Corollary \ref{Cor_equiv_path_spaces} the restriction $\lambda_1\colon \tilde B(\hat G,\mu)\ra C_{-1}(\hat G)$ is an injective $*$\hyph homomorphism preserving the actions of $G$ and $\hat G$. To show that $\lambda_1$ defines an isomorphism for the compactifications, it therefore suffices to show two more things:
\begin{enumerate}[label=(\alph*)]
\item $\lambda_1(\tilde B(G,\mu))\subset \tilde M(G,\mu)$;
\item $\lambda_1\colon \tilde B(G,\mu) \ra \tilde M(G,\mu)$ is surjective.
\end{enumerate}
Denote the regular subalgebra of $\tilde B(\hat G,\mu)$ by $\tilde\Bcal(\hat G,\mu)$. To prove (a), note that as a vector space
\begin{align*}
\tilde\Bcal(\hat G,\mu)&\cong\bigoplus_{s\in\Irr(G)} \bar\Hcal_s\ten\Hom_{\tilde\Mcal(\Rep(G),\mu)}(\unit,U_s)\\
&= \bigoplus_{s\in\Irr(G)} \bar\Hcal_s\ten(\Hom_{\tilde\Mcal(\Rep(G),\mu)}(\unit\oplus U_s, \unit\oplus U_s))^r,
\end{align*}
where $r$ denotes again the operation of restriction. Now $\Hom_{\tilde\Mcal(\Rep(G),\mu)}(\unit\oplus U_s, \unit\oplus U_s)$ is generated as a \Cstar algebra by $\nat_0(\iota\ten(\unit\oplus U_s),\iota\ten(\unit\oplus U_s))$ and $K_{\bar\mu}(\nat_{00}(\iota\ten(\unit\oplus U_s),\iota\ten(\unit\oplus U_s)))$. Lemma \ref{Lem_lambda_isomorphism_subalgebras} shows that if $T\in \bigoplus_s \bar\Hcal_s\ten \nat_0(\iota\ten(\unit\oplus U_s), \iota\ten(\unit\oplus U_s))$, then $\lambda_1(T)\in c_0(\hat G)$. Similarly if $S\in \bigoplus_s \bar\Hcal_s\ten \nat_{00}(\iota\ten(\unit\oplus U_s), \iota\ten(\unit\oplus U_s))$, then $\lambda_1(S)\in c_{00}(\hat G)$. In the latter case, write $S^r:=S_1\ten S_2^r\in \bigoplus_s \bar\Hcal_s\ten \nat_{00}(\iota\ten\unit, \iota\ten U_s)$ for the restriction applied to the second leg of $S$. Lemmas \ref{equivalent_Martin_kernels} and \ref{Cor_compatibility_lambda} imply that
\[
\lambda_1((\iota_{\Hcal}\ten K_{\bar\mu})(S)^r)=\lambda_1((\iota_\Hcal\ten K_{\bar\mu})(S^r)) = K_{\bar\mu}(\lambda_1(S^r)),
\]
which is thus an element of $K_{\bar\mu}(c_{00}(\hat G))\subset \tilde M(\hat G,\mu)$. Hence $\lambda_1$ maps the generators of $\tilde B(\hat G,\mu)$ in $\tilde M(\hat G,\mu)$, which proves (a). \\
To establish surjectivity, we reverse the argument. Clearly every $x\in c_0(\hat G)$ lies in the image of $\lambda_1$ (see Lemma \ref{Lem_lambda_isomorphism_subalgebras}). Assume that $x\in c_{00}(\hat G)$, by the same lemma there exists $S\in \bigoplus_s \bar\Hcal_s\ten \nat_{00}(\iota\ten\unit, \iota\ten U_s)$ such that $\lambda_1(S)=x$. Then $S^e\in \bigoplus_s \bar\Hcal_s\ten \nat_{00}(\iota\ten(\unit\oplus U_s), \iota\ten(\unit\oplus U_s))$ and $((\iota_{\bar\Hcal}\ten K_{\bar\mu})(S^e))^r\in \bigoplus_s \bar\Hcal_S\ten\Hom_{\Mcal(\Rep(G),\mu)}(\unit,U_s)$.
Invoking again Lemmas \ref{equivalent_Martin_kernels} and \ref{Cor_compatibility_lambda} gives us
\[
\lambda_1((\iota_{\bar\Hcal}\ten K_{\bar\mu})(S^e)^r) = \lambda_1((\iota_{\bar\Hcal}\ten K_{\bar\mu})(S)) = K_{\bar\mu}(\lambda_1(S)) = K_{\bar\mu}(x).
\]
Hence $\lambda_1$ is surjective.

From the second part of Lemma \ref{Lem_lambda_isomorphism_subalgebras} we immediately conclude that $\lambda_1$ factors through the Martin boundary.
\end{pf}

Since $R(\Ccal,\mu)$ is a \Cstar tensor category (see Proposition \ref{Prop_regular_category}) one could try to reconstruct the \Cstar algebra $R_{\varphi_\mu}$ of regular elements from $R(\Rep(G),\mu)$. There is however one problem, the algebra $R_{\varphi_\mu}$ is only a \Cstar subalgebra of $l^\infty(\hat G)$, so one cannot talk about actions in the von Neumann sense. On the other hand it is unknown whether the left $G$-action is continuous in the \Cstar sense.
But one can consider $R_{\varphi_\mu}\cap C_{-1}(\hat G)$ and show that this algebra admits a braided-commutative Yetter--Drinfeld structure.

\begin{Lem}\label{Lem_regular_bcYD-algebra}
Denote $\tilde R(\hat G,\mu):=\{x\in C_{-1}(\hat G)\,:\, x \textrm{ is } \mu\textrm{-regular}\}=C_{-1}(\hat G)\cap R_{\varphi_\mu}$. Then $\tilde R(\hat G,\mu)$ is a unital braided-commutative Yetter--Drinfeld $G$-\Cstar algebra. The left action of $G$ on $\tilde R(\hat G,\mu)$ is given by the restriction of $\alpha_l$ to $\tilde R(\hat G,\mu)$. The left $\Com[G]$-module structure $\vtr$ is defined by the restriction of \eqref{Exam_Discrete_dual_bcydalgebra_eq1} to $\tilde R(\hat G,\mu)$.
\end{Lem}
\begin{pf}
By definition $\tilde R(\hat G,\mu)\subset C_{-1}(\hat G)$ and $\alpha_l$ defines a continuous left action of $G$ on $C_{-1}(\hat G)$. Using the same argument as in Lemma \ref{Martin_YD-algebra} it thus suffices to show that $\tilde R(\hat G,\mu)$ is closed under the actions $\alpha_l$ and $\vtr$.

We deal with $\alpha_l$ first. Using the pentagon equation for $W$ it is easy to show that $\alpha_l(\hat\Delta^n(x)) = (\iota\ten\hat\Delta^n)(\alpha_l(x))$. Let $x\in\tilde R(\hat G,\mu)$, inferring \eqref{phi_invariant_eq1} yields for $n\geq m$
\begin{align*}
&\|(\iota\ten\hat\Delta^{n-1})(\alpha_l(x)) - (\iota\ten 1^{\ten n-m}\ten\hat\Delta^{m-1})(\alpha_l(x))\|_{h\ten\varphi_\mu^{\ten n}}^2 \\
&=(h\ten\varphi_\mu\ten\cdots\ten\varphi_\mu)\big(((\iota\ten\hat\Delta^{n-1})(\alpha_l(x)) - (\iota\ten 1^{\ten n-m}\ten\hat\Delta^{m-1})(\alpha_l(x)))^* \\
&\qquad\times((\iota\ten\hat\Delta^{n-1})(\alpha_l(x)) - (\iota\ten 1^{\ten n-m}\ten\hat\Delta^{m-1})(\alpha_l(x)))\big) \displaybreak[2]\\
&= (h\ten\varphi_\mu\ten\cdots\ten\varphi_\mu)\big(\alpha_l((\hat\Delta^{n-1}(x) - 1^{\ten n-m}\ten\hat\Delta^{m-1}(x))^*(\hat\Delta^{n-1}(x) - 1^{\ten n-m}\ten\hat\Delta^{m-1}(x)))\big) \displaybreak[2]\\
&= (\varphi_\mu\ten\cdots\ten\varphi_\mu)\big((\hat\Delta^{n-1}(x) - 1^{\ten n-m}\ten\hat\Delta^{m-1}(x))^*(\hat\Delta^{n-1}(x) - 1^{\ten n-m}\ten\hat\Delta^{m-1}(x))\big)\\
&=\|\hat\Delta^{n-1}(x) - 1^{\ten n-m}\ten\hat\Delta^{m-1}(x)\|_{\varphi_\mu^{\ten n}}^2,
\end{align*}
which tends to $0$ as $m,n\ra\infty$. The same argument holds if one replaces $x$ by $x^*$ and thus $\alpha_l(x)\in C(G)\ten \tilde R(\hat G,\mu)$.

For $\vtr$ we use an argument similar to the proof of Proposition \ref{Prop_regular_category} part (i). Assume that $x\in \tilde R(\hat G,\mu)$ and $t\in\Irr(G)$. Since $\mu$ is generating, let $k\geq 0$ be such that $t\in\supp(\mu^{*k})$. Select $s_1,\ldots,s_k\in\supp(\mu)$ such that $m_{s_1,\ldots,s_k}^t\geq 1$, then by Lemma \ref{constants_cnr}
\[
\varphi_t\leq \frac{d_{s_1}\cdots d_{s_k}}{d_t}\,m_{s_1,\ldots,s_k}^t (\varphi_{s_1}\ten\cdots\ten\varphi_{s_k})\hat\Delta^{k-1} \leq d_t^{-1}\,\frac{d_{s_1}\cdots d_{s_k}}{\mu(s_1)\cdots\mu(s_k)}\,\varphi_\mu^k.
\]
Write $C:= d_t^{-1}\,\frac{d_{s_1}\cdots d_{s_k}}{\mu(s_1)\cdots\mu(s_k)}$. It follows that
\begin{align*}
&(\underbrace{\varphi_\mu\ten\cdots\ten\varphi_\mu}_n\ten\varphi_t) \big(\big((\hat\Delta^{n-1}\ten\iota)\hat\Delta(x)-1^{\ten n-m}\ten(\hat\Delta^{m-1}\ten\iota)(\hat\Delta(x))\big)^*\\
&\qquad\times\big((\hat\Delta^{n-1}\ten\iota)\hat\Delta(x)-1^{\ten n-m}\ten(\hat\Delta^{m-1}\ten\iota)(\hat\Delta(x))\big)\big) \displaybreak[1]\\
&\leq C (\underbrace{\varphi_\mu\ten\cdots\ten\varphi_\mu}_n\ten (\underbrace{\varphi_\mu\ten\cdots\ten\varphi_\mu}_k)\hat\Delta^{k-1}) \big((\hat\Delta^n(x)-1^{\ten n-m}\ten\hat\Delta^m(x))^*\\
&\qquad\times(\hat\Delta^n(x)-1^{\ten n-m}\ten\hat\Delta^m(x))\big)\displaybreak[1]\\
&=C (\underbrace{\varphi_\mu\ten\cdots\ten\varphi_\mu}_{n+k}) \big((\hat\Delta^{n+k-1}(x)-1^{\ten n-m}\ten\hat\Delta^{m+k-1}(x))^*(\hat\Delta^{n+k-1}(x)-1^{\ten n-m}\ten\hat\Delta^{m+k-1}(x))\big)
\end{align*}
which, by regularity of $x$, tends to $0$ as $m,n\ra\infty$. So $\hat\Delta(x)\in \tilde R(\hat G,\mu)\ten l^\infty(\hat G)$ and thus $\tilde R(\hat G,\mu)$ is closed under $\vtr$.
\end{pf}

\begin{Lem}\label{Lem_estimate_trace_on_Bn}
The following estimates hold for $\bar\xi\ten\eta\in\bar\Hcal_s\ten\Hom_{\Rep(G)_{-n}}(\unit,U_s)$:
\begin{align*}
\|\lambda_n(\bar\xi\ten\eta)\|_{\varphi_\mu^{\ten n}}^2&\leq \|\xi\|^2\|\eta\|_{\mu^{\ten n}}^2; \\
\|\lambda_n(\bar\xi\ten\eta)^*\|_{\varphi_\mu^{\ten n}}^2&\leq d_s^2\|\xi\|^2\|\eta^*\|_{\mu^{\ten n}}^2.
\end{align*}
In particular, if $\bar\xi\ten\eta\in\bar\Hcal_s\ten\Hom_{\Rep(G)_{-1}}(\unit,U_s)$ and $n>m$, then
\begin{align*}
\big\|\hat\Delta^{n-1}(\lambda_1(\bar\xi\ten\eta))-1^{\ten n-m}\ten\hat\Delta^{m-1}& (\lambda_1(\bar\xi\ten\eta))\big\|_{\varphi_\mu^{\ten n}}^2 \\
&\leq \|\xi\|^2\big\|\hat\Delta^{n-1}(\eta)-(\iota^{\ten n-m}\ten\hat\Delta^{m-1})(\eta)\big\|_{\mu^{\ten n}}^2;\displaybreak[2]\\
\big\|(\hat\Delta^{n-1}(\lambda_1(\bar\xi\ten\eta))-1^{\ten n-m}\ten\hat\Delta^{m-1}& (\lambda_1(\bar\xi\ten\eta)))^*\big\|_{\varphi_\mu^{\ten n}}^2\\
&\leq d_s^2\,\|\xi\|^2 \big\|(\hat\Delta^{n-1}(\eta)-(\iota^{\ten n-m}\ten\hat\Delta^{m-1})(\eta))^*\big\|_{\mu^{\ten n}}^2.
\end{align*}
\end{Lem}
\begin{pf}
Let $\bar\xi\ten\eta\in\bar\Hcal_s\ten\Hom_{\Rep(G)_{-1}}(\unit,U_s)$, we compute using \eqref{thm_duality_G-Cstar_eq1}, \eqref{thm_duality_G-Cstar_eq2} and Corollary \ref{Cor_compatibility_lambda} that
\begin{align}
\|\lambda_n(\bar\xi\ten\eta)\|_{\varphi_\mu^{\ten n}}^2 &= (\varphi_\mu\ten\cdots\ten\varphi_\mu) \big(\lambda_n(\bar\xi\ten\eta)^*\lambda_n(\bar\xi\ten\eta)\big)\notag\\
&= (\varphi_\mu\ten\cdots\ten\varphi_\mu)\big(\lambda_n((\bar\xi\ten\eta)^*(\bar\xi\ten\eta))\big) \notag\displaybreak[2]\\
&= (\varphi_\mu\ten\cdots\ten\varphi_\mu) \big(\lambda_n(\pi((\bar\xi\ten\eta)^\bullet)(\bar\xi\ten\eta))\big) \notag\displaybreak[2]\\
&= (\varphi_\mu\ten\cdots\ten\varphi_\mu) \Big(\lambda_n\Big(\pi\big((\overline{\overline{\rho_s^{-1/2}\xi}}\ten(\eta^*\ten\iota)(\iota^{\ten n}\ten\bar R_s))\cdot(\bar\xi\ten\eta)\big)\Big)\Big)\notag\displaybreak[2]\\
&=\lambda_0\Big((\iota_{\bar\Hcal}\ten\tr_\mu\ten\cdots\ten\tr_\mu) \Big(\pi\big((\overline{\overline{\rho_s^{-1/2}\xi}\ten\xi})\ten(((\eta^*\ten\iota)(\iota^{\ten n}\ten\bar R_s))\ten\eta)\big)\Big)\Big)\notag\\
&=(\overline{\overline{\rho_s^{-1/2}\xi}\ten\xi})\big((\tr_\mu\ten\cdots\ten\tr_\mu) (((\eta^*\ten\iota)(\iota^{\ten n}\ten\bar R_s))\ten\eta)\big). \label{Lem_estimate_trace_on_Bn_eq1}
\end{align}
Note that since $\eta\in\Hom_{\Rep(G)_{-1}}(\unit,U_s)\subset\nat_b(\iota^{\ten n}\ten\unit, \iota^{\ten n}\ten U)$ the tensor product of natural transformations yields
\[
((\eta^*\ten\iota)(\iota^{\ten n}\ten\bar R_s))\ten\eta = (\eta^*\ten\iota_{\bar s}\ten\iota_s)(\iota^{\ten n}\ten\bar R_s\ten \iota_s)\eta.
\]
By means of Lemma \ref{Lem_categorical_path_space} write $\eta=\sum_{i=1}^k x_i\ten T_i$, where $x_i\in l^\infty(\hat G)^{\ten n}$ and $T_i\in B(\Hcal_\unit,\Hcal_{U_s})=B(\Com,\Hcal_s)$. Since $B(\Com,\Hcal_s)$ is finite dimensional, we may choose $k=\dim(U_s)$ and $T_i\colon\Com\ra\Hcal_s$, $c\mapsto c\zeta_i$, where $\{\zeta_i\}_{i=1}^{\dim(U_s)}$ forms an orthonormal basis in $\Hcal_s$ and $\langle\zeta_i,\xi\rangle =0$ if $i\neq 1$. Recall the solutions $(R_s,\bar R_s)$ of the conjugate equations for $U_s$ (see \eqref{Not_conjugate_equations_eq1}). With these choices we obtain
\begin{align}
\eqref{Lem_estimate_trace_on_Bn_eq1} &=\sum_{i,j=1}^{\dim(U_s)} (\overline{\overline{\rho_s^{-1/2}\xi}\ten\xi})\big((\varphi_\mu\ten\cdots \ten\varphi_\mu\ten\iota_{\bar U_s\ten U_s}) \notag\\
&\qquad\quad\big((x_i^*\ten T_i^*\ten\iota_{\bar s}\ten\iota_s)(\iota^{\ten n}\ten\bar R_s\ten\iota_s)(x_j\ten T_j)\big)\big) \notag\displaybreak[2]\\
&= \sum_{i,j=1}^{\dim(U_s)} (\varphi_\mu\ten\cdots\ten\varphi_\mu)(x_i^*x_j)\, (\overline{\overline{\rho_s^{-1/2}\xi}}\ten\bar\xi)\big((T_i^*\ten\iota_{\bar s}\ten\iota_s)(\bar R_s\ten\iota_s)T_j\big)\notag\displaybreak[2]\\
&= \sum_{i,j=1}^{\dim(U_s)} (\varphi_\mu\ten\cdots\ten\varphi_\mu)(x_i^*x_j)\sum_{k=1}^{\dim(U_s)} T_i^*(\rho_s^{\half}\xi_k^s)\,\big((\overline{\overline{\rho_s^{-1/2}\xi}})(\bar\xi_k^s)\big)\, \bar\xi(T_j) \notag\displaybreak[2]\\
&= \sum_{i,j=1}^{\dim(U_s)} (\varphi_\mu\ten\cdots\ten\varphi_\mu)(x_i^*x_j)\sum_{k=1}^{\dim(U_s)} T_i^*(\xi_k^s)\,\bar\xi_k^s(\xi)\,\bar\xi(T_j) \notag\\
&= \sum_{i,j=1}^{\dim(U_s)} (\varphi_\mu\ten\cdots\ten\varphi_\mu)(x_i^*x_j)\,T_i^*(\xi)\,\bar\xi(T_j). \label{Lem_estimate_trace_on_Bn_eq2}
\end{align}
By the choice of the $T_i$'s in the decomposition of $\eta$, this equals
\begin{align*}
\eqref{Lem_estimate_trace_on_Bn_eq2}&= \sum_{i,j=1}^{\dim(U_s)} (\varphi_\mu\ten\cdots\ten\varphi_\mu)(x_i^*x_j)\,\langle\zeta_i,\xi\rangle\, \langle\xi,\zeta_j\rangle\\
&= (\varphi_\mu\ten\cdots\ten\varphi_\mu)(x_1^*x_1)\,\|\xi\|^2 \displaybreak[2]\\
&\leq \|\xi\|^2 \sum_{i=1}^{\dim(U_s)} (\varphi_\mu\ten\cdots\ten\varphi_\mu)(x_i^*x_i)\, \varphi_\unit(T_i^*T_i) \notag\displaybreak[2]\\
&= \|\xi\|^2 (\tr_\mu\ten\cdots\ten\tr_\mu\ten\tr_\unit)(\eta^*\eta)\\
&=\|\xi\|^2\|\eta\|_{\mu^{\ten n}}^2.
\end{align*}

The second estimate is similar but slightly trickier. Along the same lines one can show that
\begin{equation}\label{Lem_estimate_trace_on_Bn_eq3}
\|\lambda_n(\bar\xi\ten\eta)^*\|_{\varphi_\mu^{\ten n}}^2
=(\overline{\xi\ten\overline{\rho_s^{-1/2}\xi}})\big((\tr_\mu\ten\cdots\ten\tr_\mu) ((\eta\ten\iota_{\bar s})(\eta^*\ten\iota_{\bar s})(\iota^{\ten n}\ten \bar R_s))\big).
\end{equation}
As before decompose $\eta$ by means of Lemma \ref{Lem_categorical_path_space} as $\eta=\sum_{i=1}^{\dim(U_s)}x_i\ten T_i$, but this time let $T_i\colon\Com\ra\Hcal_s$, $c\mapsto c\xi_i^s$, where $\{\xi_i^s\}_i$ are the eigenvectors of $\rho_s$. Assume for the moment that $\xi=\xi_l^s$ for some $l\in\{1,\ldots,\dim(U_s)\}$. Then in this case
\begin{align}
\eqref{Lem_estimate_trace_on_Bn_eq3}&= \sum_{i,j=1}^{\dim(U_s)} (\overline{\xi_l^s\ten\overline{\rho_s^{-1/2}\xi_l^s}}) \big((\varphi_\mu\ten\cdots\ten\varphi_\mu\ten\iota_{U_s\ten\bar U_s})\big((x_i\ten T_i\ten\iota_{\bar s}) (x_j^*\ten T_j^*\ten\iota_{\bar s})(\iota^{\ten n}\ten\bar R_s)\big)\big) \notag\\
&=\sum_{i,j=1}^{\dim(U_s)}(\varphi_\mu\ten\cdots\ten\varphi_\mu)(x_ix_j^*) \sum_{k=1}^{\dim(U_s)}(\bar\xi_l^s\ten\overline{\overline{\rho_s^{-1/2}\xi_l^s}}) \big(T_iT_j^*(\rho_s^\half\xi_k^s)\ten\bar\xi_k^s\big) \notag\displaybreak[2]\\
&=\sum_{i,j=1}^{\dim(U_s)}(\varphi_\mu\ten\cdots\ten\varphi_\mu)(x_ix_j^*) \sum_{k=1}^{\dim(U_s)} \bar\xi_l^s(T_i)\,T_j^*(\rho_s^\half\xi_k^s)\,\big((\overline{\overline{\rho_s^{-1/2}\xi_l^s}}) (\bar\xi_k^s)\big) \notag\displaybreak[2]\\
&=\sum_{i,j=1}^{\dim(U_s)}(\varphi_\mu\ten\cdots\ten\varphi_\mu)(x_ix_j^*) \sum_{k=1}^{\dim(U_s)} \langle\xi_l^s,\xi_i^s\rangle\,\langle\xi_j^s,\rho_s^\half\xi_k^s\rangle\,\langle\xi_k^s, \rho_s^{-\half}\xi_l^s\rangle \notag\\
&= (\varphi_\mu\ten\cdots\ten\varphi_\mu)(x_lx_l^*).\label{Lem_estimate_trace_on_Bn_eq4}
\end{align}
On the other hand
\begin{align}
\|\xi\|^2\|\eta^*\|^2_{\mu^{\ten n}} &= (\tr_\mu\ten\cdots\ten\tr_\mu\ten\tr_{U_s})(\eta\eta^*) \notag\\
&= \sum_{i,j=1}^{\dim(U_s)}(\varphi_\mu\ten\cdots\ten\varphi_\mu\ten\varphi_s)(x_ix_j^*\ten T_iT_j^*)\notag \displaybreak[2]\\
&= \sum_{i,j=1}^{\dim(U_s)}(\varphi_\mu\ten\cdots\ten\varphi_\mu)(x_ix_j^*) \sum_{k=1}^{\dim(U_s)}\langle\xi_k^s,\xi_i^s\langle\xi_j^s,\rho^{-1}_s\xi_k^s\rangle\rangle \notag\\
&= \sum_{i=1}^{\dim(U_s)}(\varphi_\mu\ten\cdots\ten\varphi_\mu)(x_ix_i^*)(\rho_s^{-1})_{ii}. \label{Lem_estimate_trace_on_Bn_eq5}
\end{align}
Note that $d_s=\sum_{j=1}^{\dim(U_s)}(\rho_s)_{jj}\geq (\rho_s)_{ii}$ for any $i$. Therefore comparing  \eqref{Lem_estimate_trace_on_Bn_eq4} and \eqref{Lem_estimate_trace_on_Bn_eq5} gives
\begin{equation}\label{Lem_estimate_trace_on_Bn_eq6}
\|\lambda_n(\bar\xi\ten\eta)^*\|_{\varphi_\mu^{\ten n}}^2\leq d_s\|\xi\|^2\|\eta^*\|^2_{\mu^{\ten n}}.
\end{equation}
Now we deal with general vectors $\xi\in\Hcal_s$ as follows. Write $\xi=\sum_i c_i\xi_i^s$, then $\|\xi\|^2=\sum_i |c_i|^2$. By \eqref{Lem_estimate_trace_on_Bn_eq6}, Cauchy--Schwarz and Jensen's inequality we get
\begin{align*}
\|\lambda_n(\bar\xi\ten\eta)^*\|_{\varphi_\mu^{\ten n}}^2 &\leq \Big(\sum_{i=1}^{\dim(U_s)}|c_i|\|\lambda_n(\xi_i^s\ten\eta)^*\|_{\varphi_\mu^{\ten n}}\Big)^2\\
&\leq d_s \Big(\sum_{i=1}^{\dim(U_s)}|c_i|\Big)^2 \|\eta^*\|_{\mu^{\ten n}}^2\displaybreak[2]\\
&\leq d_s \dim(U_s) \sum_{i=1}^{\dim(U_s)}|c_i|^2\, \|\eta^*\|_{\mu^{\ten n}}^2\\
&\leq d_s^2\|\xi\|^2\|\eta^*\|^2_{\mu^{\ten n}},
\end{align*}
which completes the second inequality.

The third and fourth inequalities now have become easy. Indeed, by Corollary \ref{Cor_compatibility_lambda} and the first estimate we get for $\bar\xi\ten\eta\in\bar\Hcal_s\ten\Hom_{\Rep(G)_{-1}}(\unit,U_s)$ and $n>m$ that
\begin{align*}
\big\|\hat\Delta^{n-1}(\lambda_1(\bar\xi\ten\eta))&-1^{\ten n-m}\ten\hat\Delta^{m-1} (\lambda_1(\bar\xi\ten\eta))\big\|_{\varphi_\mu^{\ten n}}^2\\
&=\|\lambda_n(\bar\xi\ten\hat\Delta^{n-1}(\eta))-1^{\ten n-m}\ten\lambda_m(\bar\xi\ten\hat\Delta^{m-1}(\eta))\|_{\varphi_\mu^{\ten n}}^2 \displaybreak[2]\\
&=\|\lambda_n(\bar\xi\ten\hat\Delta^{n-1}(\eta))-\lambda_n(\bar\xi\ten\iota^{\ten n-m}\ten\hat\Delta^{m-1}(\eta))\|_{\varphi_\mu^{\ten n}}^2 \displaybreak[2]\\
&=\|\lambda_n(\bar\xi\ten(\hat\Delta^{n-1}(\eta)-\iota^{\ten n-m}\ten\hat\Delta^{m-1}(\eta)))\|_{\varphi_\mu^{\ten n}}^2\\
&\leq \|\xi\|^2\big\|\hat\Delta^{n-1}(\eta)-\iota^{\ten n-m}\ten\hat\Delta^{m-1}(\eta)\big\|_{\mu^{\ten n}}^2.
\end{align*}

The fourth estimate is similar. We leave the details to the reader.
\end{pf}

Because of the Lemma \ref{Lem_regular_bcYD-algebra} and Theorem \ref{Thm_duality_YD} there exists a \Cstar tensor category corresponding to the braided commutative Yetter--Drinfeld $G$-\Cstar algebra $\tilde R(\hat G,\mu)$. Denote this \Cstar tensor category by $\Dcal$. So concretely, $\Dcal=\Ccal_{\tilde R(\hat G,\mu)}$ is the completion of the category with objects $\Ob(\Rep(G))$ and morphism sets
\[
\Hom_{\Dcal}(U,V):=\{T\in\tilde R(\hat G,\mu)\ten B(\Hcal_U,\Hcal_V)\,:\, V_{13}^*(\alpha_l\ten\iota_V)(T)U_{13}=1\ten T\},
\]
for $U,V\in\Ob(\Rep(G))$. Recall the \Cstar tensor category $R(\Rep(G),\mu)$ of $\mu$-regular natural transformations defined in Definition \ref{Def_cat_regular}.

\begin{Thm}\label{Thm_equivalence_regular_elements}
The category $\Dcal$ is unitarily monoidally equivalent to $R(\Rep(G),\mu)$ via a monoidal equivalence  preserving the canonical functors $\Rep(G)\ra\Dcal$ and $\Rep(G)\ra R(\Rep(G),\mu)$.
\end{Thm}
\begin{pf}
For an element $T\in l^\infty(\hat G)\ten B(\Hcal_U,\Hcal_V)$ recall the condition
\begin{equation}\label{Thm_equivalence_regular_elements_eq1}
V_{13}^*(\alpha_l\ten\iota_V)(T)U_{13}=1\ten T.
\end{equation}
The categories $\Rep(G)_{-1}$ and $\Ccal_{C_{-1}(\hat G)}$ are unitarily monoidally equivalent as $\Rep(G)$-module categories (see Corollary \ref{Cor_equiv_path_spaces}). To prove the theorem it therefore suffices to show that for any pair of objects $U,V\in\Ob(\Rep(G))$, the isomorphism of Lemma \ref{Lem_categorical_path_space} restricts to an isomorphism
\[
\Hom_{R(\Rep(G),\mu)}(U,V)\xleftarrow{\cong} \{T\in\tilde R(\hat G,\mu)\ten B(\Hcal_U,\Hcal_V)\,:\, V_{13}^*(\alpha_l\ten\iota_V)(T)U_{13}=1\ten T\}.
\]
In order to achieve this it is sufficient to show the following two statements:
\begin{enumerate}[label=(\alph*)]
\item $\eta(T)\in\Hom_{R(\Rep(G),\mu)}(U,V)$ for any $T\in\tilde R(\hat G,\mu)\ten B(\Hcal_U,\Hcal_V)$ satisfying condition \eqref{Thm_equivalence_regular_elements_eq1};
\item for any $\eta\in\Hom_{R(\Rep(G),\mu)}(U,V)$ there exists a homomorphism $T\in\tilde R(\hat G,\mu)\ten B(\Hcal_U,\Hcal_V)$ satisfying condition \eqref{Thm_equivalence_regular_elements_eq1} such that $\eta=\eta(T)$.
\end{enumerate}
(a) Assume $T\in\tilde R(\hat G,\mu)\ten B(\Hcal_U,\Hcal_V)$ satisfies \eqref{Thm_equivalence_regular_elements_eq1}. By Lemma \ref{Lem_categorical_path_space}, $\eta(T)\in\Hom_{\Rep(G)_{-1}}(U,V)$. We must show that $\eta(T)$ is $\mu$-regular. Since $B(\Hcal_U,\Hcal_V)$ is finite dimensional we can write $T=\sum_{i=1}^k x_i\ten T_i$ with $x_i\in \tilde R(\hat G,\mu)$ and $T_i\in B(\Hcal_U,\Hcal_V)$ such that $\varphi_U(T_i^*T_j)=\varphi_V(T_iT_j^*)=0$ whenever $i\neq j$ (for instance, take $T_i$ of the form $m_{\xi,\xi'}$ where $\xi$ and $\xi'$ are in a basis of $\Hcal_U$ and respectively $\Hcal_V$ in which $\rho$ acts diagonally). Let $n>m$. By Lemma \ref{Lem_properties_isom_path_space} we obtain
\begin{align}
\|\hat\Delta^{n-1}&(\eta(T))-\iota^{\ten n-m}\ten\hat\Delta^{m-1}(\eta(T))\|^2_{\mu^{\ten n}} \label{Thm_equivalence_regular_elements_eq2}\\
&= \|\eta((\hat\Delta^{n-1}\ten\iota_V)(T)-(1^{\ten n-m}\ten\hat\Delta^{m-1}\ten\iota_V)(T))\|^2_{\mu^{\ten n}} \notag \displaybreak[2]\\
&= (\varphi_\mu\ten\cdots\ten\varphi_\mu\ten\varphi_U)\big(((\hat\Delta^{n-1}\ten\iota_V)(T)- (1^{\ten n-m}\ten\hat\Delta^{m-1}\ten\iota_V)(T))^* \notag\\
&\qquad\times((\hat\Delta^{n-1}\ten\iota_V)(T)-(1^{\ten n-m}\ten\hat\Delta^{m-1}\ten\iota_V)(T))\big) \notag \displaybreak[2]\\
&= \sum_{i,j=1}^k (\varphi_\mu\ten\cdots\ten\varphi_\mu) \big((\hat\Delta^{n-1}(x_i)-1^{\ten n-m}\ten\hat\Delta^{m-1}(x_i))^* \notag\\
&\qquad\times(\hat\Delta^{n-1}(x_j)-1^{\ten n-m}\ten\hat\Delta^{m-1}(x_j))\big)\varphi_U(T_i^*T_j) \notag \displaybreak[1]\\
&= \sum_{i=1}^k \varphi_U(T_i^*T_i)\, \|\hat\Delta^{n-1}(x_i)-1^{\ten n-m}\ten\hat\Delta^{m-1}(x_i)\|^2_{\varphi_\mu^{\ten n}}. \notag
\end{align}
As each $x_i$ is $\mu$-regular, it follows that for $i=1,\ldots, k$
\[
\|\hat\Delta^{n-1}(x_i)-1^{\ten n-m}\ten\hat\Delta^{m-1}(x_i)\|^2_{\varphi_\mu^{\ten n}}\ra 0 \qquad \textrm{as } n,m\ra\infty
\]
and thus \eqref{Thm_equivalence_regular_elements_eq2} converges to zero as $m,n\ra\infty$. A similar calculation shows
\[
\|\hat\Delta^{n-1}(\eta(T)^*)-\iota^{\ten n-m}\ten\hat\Delta^{m-1}(\eta(T)^*)\|^2_{\mu^{\ten n}}\ra 0 \qquad \textrm{as } n,m\ra \infty.
\]
We conclude that $\eta(T)\in\Hom_{R(\Rep(G),\mu)}(U,V)$.\\
(b) First assume that $\eta\in\Hom_{R(\Rep(G),\mu)}(\unit, U_s)$. Pick the unique $T\in C_{-1}(\hat G)\ten B(\Com,\Hcal_s)$ satisfying \eqref{Thm_equivalence_regular_elements_eq1} such that $\eta(T)=\eta$. Given $\xi\in\Hcal_s$, observe that $(\iota\ten\bar\xi)(T)=\lambda_1(\bar\xi\ten\eta)$. From Lemma \ref{Lem_estimate_trace_on_Bn} it immediately follows that $(\iota\ten\bar\xi)(T)\in \tilde R(\hat G,\mu)$ and thus $T\in\tilde R(\hat G,\mu)\ten B(\Com,\Hcal_s)$. For general objects $U,V\in\Ob(\Rep(G))$ Frobenius reciprocity gives
\[
\Hom_{R(\Rep(G),\mu)}(U,V)\cong \Hom_{R(\Rep(G),\mu)}(\unit,\bar U\ten V) \cong \ds_{s\in\Irr(G)} m_{\bar U,V}^s \Hom_{R(\Rep(G),\mu)}(\unit, U_s),
\]
as desired.
\end{pf}

\begin{Thm}\label{Thm_equivalence_boundary_convergence}
Given a reduced compact quantum group $G$ and a generating and transient probability measure $\mu$ on $\Irr(G)=\Irr(\Rep(G))$, the random walk defined by $\mu$ on the discrete dual $l^\infty(\hat G)$ converges to the boundary if and only if the random walk defined by $\mu$ on the \Cstar tensor category $\Rep(G)$ converges to the boundary.
\end{Thm}

Most of the hard work of the proof of this theorem is already done. We will refer to requirements (i) and (ii) of Definition \ref{Def_cat_convergence_boundary} as (i)$^{\textrm{cat}}$ respectively (ii)$^{\textrm{cat}}$, to distinguish them from the corresponding properties in Conjecture \ref{conjecture_convergence}.

\medskip
\begin{pf}[ of Theorem \ref{Thm_equivalence_boundary_convergence}]
Note that $\tilde M(\hat G,\mu)\subset C_{-1}(\hat G)$. Therefore Theorem \ref{Thm_equivalence_regular_elements} implies that $\tilde M(\hat G,\mu)\subset \tilde R(\hat G,\mu)$ if and only if $\tilde M(\Rep(G),\mu)$ is a \Cstar tensor subcategory of $R(\Rep(G),\mu)$. Thus (i) and (i)$^{\textrm{cat}}$ are equivalent.

Equivalence of (ii) and (ii)$^{\textrm{cat}}$ is proved in the subsequent two lemmas.
\end{pf}

\begin{Lem}\label{Lem_ii_quantum_to_cat}
Suppose that $(G,\mu)$ satisfies {\rm (ii)}, then $(\Rep(G),\mu)$ satisfies {\rm (ii)}$^{\textrm{cat}}$.
\end{Lem}
\begin{pf}
Let $\nu\in\nat_b(\iota\ten X,\iota\ten Y)$ be $\mu$-harmonic and $\eta\in\nat_{00}(\iota\ten U,\iota\ten V)$. By Lemma \ref{Lem_lambda_isomorphism_subalgebras}, $\eta$ is of the form $\eta=\eta(T)$ for a unique $T$ satisfying \eqref{Thm_equivalence_regular_elements_eq1}. Write $T$ as $
T =\sum_{i=1}^k x_i\ten T_i\in c_{00}(\hat G)\ten B(\Hcal_U,\Hcal_V)$.
Similarly, by \cite[Thm.\ 4.1]{NeshveyevYamashita14c} also $\nu$ is of the form $\nu=\eta(S)$ for some $S$, where $S=\sum_{i=1}^m y_j\ten S_j\in H^\infty(\hat G,\mu)\ten B(\Hcal_X,\Hcal_Y)$. Note that these sums are finite, because $B(\Hcal_U,\Hcal_V)$ and $B(\Hcal_X,\Hcal_Y)$ are finite dimensional.
The Poisson boundary is invariant under the right action defined by the comultiplication (\cite[Lem.\ 2.2]{Izumi02}). Thus if we write $\hat\Delta(y_j)=y_j^{(1)}\ten y_j^{(2)}$, then $y_j^{(1)}\in H^\infty(\hat G,\mu)$. Use the description of the tensor product of natural isomorphisms given by \eqref{Lem_categorical_path_space_eq4} to obtain
\begin{align*}
\sum_{s\in\Irr(G)} d_s^2 \tr_s(\eta\ten\nu) &= \sum_{s\in\Irr(G)}\sum_{i,j} d_s^2(\varphi_s\ten\iota_V\ten\iota_Y)\big(x_iy_j^{(1)}\ten T_i\pi_U(y_j^{(2)})\ten S_j\big)\\
&= \sum_{i,j} \Big(\sum_{s\in\Irr(G)}d_s^2 \varphi_s(x_iy_j^{(1)})\Big) \big(T_i\pi_U(y_j^{(2)})\ten S_j\big) \displaybreak[2]\\
&= \sum_{i,j} \varphi_\mu^\infty\big(K_{\bar\mu}(x_i)y_j^{(1)}\big) \big(T_i\pi_U(y_j^{(2)})\ten S_j\big) \displaybreak[2]\\
&= (\varphi_\mu^\infty\ten\iota_V\ten\iota_Y)\Big(\sum_{i,j} K_{\bar\mu}(x_i)y_j^{(1)}\ten T_i\pi_U(y_j^{(2)})\ten S_j\Big)\\
&= \tr_\mu^\infty\big(K_{\bar\mu}(\eta)\ten\nu\big).
\end{align*}
In the last step we used Lemma \ref{Lem_properties_isom_path_space}.
\end{pf}

\begin{Lem}\label{Lem_ii_cat_to_quantum}
If $(\Rep(G),\mu)$ satisfies {\rm (ii)}$^\textrm{cat}$, then $(G,\mu)$ satisfies {\rm (ii)}.
\end{Lem}
\begin{pf}
Let $\bar\xi\ten\eta\in\bar\Hcal_s\ten\Hom_{\Rep(G)_{-1}}(\unit,U_s)$ be such that $\lambda_1(\bar\xi\ten\eta)\in c_{00}(\hat G)$ and let $\bar\zeta\ten\nu\in\bar\Hcal_t\ten\Hom_{\Rep(G)_{-1}}(\unit,U_t)$ be such that $\lambda_1(\bar\zeta\ten\nu)\in H^\infty(\hat G,\mu)$. Recall the tensor product of natural transformations \eqref{Def_tensor_product_nat_transf_eq0}. Using the notation and results of Corollary \ref{Cor_compatibility_lambda} we get
\begin{align}
\hat\psi(\lambda_1(\bar\xi\ten\eta)\lambda_1(\bar\zeta\ten\nu))&= \hat\psi\big(\lambda_1((\bar\xi\ten\eta)(\bar\zeta\ten\nu))\big) \notag\\
&= \sum_r d_r^2\,\varphi_r\big(\lambda_1((\overline{\xi\ten\zeta})\ten((\eta\ten\iota)(\iota\ten\nu)))\big) \notag \displaybreak[2]\\
&= \sum_r d_r^2\,\lambda_0 \circ (\iota_{\overline{\Hcal_s\ten\Hcal_t}}\ten\tr_r) \big((\overline{\xi\ten\zeta})\ten(\eta\ten\nu)\big) \notag \displaybreak[2]\\
&= \sum_r d_r^2\lambda_0((\overline{\xi\ten\zeta})\ten\tr_r(\eta\ten\nu)) \notag\\
&= \lambda_0\Big((\overline{\xi\ten\zeta})\ten \Big(\sum_r d_r^2\tr_r(\eta\ten\nu)\Big)\Big).\label{Lem_ii_cat_to_quantum_eq2}
\end{align}
By assumption $(\Rep(G),\mu)$ satisfies (ii)$^{\textrm{cat}}$, so Corollary \ref{Cor_compatibility_lambda} gives
\begin{align*}
\eqref{Lem_ii_cat_to_quantum_eq2}&= \lambda_0\big((\overline{\xi\ten\zeta})\ten \tr_\mu^\infty(K_{\bar\mu}(\eta)\ten\nu)\big)\\
&= \lambda_0\circ(\iota_{\overline{\Hcal_s\ten\Hcal_t}}\ten\tr_\mu^\infty) \big((\overline{\xi\ten\zeta})\ten((K_{\bar\mu}(\eta)\ten\iota)(\iota\ten\nu))\big)\\
&= \varphi_\mu^\infty\big(\lambda_1((\overline{\xi\ten\zeta}) \ten((K_{\bar\mu}(\eta)\ten\iota)(\iota\ten\nu)))\big)\\
&= \varphi_\mu^\infty\big(K_{\bar\mu}(\lambda_1(\bar\xi\ten\eta))\lambda_1(\bar\zeta\ten\nu)\big).
\end{align*}
By linearity it follows that for any $x\in\lambda_1\big(\ds_s\bar\Hcal_s\ten\nat_{00}(\unit,U_s)\big)$ and any $h=\lambda_1(T)$ where $T\in\ds_s\bar\Hcal_s\ten\nat_b(\iota,\iota\ten U_s)$ with $(\iota\ten P_\mu)(T)=T$, it holds that
\begin{equation}\label{Lem_ii_cat_to_quantum_eq1}
\hat\psi(xh) = \varphi_\mu^\infty(K_{\bar\mu}(x)h).
\end{equation}
By density and strong$^*$-continuity we obtain that \eqref{Lem_ii_cat_to_quantum_eq1} holds for any $x\in c_{00}(\hat G)$ and any $\mu$-harmonic element $h\in H^\infty(\hat G,\mu)$.
\end{pf}

\bibliographystyle{alpha}
\bibliography{C:/Documenten/PhD/Papers/References/references}

\begin{thebibliography}{DRVV10}

\bibitem[Ban96]{Banica96}
T.~Banica.
\newblock Th\'eorie des repr\'esentations du groupe quantique compact libre
  {${\rm O}(n)$}.
\newblock {\em C. R. Acad. Sci. Paris S\'er. I Math.}, 322(3):241--244, 1996.

\bibitem[BDRV06]{BichonDeRijdtVaes06}
J.~Bichon, A.~De~Rijdt, and S.~Vaes.
\newblock Ergodic coactions with large multiplicity and monoidal equivalence of
  quantum groups.
\newblock {\em Comm. Math. Phys.}, 262(3):703--728, 2006.

\bibitem[Bia91]{Biane91}
P.~Biane.
\newblock Quantum random walk on the dual of {${\rm SU}(n)$}.
\newblock {\em Probab. Theory Related Fields}, 89(1):117--129, 1991.

\bibitem[Bia92a]{Biane92a}
P.~Biane.
\newblock \'{E}quation de {C}hoquet-{D}eny sur le dual d'un groupe compact.
\newblock {\em Probab. Theory Related Fields}, 94(1):39--51, 1992.

\bibitem[Bia92b]{Biane92b}
P.~Biane.
\newblock Minuscule weights and random walks on lattices.
\newblock In {\em Quantum probability \& related topics}, QP-PQ, VII, pages
  51--65. World Sci. Publ., River Edge, NJ, 1992.

\bibitem[Bia94]{Biane94}
P.~Biane.
\newblock Th\'eor\`eme de {N}ey-{S}pitzer sur le dual de {${\rm SU}(2)$}.
\newblock {\em Trans. Amer. Math. Soc.}, 345(1):179--194, 1994.

\bibitem[BR97]{BratteliRobinson97}
O.~Bratteli and D.~W. Robinson.
\newblock {\em Operator algebras and quantum statistical mechanics. 2}.
\newblock Texts and Monographs in Physics. Springer-Verlag, Berlin, second
  edition, 1997.
\newblock Equilibrium states. Models in quantum statistical mechanics.

\bibitem[CE77]{ChoiEffros77}
M.~D. Choi and E.~G. Effros.
\newblock Injectivity and operator spaces.
\newblock {\em J. Functional Analysis}, 24(2):156--209, 1977.

\bibitem[Col04]{Collins04}
B.~Collins.
\newblock Martin boundary theory of some quantum random walks.
\newblock {\em Ann. Inst. H. Poincar\'e Probab. Statist.}, 40(3):367--384,
  2004.

\bibitem[Con75]{Connes75}
A.~Connes.
\newblock On hyperfinite factors of type {${\rm III}_{0}$} and {K}rieger's
  factors.
\newblock {\em J. Functional Analysis}, 18:318--327, 1975.

\bibitem[DCY13]{DeCommerYamashita13}
K.~De~Commer and M.~Yamashita.
\newblock Tannaka-{K}re\u\i n duality for compact quantum homogeneous spaces.
  {I}. {G}eneral theory.
\newblock {\em Theory Appl. Categ.}, 28:No. 31, 1099--1138, 2013.

\bibitem[Doo59]{Doob59}
J.~L. Doob.
\newblock Discrete potential theory and boundaries.
\newblock {\em J. Math. Mech.}, 8:433--458; erratum 993, 1959.

\bibitem[DRVV10]{DeRijdtVanderVennet10}
A.~De~Rijdt and N.~Vander~Vennet.
\newblock Actions of monoidally equivalent compact quantum groups and
  applications to probabilistic boundaries.
\newblock {\em Ann. Inst. Fourier (Grenoble)}, 60(1):169--216, 2010.

\bibitem[Hun60]{Hunt60}
G.~A. Hunt.
\newblock Markoff chains and {M}artin boundaries.
\newblock {\em Illinois J. Math.}, 4:313--340, 1960.

\bibitem[Izu02]{Izumi02}
M.~Izumi.
\newblock Non-commutative {P}oisson boundaries and compact quantum group
  actions.
\newblock {\em Adv. Math.}, 169(1):1--57, 2002.

\bibitem[Izu12]{Izumi12}
M.~Izumi.
\newblock {$E_0$}-semigroups: around and beyond {A}rveson's work.
\newblock {\em J. Operator Theory}, 68(2):335--363, 2012.

\bibitem[KS97]{KlimykSchmudgen97}
A.~Klimyk and K.~Schm{\"u}dgen.
\newblock {\em Quantum groups and their representations}.
\newblock Texts and Monographs in Physics. Springer-Verlag, Berlin, 1997.

\bibitem[KT99]{KustermansTuset99}
J.~Kustermans and L.~Tuset.
\newblock A survey of {$C^*$}-algebraic quantum groups. {I}.
\newblock {\em Irish Math. Soc. Bull.}, 43:8--63, 1999.

\bibitem[Mar41]{Martin41}
R.~S. Martin.
\newblock Minimal positive harmonic functions.
\newblock {\em Trans. Amer. Math. Soc.}, 49:137--172, 1941.

\bibitem[{Nes}13]{Neshveyev13}
S.~{Neshveyev}.
\newblock {Duality theory for nonergodic actions}.
\newblock {\em ArXiv e-prints}, March 2013.

\bibitem[NT03]{NeshveyevTuset03}
S.~Neshveyev and L.~Tuset.
\newblock Quantum random walks and their boundaries.
\newblock In {\em Analysis of (Quantum) Group Actions on Operator Algebras},
  volume 1332, pages 57--70. RIMS Kokyuroku, Kyoto, 2003.

\bibitem[NT04]{NeshveyevTuset04}
S.~Neshveyev and L.~Tuset.
\newblock The {M}artin boundary of a discrete quantum group.
\newblock {\em J. Reine Angew. Math.}, 568:23--70, 2004.

\bibitem[NT13]{NeshveyevTuset13}
S.~Neshveyev and L.~Tuset.
\newblock {\em Compact quantum groups and their representation categories},
  volume~20 of {\em Cours Sp\'ecialis\'es}.
\newblock Soci\'et\'e Math\'ematique de France, Paris, 2013.

\bibitem[NY14a]{NeshveyevYamashita14c}
S.~Neshveyev and M.~Yamashita.
\newblock Categorical duality for {Y}etter-{D}rinfeld algebras.
\newblock {\em Doc. Math.}, 19:1105--1139, 2014.

\bibitem[NY14b]{NeshveyevYamashita14a}
S.~{Neshveyev} and M.~{Yamashita}.
\newblock {Poisson boundaries of monoidal categories}.
\newblock {\em ArXiv e-prints}, May 2014.

\bibitem[Pod87]{Podles87}
P.~Podle{\'s}.
\newblock Quantum spheres.
\newblock {\em Lett. Math. Phys.}, 14(3):193--202, 1987.

\bibitem[Rev84]{Revuz84}
D.~Revuz.
\newblock {\em Markov chains}, volume~11 of {\em North-Holland Mathematical
  Library}.
\newblock North-Holland Publishing Co., Amsterdam, second edition, 1984.

\bibitem[Tak03]{Takesaki03III}
M.~Takesaki.
\newblock {\em Theory of operator algebras. {III}}, volume 127 of {\em
  Encyclopaedia of Mathematical Sciences}.
\newblock Springer-Verlag, Berlin, 2003.
\newblock Operator Algebras and Non-commutative Geometry, 8.

\bibitem[VD96]{VanDaele96}
A.~Van~Daele.
\newblock Discrete quantum groups.
\newblock {\em J. Algebra}, 180(2):431--444, 1996.

\bibitem[VDW96]{VanDaeleWang96}
A.~Van~Daele and S.~Wang.
\newblock Universal quantum groups.
\newblock {\em Internat. J. Math.}, 7(2):255--263, 1996.

\bibitem[Voi11]{Voigt11}
C.~Voigt.
\newblock The {B}aum-{C}onnes conjecture for free orthogonal quantum groups.
\newblock {\em Adv. Math.}, 227(5):1873--1913, 2011.

\bibitem[VV07]{VaesVergnioux07}
S.~Vaes and R.~Vergnioux.
\newblock The boundary of universal discrete quantum groups, exactness, and
  factoriality.
\newblock {\em Duke Math. J.}, 140(1):35--84, 2007.

\bibitem[VVV08]{VaesVanderVennet08}
S.~Vaes and N.~Vander~Vennet.
\newblock Identification of the {P}oisson and {M}artin boundaries of orthogonal
  discrete quantum groups.
\newblock {\em J. Inst. Math. Jussieu}, 7(2):391--412, 2008.

\bibitem[Woe00]{Woess00}
W.~Woess.
\newblock {\em Random walks on infinite graphs and groups}, volume 138 of {\em
  Cambridge Tracts in Mathematics}.
\newblock Cambridge University Press, Cambridge, 2000.

\end{thebibliography}
\end{document}